\documentclass[draft]{amsart}
\usepackage{amssymb,latexsym}
\usepackage[mathscr]{eucal}
\usepackage{verbatim}

\newtheorem*{theorem*}{Theorem}
\newtheorem*{corollary*}{Corollary}
\newtheorem{theorem}{Theorem}[section]
\newtheorem{corollary}[theorem]{Corollary}
\newtheorem{lemma}[theorem]{Lemma}
\newtheorem{proposition}[theorem]{Proposition}
\newtheorem{remark}[theorem]{Remark}
\newtheorem{definition}[theorem]{Definition}
\newtheorem{example}[theorem]{Example}

\DeclareMathOperator{\card}{card}
\DeclareMathOperator{\cl}{cl}
\DeclareMathOperator{\tr}{Tr}

\DeclareMathOperator{\slim}{s- lim}
\DeclareMathOperator{\diag}{diag}
\DeclareMathOperator{\co*}{c_o^*}

\def\sideremark#1{\ifvmode\leavevmode\fi\vadjust{\vbox to0pt{\vss
\hbox to 0pt{\hskip\hsize\hskip1em
\vbox{\hsize2cm\tiny\raggedright\pretolerance10000
\noindent#1\hfill}\hss}\vbox to8pt{\vfil}\vss}}}

\newcommand{\1}{\mspace{2mu}}
\newcommand{\be}{\begin{equation}\label}
\newcommand{\ee}{\end{equation}}
\newcommand{\bq}{\begin{equation*}}
\newcommand{\eq}{\end{equation*}}
\newcommand{\ba}{\begin{align*}}
\newcommand{\ea}{\end{align*}}
\newcommand{\bp}{\begin{proof}}
\newcommand{\ep}{\end{proof}}
\newcommand{\bL}{\begin{lemma}\label}
\newcommand{\eL}{\end{lemma}}
\newcommand{\bP}{\begin{proposition}\label}
\newcommand{\eP}{\end{proposition}}
\newcommand{\bC}{\begin{corollary}\label}
\newcommand{\eC}{\end{corollary}}
\newcommand{\bT}{\begin{theorem}\label}
\newcommand{\eT}{\end{theorem}}
\newcommand{\bR}{\begin{remark}\label}
\newcommand{\eR}{\end{remark}}
\newcommand{\bD}{\begin{definition}\label}
\newcommand{\eD}{\end{definition}}

\newcommand{\bE}{\begin{example}\label}
\newcommand{\eE}{\end{example}}

\begin{document}

\title[Infinite Dimensional Schur-Horn Theorem]%
{An infinite dimensional Schur-Horn theorem and \\majorization theory 
with applications to operator ideals}
\author{Victor Kaftal}
\address{University of Cincinnati\\
          Department of Mathematics\\
          Cincinnati, OH, 45221-0025\\
          USA}
\email{victor.kaftal@math.uc.edu}
\author{Gary Weiss}
\email{gary.weiss@math.uc.edu}

\date{\today}

\begin{abstract}
The main result of this paper is the extension of the Schur-Horn Theorem to infinite sequences: 
For two  nonincreasing nonsummable sequences $\xi$ and $\eta$ that converge to 0, 
there exists a positive compact operator $A$ with eigenvalue list $\eta$ and diagonal sequence $\xi$ if and only if $\sum_{j=1}^n\xi_j \le \sum_{j=1}^n\eta_j$ for every $n$ if and only if $\xi = Q\eta$ for some orthostochastic matrix $Q$. 
When  $\xi$ and $\eta$ are summable, requiring additionally equality of their infinite series obtains the same conclusion, extending a theorem by Arveson and Kadison. 
Our proof depends on the construction and analysis of an infinite product of T-transform matrices. 
Further results on majorization for infinite sequences providing ``intermediate" sequences generalize known results from the finite case. 
Majorization properties and invariance under various classes of stochastic matrices are then used to characterize arithmetic mean closed operator ideals.

\end{abstract}

\maketitle


\section{Introduction}\label{S: 1}

The Schur-Horn Theorem characterizes the diagonals of a (finite) selfadjoint matrix  in terms of sequence majorization, that is, the order relation $\xi\preccurlyeq \eta$ for $\xi, \eta \in \mathbb R^N$ given by   $ \sum_{j=1}^n \xi^*_j \le \sum_{j=1}^n \eta^*_j $ for $1\le n\le N$ and $ \sum_{j=1}^N \xi_j = \sum_{j=1}^N \eta_j $, where $\xi^*, \eta^*$ denote the monotone nonincreasing rearrangement of  $\xi, \eta$.
The theory of majorization arose during the early part of the 20th century from a number of apparently unrelated  topics: wealth distribution  (Lorenz \cite {Lm1905}), inequalities involving convex functions (Hardy, Littlewood and P\'olya \cite {HLP52}), convex combinations of permutation matrices (Birkhoff \cite{Bg46}), and more central to our interests herein, doubly stochastic matrices and the relation between eigenvalue lists and diagonals of selfadjoint matrices: 

\bT{T:1.1} Let $ \xi, \eta \in \mathbb R^N$. 
\item[(i)]  \textbf{Hardy, Littlewood and P\'olya  Theorem} \cite {HLP52}. $\xi\preccurlyeq \eta$ if and only $\xi=P\eta$ for some doubly stochastic matrix $P$.
\item[(ii)]  \textbf{Horn Theorem} \cite [Theorem 4]{Horn}.  $\xi\preccurlyeq \eta$  if and only if $\xi = Q\eta$ for some  orthostochastic matrix $Q$, i.e., the Schur-square of  an orthogonal matrix ($Q_{ij}=(U_{ij})^2~ \forall ~ i,j$ for some unitary matrix  $U$ with real entries).
\item[(iii)] \textbf{Schur-Horn Theorem} \cite {Si23},\cite{Horn}. Given a selfadjoint $N\times N$ matrix $A$ having eigenvalue list $\eta$, there is a basis in which $A$ has diagonal entries $\xi$ if and only if $\xi\preccurlyeq \eta$.
\eT

The sufficiency part of the Schur-Horn Theorem is due to Schur and the necessity follows immediately from the Horn Theorem. The main goal of this paper is to extend to infinite dimension the Horn Theorem and hence  the Schur-Horn Theorem. 

The notion of  majorization extends seamlessly to infinite sequences that admit a monotone nonincreasing  rearrangement and in this paper, to avoid always having to pass to monotone rearrangements, we will focus on sequences decreasing monotonically to 0 and will denote by  $\co*$  their positive cone and  by $(\ell^1)^*$ the subcone of summable decreasing sequences. (We note explicitly that $\co*$ and $(\ell^1)^*$ do not mean herein the duals of $\text{c}_{\text{o}}$ and $\ell^1$.)  
Even for finite sequences, the terminology and notations describing majorization vary considerably in the literature. In this paper, we will use the following notations:

\bD{D:1.2}
For $\xi, \eta \in \co*$ we say that
\begin{itemize}
\item 
$\eta$ \textit{majorizes} $\xi$ ($\xi \prec \eta$) if $ \, \sum_{j=1}^n \xi_j \le \sum_{j=1}^n \eta_j $ for every $n\in \mathbb N; $
\item 
 $\eta$ \textit{strongly majorizes} $\xi$ ($\xi\preccurlyeq \eta) 
\text{ if } \xi \prec \eta  \, \text { and } {\varliminf  }\sum_{j=1}^n ( \eta_j -\xi_j) = 0$;
\item $\eta$ \textit{block majorizes} $\xi$  ($\xi\prec_b  \eta)  
\text { if } \xi \prec \eta  \, \text { and } \sum_{j=1}^{n_k} \xi_j = \sum_{j=1}^{n_k} \eta_j$ 
for some sequence $ \mathbb N \ni n_k \uparrow \infty$.
\end{itemize}
\noindent For $\xi, \eta \in (\ell^1)^*$ we say that
\begin{itemize}
\item $\eta$ \textit{majorizes at infinity} $\xi$  ($ \xi \prec_\infty\eta$) 
if $ \sum_{j=n}^\infty \xi_j \le \sum_{j=n}^\infty \eta_j$ for every $ n\in \mathbb N;$
\item  $\eta$ \textit{strongly majorizes at infinity} $\xi$  ($ \xi \preccurlyeq_\infty \eta$) 
if $\xi \prec_\infty\eta $ and 
$ \sum_{j=1}^\infty  \xi_j = \sum_{j=1}^\infty\eta_j$.
\end{itemize}
If $\xi, \eta \in \text{c}_{\text{o}}^+ $ (resp. $ (\ell^1)^+$), then we say that any of the above relations hold for $\xi$ and $\eta$ if they hold for their monotone rearrangements $\xi^*$ and $\eta^*$. 
\eD

For nonsummable monotone decreasing sequences, the condition 
${\varliminf}\sum_{j=1}^n (\eta_j - \xi_j) = 0$
 retains many of  the key properties of "equality at the end"  for finite and for summable sequences  (e.g., see Section \ref {S: 7}) and it will prove  crucial for our extension of the Schur-Horn Theorem.

Majorization at infinity, aka ``tail majorization,"  was first introduced and studied  for finite sequences and 
appears in \cite{AGPS87}, \cite {nK89}, \cite{DFWW}, \cite {mW02} among others. 
It holds particular relevance for this paper as  it provides the natural characterization of the notion of arithmetic mean closure at infinity for operator ideals contained in the trace class (see Theorem \ref {T:8.3} and \cite{vKgW02}-\cite {vKgW07-OT21}).
 
Block-majorization is both a natural way to bring the results of finite majorization theory to bear on infinite sequences and it also arises naturally in Sections \ref {S:4},\ref{S: 7}, and \ref{S: 8}. See in particular the proof of our extension of a known finite intermediate sequence result:
\begin{itemize}
\item If $\xi,\eta\in \co*$ and $\xi\prec\eta$, then there are $\zeta, \rho \in \co*$ for which $\xi\preccurlyeq \zeta \le \eta$ and $\xi \le \rho\preccurlyeq \eta$ (Theorem \ref {T:7.4}).
\end{itemize} 
and of a similar result  for majorization at infinity (Theorem \ref {T:7.7}). 

In 1964, in two papers that are not nearly as well-known as they deserve and with two almost disjoint approaches, Markus \cite {aM64} and Gohberg and Markus  \cite {GiMa64}  found infinite dimensional versions of  the Hardy, Littlewood and P\'olya Theorem \ref{T:1.1}(i) and an extension to the summable case of the Horn Theorem  \cite [Theorem 4]{Horn}, (Theorem \ref {T:1.1} (ii)).   More recently, Arveson and Kadison obtained  other characterizations in \cite {AK02} using still different methods. 
\begin{alignat}{4}
\text{If } \xi, \eta &~ \in \co*, \text{ ~then} \phantom{pppppppppppppppppppppppppppppppppppppppppppppppppppppppppppppppppppppppppppppp}\notag\\
&\xi \prec \eta  \Longleftrightarrow \xi = Q\eta \text{ for some  substochastic  $Q$ ($Q$ row-stochastic if  $\xi_n>0 ~\forall n$)  \cite[Lemma 3.1]{aM64}  \label{e:1}}\\
&\xi \prec \eta \Longleftrightarrow \xi = Q\eta \text{ with } Q_{ij}=|W_{ij}|^2, \text{for  some co-isometry $W$ \cite [Proposition III, pg 205]{GiMa64}}\notag\\
&\xi \prec \eta \Longleftrightarrow \xi = Q\eta \text{ with } Q_{ij}=|L_{ij}|^2, \text{for  some contraction $L$ \cite [Theorem 4.2]{AK02}}\notag\\
\text{If  } \xi,  \eta &~ \in (\ell^1)^+, \text{~ then}\notag\\
& \xi \preccurlyeq \eta  \Longleftrightarrow \xi = Q\eta \text{ with } Q_{ij}=|U_{ij}|^2, \text{for  some unitary $U$ \cite [Theorem 1]{GiMa64}}\label{e:2}\\
& \xi \preccurlyeq \eta  \Longleftrightarrow \xi = Q\eta \text{ with } Q_{ij}=|W_{ij}|^2, \text{for  some isometry $W$ \cite [Theorem 4.1]{AK02}}\notag
 \end{alignat}
 
By reformulating matricially the proof of  Markus's \cite [Lemma 3.1] {aM64} and slightly tightening it (see Remark \ref {R:3.8}),  we can identify a sequence of orthogonal matrices underlying the construction. An analysis of their properties and infinite products permits us to obtain: 
\begin{itemize}
 \item If $\xi, \eta \in \co*, ~\xi_n>0~ \forall n, \text{ and } \xi \prec \eta$, then there is a  canonical co-isometry with real entries $W(\xi, \eta)$ for which $\xi = Q(\xi, \eta)\mspace{1mu}\eta$  with $Q(\xi, \eta)_{ij}= (W(\xi, \eta)_{ij})^2$ (Theorem \ref {T:3.7}).
\item 
If $\xi, \eta \in \co*$, then $\xi \preccurlyeq \eta \Longleftrightarrow  \xi = Q\eta$ with $Q_{ij}=|U_{ij}|^2$ for  some orthogonal matrix $U$  (Theorem \ref {T:3.9}).
\end{itemize}
Not surprisingly,  this construction applied to finite sequences provides another proof of the (finite) Horn Theorem \ref {T:1.1}(ii).

The canonical matrix  $Q(\xi, \eta)$ is obtained as an infinite product of T-transforms (see Section \ref{S:3} for details) and is therefore completely determined by the sequence $\{m_k,t_k\} $ where $m_k$ is the matrix size and $0< t_k\le1$ is the ``convex coefficient"  of the k-th transform.
In Section \ref{S:4} we further analyze   this double sequence and, more precisely the set $\{t_k \mid m_k =1\}$, in order to link properties of the majorization $\xi \prec \eta$ with  properties of the corresponding canonical matrix $Q(\xi, \eta)$:
\begin{itemize}
\item $\xi \preccurlyeq \eta $ if and only if  $Q(\xi, \eta)$ is orthostochastic  (i.e.,  if $W(\xi, \eta)$ is orthogonal) if and only if \linebreak $\sum \{t_k \mid m_k =1\} = \infty$ (Theorem \ref {T:4.7}). 
\item
$\xi\prec_b  \eta$ if and only  $Q(\xi, \eta)$ is the direct sum of finite orthostochastic matrices if and only if   \linebreak card$\{k \mid t_k= m_k =1\} =\infty$ (Proposition \ref {P:4.4}).
\end{itemize}
Notice that if  $\xi, \eta \in \co*$ and $\xi= Q\eta $ for some orthostochastic matrix $Q$, 
then by \eqref{e:1} it follows that $\xi \prec \eta$, 
but in general it does not follow that $\xi\preccurlyeq \eta$. 
In fact, a main results of this paper is:
\begin{itemize} \item If $\xi, \eta \in \co*$, $\xi\not\in (\ell^1)^*$ and $\xi\prec\eta$,  then there is an orthostochastic matrix $Q$ for which $\xi=Q\eta$ (Theorem \ref {T:5.3}). \end {itemize} 
Section \ref {S: 5} is devoted to the proof of this theorem by showing  that a pair of $\co*$-sequences $\xi\prec\eta$ with $\xi$ nonsummable and 
$\xi\not\preccurlyeq\eta$ can be decomposed into ``mutually  orthogonal"  pairs of infinite subsequences for which strong majorization holds 
(Lemmas \ref {L:5.1}-\ref {L:5.2}) and then invoking Theorem \ref {T:3.9} to obtain an  orthostochastic matrix for each pair and direct summing them.
Together,  Theorems \ref {T:3.9} and \ref {T:5.3} provide the following infinite dimensional extension of the Horn Theorem (Remark \ref{R:5.4}):
\begin{itemize}
\item If  $\xi, \eta \in \co*$ then $\xi=Q\eta $ for some orthostochastic matrix $Q \Longleftrightarrow \begin{cases} \xi \prec \eta &\text{ if } \xi \not\in \ell^1\\
 \xi \preccurlyeq \eta &\text{ if } \xi \in \ell^1.
\end{cases}$
\end{itemize} 
To apply the Horn Theorem to positive compact operators, notice first that the  eigenvalue list with multiplicity (which  is the sequence $s(A)$ of s-numbers of $A\in K(H)^+$) ``ignores" the nullspace of $A$ (e.g., see (\ref {e:32})) and hence it characterizes the partial isometry orbit $\mathscr V(A):=\{VAV^*\mid V \text{ partial isometry}, V^*VA=A\}$ of $A$ rather than, as in the finite rank case, the unitary orbit   $\mathscr U(A)$ of $A$.  Then if we fix an orthonormal basis of the Hilbert space $H$ and denote by $E$ the canonical conditional expectation on  the corresponding atomic masa $\mathscr D$  (i.e., the operation of ``taking the main diagonal"), we obtain the following infinite dimension extension of the Schur-Horn Theorem for 
positive compact operators:
\begin{itemize}
\item 
$E(\mathscr V(A))=\begin{cases}
  \{B\in \mathscr D\cap K(H)^+\mid s(B) \prec s(A) \}\setminus \mathscr L_1 \quad &\text{if } \tr(A)= \infty\\
   \{B\in \mathscr D\cap K(H)^+ \mid s(B) \preccurlyeq s(A)\} &\text{if }  \tr(A)< \infty.\\
  \end{cases}
$ (Proposition \ref{P:6.4} ). 
\end{itemize}

If $A$ has finite rank, then $\mathscr U(A)=\mathscr V(A)$ and if $A\in K(H)^+$ has dense range, i.e., $R_A=I$ then 

\begin{itemize}
\item $E(\mathscr U(A))=\begin{cases}
  \{B\in \mathscr D\cap K(H)^+\mid s(B) \prec s(A), R_B=I \}\setminus \mathscr L_1 \quad &\text{if } \tr(A)= \infty\\
   \{B\in \mathscr D\cap K(H)^+ \mid s(B) \preccurlyeq s(A), R_B=I\} &\text{if }  \tr(A)< \infty.\\
  \end{cases}$ (Proposition \ref {P:6.6} )
\end{itemize}

For positive compact operators with infinite rank some sufficient conditions and some necessary conditions for membership in $E(\mathscr U(A))$ are presented in Propositions  \ref {P:6.6} and \ref {P: 6.10}. 
Our work extends some of the results of Gohberg and Marcus in \cite {GiMa64} and Arveson and Kadison in \cite {AK02}. There are only limited overlaps between our paper and those  by A. Neumann \cite {Na99} and by Antezana,  Massey,  Ruiz, and Stojanoff \cite{AMRS} as these authors characterize the closures of the expectation of the unitary orbit of a selfadjoint not necessarily compact operator while we do not take closures. The connections with  these three papers are further discussed in Section \ref{S: 6} where we also answer a couple of questions of Neumann.  The following,  in the case of sequences $\xi,\eta \in \co*$.  compares  these different results to Proposition \ref{P:6.4} .   

\begin{alignat*}{3}
&\text{If } \xi &&\not \in (\ell^1)^*, \text{ then } \xi\prec \eta  \Leftrightarrow  \diag \xi \in 
&&\begin{cases}  
\overline{E(\mathscr U( \diag \eta ))}^{||.||}  &\text{(\cite [Theorem 3.13] {Na99})}\\
E\{L \diag \eta \1L^* \mid L\in B(H)_1\} & \text{(\cite[Theorem 4.2]  {AK02})}\\
E(\mathscr U( \diag \eta )) &\text{(Proposition \ref{P:6.4} )}
\end{cases} \\
&\text{If } \xi &&\in (\ell^1)^*, \text{ then } \xi\preccurlyeq \eta  \Leftrightarrow \diag \xi \in 
&&\begin{cases} 
 \overline{E(\mathscr U( \diag\eta ))}^{||.||_1}   &\qquad \qquad \quad \text{(\cite[Proposition 3.13] {AMRS})}\\
E(\overline{\mathscr U( \diag \eta)}^{||.||_1} ) &\qquad \qquad \quad \text{(\cite[Theorem 4.2]{AK02})}\\
E(\mathscr U( \diag \eta))  &\qquad \qquad \quad \text{(\cite[Theorem 1]{GiMa64}, Proposition \ref{P:6.4} )}
\end{cases} \notag
\end{alignat*}

Finally, in Section \ref {S: 8} we apply results on majorization for infinite sequences to operator ideals 
(two-sided ideals of $B(H)$). We characterize arithmetic mean closed ideals and arithmetic mean at infinity closed ideals (see Section \ref {S: 8} for the definitions) in terms of the diagonals of their positive elements (Theorems \ref {T:8.1}, \ref{T:8.3}) and in terms of invariance under the action of various classes of (sub) stochastic matrices (Corollary \ref {C:8.7}, Theorem \ref {T:8.9}).
This paper began as part of a long-term project investigating arithmetic mean ideals and arithmetic mean at infinity ideals \cite{vKgW02}-\cite{vKgW07-OT21}.

\section{Notations and  preliminaries on  stochastic matrices} \label{S:2} 

Let $\co*$ denote the cone of nonnegative monotone nonincreasing sequences converging to 0 and
$(\ell^1)^*$ the  cone of nonnegative monotone nonincreasing summable sequences. 
(Again, notice that $\co*$ and $(\ell^1)^*$ here do not denote the duals of $\text{c}_{\text{o}}$ and $\ell^1$.) If $\xi\in (\text{c}_\text{o})^+$, denote by $\xi^*$ its nonincreasing rearrangement.

For every sequence $\xi = \,<\xi_1,\xi_2, \ldots>$ and every $n = 0, 1,  \ldots\,$, denote by $\xi^{(n)}$ the \textit{truncated} sequence 
$\xi^{(n)}=\,<\xi_{n+1},\xi_{n+2}, \ldots>$ and by $\xi\chi[1,n]$ the sequence $<\xi_1, \xi_2,\ldots, \xi_{n}, 0, \ldots>$. 
We will  of course identify $\xi\chi[1,n]$ with a vector in $\mathbb R^n$  and conversely, embed $\mathbb R^n$ into  $\text{c}_{\text{o}}$ by completing finite sequences with infinitely many zeros. 

When applying  the majorization notations of Definition \ref {D:1.2} to finite sequences, we caution the reader again that  what we call majorization ($\xi \prec \eta$, i.e., $\sum_1^k \xi_j \le \sum_1^k \eta_j$ for all $1\le k \le n$) is often called weak majorization, 
and what we call strong majorization ($\xi \preccurlyeq \eta$, i.e., $\xi \prec \eta$ with $\sum_{j=1}^n \xi_j = \sum_{j=1}^n \eta_j$) is mostly called majorization, although with no universal agreement about notations or even about the direction of the inequalities (see \cite [Remark, page 198] {HJ85}).
For the theory of majorization of finite sequences we refer the reader to  Marshall and Olkin  \cite {MO79}.

Immediate consequences of Definition \ref {D:1.2} are:
\begin{align}\label{e:3}
&\text{If  $\xi, \eta \in \co*$, then}\quad \xi\prec_b  \eta \Rightarrow \xi\preccurlyeq \eta  \Rightarrow \xi \prec \eta.\\
&\text{If  $\xi, \eta \in (\ell^1)^*$, then} \quad \xi\prec \eta \text{ and }\sum_{j=1}^\infty  \xi_j = \sum_{j=1}^\infty\eta_j\Leftrightarrow \xi\preccurlyeq \eta\, \Leftrightarrow \, \eta \preccurlyeq_\infty\xi.\notag
\end{align}

Once we have fixed an orthonormal basis $\{e_k\}$ for a complex separable infinite-dimensional Hilbert space $H$, i.e., once we have identified $H$ with $\ell^2$, we will also identify infinite matrices with operators and will use these terms interchangeably. E.g., when we apply a Hilbert space operator to a sequence in $\co*$, what we mean is  that we apply the corresponding matrix to that sequence -- which for substochastic matrices is always possible (e.g., see Remark \ref {R:2.2}).  Also, for typographical reasons we are not going to distinguish between row and column vectors, e.g., if $P$ is a matrix, we shall write $P<\xi_1, \xi_2, \ldots>$ in lieu the  more precise $P<\xi_1, \xi_2, \ldots>^{\text{T}}$.

$K(H)$ denotes the ideal of compact operators and $\mathscr L_1$ the trace class ideal, with $tr$ denoting the usual trace.    Given a  compact operator $A\in K(H)$, the sequence $s(A)\in \co*$ of its s-numbers (singular numbers) consists of the eigenvalues of $(A^*A)^{1/2}$ in monotone order, with repetition according to multiplicity, and with infinitely many zeros added in case $A$ has finite rank. In particular, if $A\ge 0$ has infinite rank, then $s(A)$ is precisely the ``eigenvalue list" of $A$.

Given a sequence $\rho \in \ell^\infty$, we denote by $\diag \rho$ the diagonal matrix having $\rho$ as its main diagonal. Given an operator $A\in B(H)$, we denote by $E(A)$ the diagonal matrix having as diagonal the main diagonal of $A$, i.e.,  $E: B(H) \to \mathscr D$  is the normal faithful and trace preserving conditional expectation from $B(H)$ onto the masa $\mathscr D$ of  diagonal operators. \\

Stochastic matrices play a key role in majorization theory of finite sequences (e.g, see Theorem \ref{T:1.1}(i)). A similar  but necessarily more complex role is played in the case of infinite sequences.

\bD{D:2.1}
A matrix $P$ with nonnegative entries is called
\begin{itemize}
\item substochastic if its row and column sums are bounded by 1;
\item column-stochastic if it is substochastic with column sums equal to 1;
\item row-stochastic if it is substochastic with row sums equal to 1;
\item  doubly stochastic if it is both column and row-stochastic;
\item block stochastic if it is the direct sum of finite doubly stochastic matrices.
\end{itemize}
\eD

\bR{R:2.2}
\item [(i)] Contrary to the finite case a (square) matrix can be column-stochastic without being row stochastic and vice versa.
\item [(ii)] We can apply a substochastic matrix $P$ to any sequence $\rho \in \ell^\infty$, where by $P\rho$ we just mean the sequence $ < \sum_{j=1}^\infty P_{ij}\rho_j>$.
\item[(iii)] If $\rho\in \text{c}_{\text{o}}$ and $P$ is substochastic, then $P\rho \in  \text{c}_{\text{o}}$.
\item[(iv)] By Schur's test (e.g., see \cite [Problem 45]{pH82}) substochastic matrices viewed as Hilbert space operators are contractions.  
\eR

An important class of stochastic matrices is the one obtained as ``Schur-squares" of contractions. To be more precise, we should call them  the Schur product of a contraction by its complex conjugate matrix, but in most cases we consider only matrices with real entries. The Schur product is also called  Hadamard product or entrywise product. The relevance of these stochastic matrices is clear from the following well-known lemma whose verification is straightforward.

\bL{L:2.3} Let $\xi, \eta \in \co*$ and let $Q_{ij} = |L_{ij}|^2$ for all $i,j$ for some contraction $L\in B(H)$. Then  \linebreak $\diag \xi = E(L\diag\eta \, L^*)$ if and only if $\xi = Q \eta$.
\eL

\bL{L:2.4}  Let $Q_{ij} = |L_{ij}|^2$ for all $i,j$ for some contraction $L\in B(H)$. Then
\item[(i)]  $Q$ is substochastic. 
\item[(ii)] $Q$ is column-stochastic if and only if $L$ is an isometry.  
\item[(ii$'$)] $Q$  is row-stochastic if and only if $L$ is a co-isometry. 
\eL
\bp Notice first that 
\be{e:4}
\sum_{j=1}^\infty Q_{ij}= \sum_{j=1}^\infty L_{ij}L^*_{ji} = ||L^*e_i||^2\le 1 \quad \text{for every  } i,
\ee
and similarly
 \be{e:5}\sum_{i=1}^\infty Q_{ij} = ||Le_j||^2\le 1\quad \text{for every  } j.  
 \ee
\item[(i)] Immediate from (\ref {e:4}) and (\ref {e:5}).
\item[(ii)] Sufficiency is immediate from (i) and (\ref {e:5}). 
Conversely, assume that $Q$ is column-stochastic  and  hence  $||Le_j||=1$ for all $j$ by (\ref {e:5}). 
Then $(L^*Le_j,e_j)=1$ for all $j$ and thus $E(I-L^*L)=0$. 
Since $E$ is faithful and $I-L^*L\ge 0$ because $L$ is a contraction, it follows that $L^*L=I$.
\item[(ii$'$)] Apply (ii) to $L^*$.
\ep

\bD{L:2.5}  If $Q_{ij} = |L_{ij}|^2$ for some contraction $L$, then we say that $Q$ is  isometry stochastic (resp. co-isometry stochastic, unistochastic, orthostochastic) if $L$ is an isometry (resp. co-isometry, unitary, orthogonal matrix, i.e., a unitary matrix with real entries). 
If $L$ is the direct sum of finite unitary (resp. orthogonal) matrices, we say that $Q$ is block unistochastic (resp. block orthostochastic.)
\eD

\bR{L:2.6} 
\item[(i)] The terminology ``orthostochastic'' goes back at least to Horn (cf.  \cite {Horn}). 
When the entries of the unitary matrix are not necessarily real, its Schur-square is called  unitary stochastic in \cite {MO79},  Pythagorean in  \cite [Section 4] {Kr02b},  and orthostochastic  in \cite{HJ91}, although unistochastic appears to be the more common term now.
\item[(ii)]  $Q$ can be the Schur-square of different contractions, e.g.,  $Q_{ij} = |L^{(1)}_{ij}|^2=| L^{(2)}_{ij}|^2$; but $L^{(1)}$ is an isometry, a co-isometry, a unitary if and only if so is $L^{(2)}$, respectively.  Of course, $L^{(1)}$ may have real entries, while $L^{(2)}$ does not.
\item[(iii)] $L$ does not  need to be a contraction for $Q$ to be doubly stochastic, 
e.g., consider  the $4 \times 4$ matrix $L$ with constant entries $\frac{1}{2}$.
\item[(iv)] As remarked by Horn \cite {Horn}, every $2 \times 2$ doubly stochastic matrix is necessarily orthostochastic but 
\[
\frac{1}{2}
\begin{pmatrix}
1&1&0\\
1&0&1\\
0&1&1
\end{pmatrix}
\]
is doubly stochastic but not  unistochastic (see \cite[p. 39]{MO79} for more examples).
\item[(v)] Let $P$ be a matrix and $\Pi$ be a permutation matrix. 
Then $\Pi P$ is substochastic (resp. row-stochastic, column-stochastic, isometry stochastic, co-isometry stochastic, unistochastic, orthostochastic) precisely when $P$ is. 
However, permutations do not preserve block stochasticity.
\eR

An immediate consequence  of Theorem \ref {T:1.1}(i) and (ii) is that
\begin{align}\label{e:6}
\text{if }\xi, \eta \in \co*,\text{ then }\xi\prec_b  \eta   & \Leftrightarrow \xi = Q\eta \text{ for a block orthostochastic matrix }Q\\
& \Leftrightarrow \xi = Q\eta \text{ for a block stochastic matrix }Q.\notag
\end{align}

Another simple application of the Horn Theorem, which we will need in Theorem  \ref {T:3.9}, is to the case when $\eta$ is a sequence with finite support. This result generalizes \cite [Theorem 13]{Kr02b} (see also \cite[Theorem 4.7]{AMRS}.)

\bL{L:2.7} If  $\xi,\eta\in\co*$ , $\xi\preccurlyeq\eta$, and $\eta$ has finite support, then  $\xi=Q\eta$ for some  orthostochastic matrix $Q$.
\eL
\bp
The case when $\eta=0$ being trivial, let $n$ be the largest integer for which $\eta_n>0$. If $n=1$, then let  $U$ be an orthogonal matrix that has as its first column the unit vector 
$\begin{pmatrix}
\sqrt{\frac{\xi_1}{\eta_1}},
\sqrt{\frac{\xi_2}{\eta_1}},
\sqrt{\frac{\xi_3}{\eta_1}}, 
\cdots
\end{pmatrix}^{\text {T}}$
and let $Q_{ij}:= U_{ij}^2$. Then $Q\eta = \xi$.  In the case that  $n>1$,  we have $\sum_{j=1}^\infty \xi_j = \sum_{j=1}^n \eta_j >  \sum_{j=1}^{n-1} \eta_j$. Let $m$ be the index for 
which $\sum_{j=1}^{m-1} \xi_j \le  \sum_{j=1}^{n-1} \eta_j <  \sum_{j=1}^{m} \xi_j$ and let $\alpha:= \sum_{j=1}^{m} \xi_j- \sum_{j=1}^{n-1} \eta_j$. Then $m\ge n$, $0 < \alpha = \eta_n-\sum_{j=m+1}^\infty\xi_j \le \eta_n$ and 
\[
<\xi_1,\xi_2,\ldots,\xi_m>\,\preccurlyeq\,<\eta_1,\eta_2, \ldots, \eta_{n-1}, \alpha, 0,\dots,0>.
\]
By applying the Horn Theorem  (see Theorem \ref {T:1.1} (ii))  to the above two vectors of $ \mathbb (R^m)^+$, we find an $m\times m$ orthostochastic matrix $Q^o$ for which 
\[
Q^o<\eta_1,\eta_2, \ldots, \eta_{n-1}, \alpha, 0, \dots,0>\,=\,<\xi_1,\xi_2,\ldots,\xi_m>
\]
 and let $U^o$ be an orthogonal matrix for which $Q_{ij}^o=(U_{ij}^o)^2$.  In particular, the first $n$ columns $U_1, U_2, \ldots,U_n$ of $U^o$ are orthonormal.  Denote by  $Q_1, Q_2, \ldots,Q_n$ their Schur-squares, i.e., the first $n$ columns of $Q^o$. Therefore the $\ell^2$ vectors
\[
\begin{pmatrix}
U_1\\
0\\
0\\
\vdots
\end{pmatrix},
\begin{pmatrix}
U_2\\
0\\
0\\
\vdots
\end{pmatrix},
\dots , 
\begin{pmatrix}
U_{n-1}\\
0\\
0\\
\vdots
\end{pmatrix},
\begin{pmatrix}
\sqrt\frac{\alpha}{\eta_n}U_n\\
\sqrt\frac{\xi_{m+1}}{\eta_n}\\
\sqrt\frac{\xi_{m+2}}{\eta_n}\\
\vdots
\end{pmatrix}
\]
are also orthonormal. Complete them to an orthonormal basis of $\ell^2$  with real entry vectors and denote by $U$ the orthogonal matrix having as columns these vectors and by $Q$ the orthostochastic matrix  $Q_{ij}=U_{ij}^2$. Then 
\[
Q=\begin{pmatrix}
Q_1&\dots&Q_{n-1}&\frac{\alpha}{\eta_n}Q_n&*&*&\dots\\
0&\dots&0&\frac{\xi_{m+1}}{\eta_n}&*&*&\dots\\
0&\dots&0&\frac{\xi_{m+2}}{\eta_n}&*&*&\dots\\
\vdots&\vdots&\vdots&\vdots&\vdots&\vdots&\vdots
\end{pmatrix}
\]
and hence
\[
Q\eta=
\begin{pmatrix}
\begin{pmatrix}
Q_1&\dots&Q_{n-1}&\frac{\alpha}{\eta_n}Q_n
\end{pmatrix}
\begin{pmatrix}
\eta_1\\
\eta_2\\
\vdots\\
\eta_n\\
\end{pmatrix}\\
\begin{pmatrix}
0&\dots&0&\frac{\xi_{m+1}}{\eta_n}.\\
0&\dots&0&\frac{\xi_{m+2}}{\eta_n}\\
\vdots&\vdots&\vdots&\vdots
\end{pmatrix}
\begin{pmatrix}
\eta_1\\
\eta_2\\
\vdots\\
\eta_n\\
\end{pmatrix}\\
\end{pmatrix}\\
=
\begin{pmatrix}
\begin{pmatrix}
Q_1&\dots&Q_{n-1}&Q_n
\end{pmatrix}
 \begin{pmatrix}\\
\eta_1\\
\eta_2\\
\vdots\\
\eta_{n-1}\\
\alpha\\
\end{pmatrix}\\
\begin{pmatrix}
\xi_{m+1}\\
\xi_{m+2}\\
\vdots\\
\end{pmatrix}\\
\end{pmatrix}
=\xi.
\]
\ep

The following lemma is a  key ``bridge" between  properties of majorization and properties of stochastic matrices. 

\bL{L:2.8}  If $P$ is a  substochastic matrix for which ${\varliminf}_{n}\sum_{i=1}^n (\eta_i -(P\eta)_i )=0$  for some $\eta\in\co*$ with $\eta_n > 0$ for all $n$, then $P$ is column-stochastic.
\eL
\bp
\begin{align*}
\sum_{i=1}^n (\eta_i -(P\eta)_i )&= \sum_{j=1}^n\big(1- \sum_{i=1}^n P_{ij}\big)(\eta_j-\eta_n)
+\sum_{j=n+1}^{\infty}\sum_{i=1}^n P_{ij}(\eta_n-\eta_j) +
\big(n-  \sum_{i=1}^n \sum_{j=1}^{\infty}P_{ij}\big)\eta_n\\
&\ge  (1- \sum_{i=1}^n P_{ij})(\eta_j-\eta_n)\ge 0\quad \text{for all  $n \ge  j.$}
\end{align*}
Thus for all $j$,
\[
0= {\varliminf}_{_n}\sum_{i=1}^n (\eta_i -(P\eta)_i )\ge {\varliminf}_{_n}\ (1- \sum_{i=1}^n P_{ij})(\eta_j-\eta_n) = (1- \sum_{i=1}^\infty P_{ij})\eta_j \ge 0,
\]
hence $ \sum_{i=1}^\infty P_{ij}=0$.
\ep

\bR{R:2.9}\item[(i)] The first line of the proof is based  on Ostrowski's decomposition \cite{Oa52} and shows that $\sum_{i=1}^n (\eta_i -(P\eta)_i )\ge0$ whether $P\eta$ is monotone or not.  It  was used by Markus to prove that $P\eta \prec \eta$ in \cite[Lemma 3.1]{aM64}. 
\item[(ii)] If  $P\eta$ is monotone, then the condition $P\eta\preccurlyeq \eta$ is equivalent to ${\varliminf}_{_n}\sum_{i=1}^n (\eta_i -(P\eta)_i )=0$  and hence implies that $P$ is column-stochastic.
\item[(iii)] In the case that $\zeta\in( \ell^1)^+$ and $\sum_{j=1}^\infty (P\zeta)_j =\sum_{j=1}^\infty \zeta_j $, it is immediate to verify that $\sum_{=1}^\infty P_{jn}=1$ for every $n$ which $\zeta_n\ne 0$. The same conclusion can be obtained by the operator theoretic argument in \cite[Theorem 4.1] {AK02} in the case that $P_{ij}=|L_{ij}|^2$ for some contraction $L$.  For the reader's convenience, a sketch of the argument is that  then $\tr(\diag \zeta)= \tr (E(L\diag\zeta\, L^*))$ and hence $\tr\Big(E\big((\diag \zeta)^\frac{1}{2}(I-L^*L)(\diag \zeta)^\frac{1}{2}\big)\Big)=0$ Since $I-L^*L\ge0$ and $\tr$ and $E$ are faithful,  it follows that  $(\diag \zeta)^\frac{1}{2}(I-L^*L)(\diag \zeta)^\frac{1}{2}=0$   and hence $||Le_n||=1$ for every $n$ for which $\zeta_n\ne 0$.
\item[(iv)] Notice that if $P$ is a substochastic matrix for which $ P\eta\preccurlyeq \eta$ for some $\eta\in\co*$ with $\eta_n > 0$ for all $n$, $P$ can fail to be row-stochastic as is the case for 
\[
P:=
\begin{pmatrix}
1/2&0&0&0&\dots\\
1/2&0&0&0&\dots\\
0&1/2&0&0&\dots\\
0&1/2&0&0&\dots\\
0&0&1/2&0&\dots\\
0&0&1/2&0&\dots\\
\vdots&\vdots&\vdots&\vdots&\vdots
\end{pmatrix}
\]
and $\eta$ any summable sequence with  $\eta_n > 0$ for all $n$. 
\eR
For summable sequences the converse of Lemma \ref {L:2.8}   holds.

\bL{L:2.10} Let $\xi, \eta \in \co*$ and $\xi = P\eta$ for some  column-stochastic matrix $P$.  If  $\xi \in (\ell^1)^*$ (resp. $\eta\in(\ell^1)^*$), then $\eta \in (\ell^1)^*$ (resp.  $\xi \in (\ell^1)^*$) and $\xi \preccurlyeq \eta$.
\eL
\bp
 We know from \eqref{e:1} that $\xi\prec\eta$. Moreover,
\[
\sum_{i=1}^\infty \xi_i= \sum_{i=1}^\infty (P\eta)_i=\sum_{i=1}^\infty \sum_{j=1}^\infty P_{ij}\eta_j = \sum_{j=1}^\infty \sum_{i=1}^\infty P_{ij}\eta_j = \sum_{j=1}^\infty \eta_j,
\]
thus $\xi \in (\ell^1)^*$  if and only if $\eta \in (\ell^1)^*$  and $\xi \preccurlyeq \eta$.
\ep

Without the condition of summability, the implication in Lemma \ref{L:2.10} can fail. 
In fact, the following example shows that it can\,fail even for an orthostochastic matrix
as seen by modifying the matrix in Remark \ref {R:2.9}(iv)  as follows. Let $\omega$ denote the harmonic sequence, i.e.,  $\omega:=\,<\frac{1}{n}>$. 

\begin{example}\label{E:2.11}An orthostochastic matrix $Q$ for which  $Q\omega \not \preccurlyeq \omega$.
\end{example}
\bp
Partition $\mathbb N$ into two infinite strictly increasing sequence $ \{n_k\}$ and $\{m_k\}$ for which 
 \linebreak ${\varliminf} \Big(\sum_{j=k+1}^{2k} \frac{1}{n_j}- \sum_{j=1}^{k}\frac{1}{m_j}\Big) > 0$. 
For instance, this can be achieved by setting $m_k:=(k+1)^2$ and listing $\mathbb N \setminus \{m_k\}$ as $\{n_k\}$. \\

Defining
\[
U_{ij}=\begin{cases} \frac{1}{\sqrt{2}}\quad &\text{if }\,  i=2k-1\,,j=n_k, \, k\in \mathbb N\\
 \frac{1}{\sqrt{2}}\quad &\text{if }\,  i=2k-1 ,j=m_k , \,  k\in \mathbb N\\
 \frac{1}{\sqrt{2}}\quad &\text{if }\,  i=2k ,j=n_k, \, k\in \mathbb N\\
 - \frac{1}{\sqrt{2}}\quad &\text{if }\,  i=2k ,j=m_k,\,   k\in \mathbb N\\
0&\text{otherwise}
\end{cases}.
\]
it is easy to see that $U$ is an orthogonal matrix. Let $Q$ be the Schur-square of $U$, i.e., $Q_{ij}:=U_{ij}^2$. Then a simple computation shows that for every $k\in\mathbb N$,
\[
(Q\omega)_{2k -1} = (Q\omega)_{2k } = \frac{1}{2}\big(\frac{1}{n_k}+\frac{1}{m_k}\big)
\]
and hence $(Q\omega)_{2k}  > (Q\omega)_{2k +1}$, 
that is, $Q\omega$ is monotone nonincreasing. Moreover,
\[
\sum_{j=1}^{2k}\omega_j- \sum_{j=1}^{2k}(Q\omega)_j = 
\sum_{j=1}^{2k}\frac{1}{j}- \sum_{i=1}^{k}\big(\frac{1}{n_i}+\frac{1}{m_i}\big)\\
\ge \sum_{j=k+1}^{2k} \frac{1}{n_i}- \sum_{j=1}^{k}\frac{1}{m_i}.
\]
Similarly, 
\begin{align*}
\sum_{j=1}^{2k+1}\omega_j- \sum_{j=1}^{2k+1}(Q\omega)_j &= 
\sum_{j=1}^{2k+1}\frac{1}{j}- \sum_{i=1}^{k}\big(\frac{1}{n_i}+\frac{1}{m_i}\big) -\frac{1}{2}( \frac{1}{n_{k+1}}+\frac{1}{m_{k+1}})\\
&\ge  \sum_{j=k+1}^{2k} \frac{1}{n_i}- \sum_{j=1}^{k}\frac{1}{m_i} -\frac{1}{2}( \frac{1}{n_{k+1}}+\frac{1}{m_{k+1}}).
\end{align*}
Thus $
{\varliminf} \left(\sum_{j=1}^n\omega_j- \sum_{j=1}^n(Q\omega)_j\right) \ge {\varliminf} \Big(\sum_{j=k+1}^{2k} \frac{1}{n_j}- \sum_{j=1}^{k}\frac{1}{m_j}\Big) > 0 $, i.e.,  $Q\omega \not \preccurlyeq \omega$.
\ep

\noindent  We will see from Theorem \ref{T:5.3} that for any nonsummable sequence $\eta$ we can choose an orthostochastic matrix $Q$ for which $Q\eta= \frac{1}{2}\eta$ and hence  $Q\eta \not \preccurlyeq \eta$.

We know of no simple condition that characterizes substochastic matrices for which $P\eta \preccurlyeq \eta$ for all $\eta\in\co*$. Notice that since $P\eta $ is not necessarily monotone,  by the latter condition we mean $(P\eta)^* \preccurlyeq \eta$ for the monotone  rearrangement $(P\eta)^*$ of $P\eta $. A sufficient condition is that  $P$ is block stochastic, i.e., the direct sum of  doubly stochastic finite matrices. A  more general sufficient condition is provided by the following lemma.

\bL{L:2.12}
If $P$ is a substochastic matrix and $ \varliminf \, (n-\sum_{i,j=1}^nP_{ij})=0$, then
\item[(i)] $P$ is doubly stochastic;
\item[(ii)]  $P\eta \preccurlyeq \eta$  for every $\eta\in\co*$.
\eL
\bp
\item[(i)]
\[
 n-\sum_{i,j=1}^n P_{ij}=  \sum_{i=1}^n (1- \sum_{j=1}^n P_{ij}) \ge \sum_{i=1}^n (1- \sum_{j=1}^\infty P_{ij})
 \]
thus
\[
0= {\varliminf} \,(n-\sum_{i,j=1}^n P_{ij})\ge \sum_{i=1}^\infty (1- \sum_{j=1}^\infty P_{ij}). 
 \]
Then because $P$ is substochastic, $\sum_{j=1}^\infty P_{ij} =1$ for all $i$.
Similarly,  $\sum_{i=1}^\infty P_{ij} = 1$ for all $j$.
\item[(ii)] Since $\sum_{i=1}^n (P\eta)_i \le  \sum_{i=1}^n (P\eta)_i^*$ for every $n$,
\begin{align*}
\sum_{i=1}^n (\eta_i-(P\eta)^*_i)&\le \sum_{i=1}^n (\eta_i-(P\eta)_i) \\
&=  \sum_{j=1}^n (1- \sum_{i=1}^n P_{ij})\eta_j - \sum_{i=1}^n\sum_{j=n+1}^\infty P_{ij}\eta_j \\
&\le  \sum_{j=1}^n (1- \sum_{i=1}^n P_{ij})||\eta||_\infty \\
&= (n-\sum_{i,j=1}^nP_{ij})||\eta||_\infty,
\end{align*}
hence ${\varliminf}\, \sum_{i=1}^n (\eta_i-(P\eta)_i)^* \le {\varliminf} \,(n-\sum_{i,j=1}^nP_{ij})||\eta||_\infty =0.$ 
Thus $(P\eta)^* \preccurlyeq \eta$, i.e., by Definition \ref {D:1.2}, $P\eta \preccurlyeq \eta$.

\ep

\bR{R:2.13}
The condition $ {\varliminf} \,(n-\sum_{i,j=1}^nP_{ij})=0$ is not necessary for (ii). 
For instance, (ii) holds trivially for every permutation matrix $\Pi$ since  $(\Pi\eta)^*=\eta$ for every $\eta$, but it is easy to find a permutation matrix $\Pi$ for which ${\varliminf} (n-\sum_{i,j=1}^n\Pi_{ij})>0 $. 
\eR

\section{The canonical co-isometry of a majorization} \label{S:3}
We start with some historical notes about the link between majorization and stochastic matrices.

Muirhead  \cite {Mu1903} for the case of integer-valued finite sequences and then Hardy, Littlewood and P\'olya \cite [p. 47]{HLP52} for the case of real-valued finite sequences proved that for all  $\xi,\eta\in \mathbb R^N$ with $\xi \preccurlyeq \eta $, there is  a doubly stochastic matrix $P$ with $\xi=P\eta$. $P$ was obtained as a finite product  of T-transforms, i.e. matrices of the form $tI+(1-t)\Pi$ with $\Pi$ a transposition and $0\le t<1$.  The T-transforms  were chosen so to reduce at each step  the Hamming distance (i.e., the number of positions where two sequences differ) between the sequence $\xi$ and the iterated transform of $\eta$.  Notice that while individual T-transforms are orthostochastic, the product of two T-transforms can fail to be even unistochastic (\cite [Chapter II, Section G] {MO79}).

In 1952, Horn proved that the matrix $P$ can be chosen to be orthostochastic by using a different method based on convexity arguments and a technically difficult proof \cite[Theorem 4]{Horn}.

A  proof based on properties of determinants was given a few years later by Mirsky \cite[Theorem 2] {Mi58}.

After a four decade hyatus, a  proof based on composition of Givens rotations (special permutations of T-transforms) was obtained by Casazza and Leon  in the Appendix of \cite {CL02}.

More recently,  Arveson and Kadison \cite [Theorem 2.1]{AK02} gave an elegant  proof of the Horn theorem showing that $P$ can be chosen to be unistochastic (see also \cite[Lemma 5 and Theorem 6]{Kr02a}). Reformulating their result in our terminology, they showed that $\xi$ is obtained by applying to $\eta$ a finite number of T-transforms and that by properly choosing unitary matrices whose Schur-squares are those T-transforms,  the Schur-square of their product (a unistochastic matrix by definition)  applied to $\eta$ also yields $\xi$.  

Another recent proof was obtained by Kornelson and Larson  \cite [Theorem 2] {KkLd04}. More precisely, they proved the equivalent statement that  every positive finite rank operator $B$ with eigenvalue list $ \eta$ can be decomposed as the linear combination $B=\sum_{j=1}^k\xi_j P_j$ of  rank-one projections (not required to be mutually orthogonal)  with the given monotone nonincreasing coefficient sequence $\xi:= \, <\xi_1,  \xi_2 ,  \ldots   \xi_k, 0\ldots>$,  if and only if (in our notations) $\xi\preccurlyeq \eta$.

This link between majorization and stochastic matrices  was partially extended to the infinite case by Markus in \cite[Lemma 3.1]{aM64} (see \eqref {e:1}) and the Schur-Horn Theorem was extended to summable sequences by Gohberg and Markus in \cite [Lemma 1] {GiMa64} (see  \eqref{e:2}) based on \cite[Theorem 1]{GiMa64}.  The latter proof depended crucially on the summability of the sequence, so we focus on the former proof.

 At the core of Markus's proof, although he did not employ this terminology nor exhibit explicitly the matrices, is the construction, for every $\xi, \eta \in \co*$ with  $\xi \prec \eta $,  of an infinite sequence of permutations of \mbox{T-transforms} whose product is a substochastic matrix $Q$ for which $\xi = Q\eta$. Furthermore,  a remark in his proof states that when  $\xi_n >0$ for all $n$, the matrix $Q$ is row-stochastic. 
In this section we revisit and slightly tighten   the Markus construction for the case when  $\xi_n >0$ for all $n$ (see Remark \ref {R:3.8}) and prove that it provides a co-isometry stochastic matrix (Theorem \ref {T:3.7}) and  in the strong majorization case, an orthostochastic matrix (Theorems \ref {T:3.9}, \ref {T:4.7}). Not surprisingly, this construction restricted to finite sequences yields another proof of the Horn Theorem (Remark \ref {R:4.3})

For every integer $m\ge 1$ and $0< t\le1$, define the $m+1 \times m+1$ orthogonal matrix

\begin{equation}\label{e:7}
V(m,t):=
\begin{pmatrix}
0&0&\dots&0&\sqrt {t}&-\sqrt{1-t}\\
1&0&\dots&0&0&0&\\
0&1&\dots&0&0&0&\\
\vdots&\vdots&\vdots&\vdots&\vdots&\vdots\\
0&0&\dots&1&0&0\\
0&0&\dots&0&\sqrt {1-t}&\sqrt{t}\\
\end{pmatrix}.
\end{equation}
Notice that the Schur-square of the matrix $V(m,t)$ is the product of the permutation matrix $\Pi_1$ that sends $<1,2,\ldots, m,m+1>$ to $<m, 1,2,\ldots, m-1,m+1>$ and of the T-transform $tI+(1-t)\Pi_2$ where  $\Pi_2$ is the transposition that interchanges $m$ and $m+1$.  

Given a sequence $\{m_k,t_k\}$ where $m_k\in \mathbb N$ and  $0< t_k\le1$, define for every $n\in \mathbb N$
\be{e:8}
W^{(n)}:= \big(I_{n-1}\oplus V(m_n,t_n)\oplus I_\infty\big)\big(I_{n-1}\oplus V(m_n,t_n)\oplus I_\infty\big)\cdots \big(I_1\oplus V(m_2,t_2)\oplus I_\infty\big)\big( V(m_1,t_1)\oplus I_\infty\big)
\ee
where $I_n$ denotes the $n\times n$ identity matrix for $1\le n\le\infty$ and $I_0$ is simply dropped. 
Define also  $R^{(n)}$ to be the Schur-square of $V(m_n,t_n)\oplus I_\infty$, and let 
\be {e:9} 
Q^{(n)}: = \big(I_{n-1}\oplus R^{(n)}\big)\big(I_{n-2}\oplus R^{(n-1)}\big)\cdots  \big(I_{1}\oplus R^{(2)}\big) R^{(1)}.
\ee
Being a product of orthogonal matrices, all  $W^{(n)}$ are also orthogonal. Denote by  $P_n:= I_n\oplus \,0$  the projection on span$\{e_1,\ldots, e_n\}$. 

\bP{P:3.1}
Let  $\{m_k,t_k\}$ be a sequence with $m_k\in \mathbb N$ and  $0< t_k\le1$. Then 
\item[(i)] The sequence of operators $W^{(n)}$ converges in the weak operator topology to a co-isometry 
$W(\{m_k,t_k\})$ and $P_nW(\{m_k,t_k\})= P_nW^{(n)}$ for every $n$. 
The convergence is in the strong operator topology if and only if $W(\{m_k,t_k\})$ is orthogonal.
\item[(ii)] The sequence of operators $Q^{(n)}$ converges in the weak operator topology to a row-stochastic operator $Q(\{m_k,t_k\})$ and $P_nQ(\{m_k,t_k\})= P_nQ^{(n)}$ for every $n$. 
\eP
\bp
\item[(i)]  From (\ref {e:8}) we have for  all integers $j > n$,
\[
 W^{(j)}= (I_{j-1}\oplus V(m_j,t_j)\oplus I_\infty)\cdots(I_n\oplus V(m_{n+1},t_{n+1})\oplus I_\infty) W^{(n)}
 \]
and hence
 \be{e:10} 
 P_nW^{(j)} = P_nW^{(n)}
 \quad \text{for all } j\ge n. 
\ee
As a consequence, $\big(( W^{(j)}- W^{(i)})x, y \big) = \big(( W^{(j)}- W^{(i)})x,P_n^\perp y\big)$ 
for all $x, y \in H$ and $i,j \ge n$. Thus the sequence of orthogonal matrices  $\{W^{(j)}\}$ is weakly Cauchy and hence converges weakly  to a contraction $W(\{m_k,t_k\})$ with real entries.  Set $W:=W(\{m_k,t_k\})$. From (\ref {e:10}), it follows that   $P_nW = P_nW^{(n)}$ for all $n$, that is, the first $n$ rows of the matrix $W^{(n)}$ stabilize: 
$W_{ij}= W^{(n)}_{ij} $ for all $n,j$ and all $i\le n.$ 
Therefore 
\[
WW^*=\slim P_nWW^*P_n= \slim P_nW^{(n)}(W^{(n)})^*P_n= \slim P_n = I,
\]
i.e., $W$ is a co-isometry.
Since $ ||(W-W^{(n)})x||^2 \to ||x||^2-||Wx||^2$ for all $x\in H$, the sequence $W^{(n)}$ converges strongly to $W$ if and only if $W$ is orthogonal.
\item[(ii)] Same proof.
\ep

\bR{R:3.2}
Proposition \ref {P:3.1} holds even if we allow $t_k=0$. However, in order to obtain the uniqueness of the sequence $\{m_k,t_k\}$ in the construction in Theorem \ref {T:3.7}, we will have to assume there that $t_k >0$ for all $k$. In addition, this assumption will simplify some of the proofs.
\eR

A case where it is simple to find the form of $W(\{m_k,t_k\}) $ and $Q(\{m_k,t_k\}) $  is when  $m_k=1$ for all $k$.

\begin{example}\label{E:3.3}
 \begin{align*} &W(\{1,t_k\})=\\
&\begin{pmatrix}
\sqrt{t_1}&-\sqrt{1-t_1}&0&0&\dots&\\
\sqrt{t_2(1-t_1)}&\sqrt{t_2t_1}&-\sqrt{1-t_2}&0&\dots\\
\sqrt{t_3(1-t_2)(1-t_1)}&\sqrt{t_3(1-t_2)t_1}&\sqrt{t_3t_2}&-\sqrt{1-t_3}&\dots\\
\vdots&\vdots&\vdots&\vdots&\vdots\\
\sqrt{t_k\prod_{i=1}^{k-1}(1-t_i)}&\sqrt{t_kt_1\prod_{i=2}^{k-1}(1-t_i)}&\sqrt{t_kt_2\prod_{i=3}^{k-1}(1-t_i)}&\sqrt{t_kt_{3}\prod_{i=4}^{k-1}(1-t_i)}&\dots\\
\vdots&\vdots&\vdots&\vdots&\vdots\\
\end{pmatrix}
\end{align*}
\[Q(\{1,t_k\})=
\begin{pmatrix}
t_1&1-t_1&0&0&\dots&\\
t_2(1-t_1)&t_2t_1&1-t_2&0&\dots\\
t_3(1-t_2)(1-t_1)&t_3(1-t_2)t_1&t_3t_2&1-t_3&\dots\\
\vdots&\vdots&\vdots&\vdots&\vdots\\
t_k\prod_{i=1}^{k-1}(1-t_i)&t_kt_1\prod_{i=2}^{k-1}(1-t_i)&t_kt_2\prod_{i=3}^{k-1}(1-t_i)&t_kt_{3}\prod_{i=4}^{k-1}(1-t_i)&\dots\\
\vdots&\vdots&\vdots&\vdots&\vdots\\
\end{pmatrix}
\]
\vspace{1cm}
\end{example}
We see in this case that $Q(\{1,t_k\})$ is the Schur-square of $W(\{1,t_k\})$.  For the general case we first state a couple of elementary lemmas leaving their proof to the reader.

\bL{L:3.4} Let $A$ and $B$ be two bounded matrices with the property that  for every $i,j$ there is at most one index $k$ for which $A_{ik}B_{kj}\ne 0$ and let $A'$ and $B'$ be the Schur-square of $A$ and $B$ respectively. Then $A'B'$ is the Schur-square of $AB$. If furthermore $A$ and $B$ are orthogonal (resp. unitary) then $A'B'$ is orthostochastic (resp. unistochastic). In particular,  if $Q$ is orthostochastic (resp. unistochastic) and $\Pi$ is a permutation, then $\Pi Q$ is  orthostochastic (resp. unistochastic).
\eL

Next,  we consider a simple case where this sufficient condition is satisfied.

\bL{L:3.5}Let  $A$ and $B$ be bounded matrices, let $n\in \mathbb N$,  and assume that  for every  $j$ there is at most one index $i> n$ for which  $B_{ij}\ne 0$. 
\item[(i)] For every $i$ and $j$ there is at most one index $k$ for which $(I_n\oplus A)_{ik}B_{kj}\ne 0$.
\item [(ii)] If for every  $j$ there is at most one index $i> 1$ (resp. $i> 0$) for which $A_{ij}\ne 0$, then for every  $j$ there is at most one index $i> n+1$ (resp. $i > n$)  for which $\big((I_{n}\oplus A)B\big)_{ij}\ne 0$.
\item[(iii)] Let $A(n)$ be a sequence of bounded matrices for which for every  $n$ and $j$, $A(n) _{ij}\ne 0$ for at  most one index $i> 1$ and let  $Q= \big(I_{n-1}\oplus A(n)\big)\big(I_{n-2}\oplus A(n-1)\big)\cdots  \big(I_{1}\oplus A(2)\big) A(1)$. Then for every $j$, $Q_{ij}\ne0$ for at most one index $i>n$.
\eL

Given a sequence $\{m_k,t_k\}$ where $m_k\in \mathbb N$  and  $0< t_k\le1$, then  for every $n, j$,  at most one of the entries $W^{(n)}_{ij}$ for $i>n$ is nonzero. 

\bP{P:3.6}
Let  $\{m_k,t_k\}$ be a sequence with $m_k\in \mathbb N$ and  $0< t_k\le1$ and let $W^{(n)}$ and $Q^{(n)}$ be as  in \eqref {e:8} and \eqref {e:9}. Then 
\item [(i)] $Q^{(n)}$ is the Schur-square of $W^{(n)}$ for every $n$ and for every $n$ and $j$ there is at most one index $i>n$ for which $Q^{(n)}_{ij}\ne 0.$
\item [(ii)] $Q(\{m_k,t_k\})$ is the Schur-square of $W(\{m_k,t_k\})$. 
\eP

\bp
\item[(i)]  We reason by induction. $Q^{(1)}$ is by definition the Schur-square of $W^{(1)}$.  Assume that $Q^{(n)}$ is the Schur-square of $W^{(n)}$.  Now $W^{(n+1)}=\big(I_{n}\oplus V(m_{n+1}, t_{n+1})\oplus I_\infty \big) W^{(n)}$. Since for every factor $V(m_k,t_k)\oplus I_\infty$ and every $j$ there is at most one index $i>1$ for which $\big(V(m_k,t_k)\oplus I_\infty\big)_{ij}\ne 0$, it follows from Lemma \ref {L:3.5}(iii)  that for every $n$ and $j$, $W^{(n)}_{ij}\ne0$ for at most one index $i>n$. Thus  by Lemma \ref {L:3.5}(i),  for every $i$ and $j$, $(I_{n}\oplus V(m_{n+1}, t_{n+1})\oplus I_\infty \big) _{ik}W^{(n)}_{kj}\ne 0$ for at most one index $k$. But then, the product of  $Q^{(n+1)}= \big(I_{n}\oplus R^{(n+1)}\big)Q^{(n)}$ of the Schur-square of $I_{n}\oplus V(m_{n+1}, t_{n+1})\oplus I_\infty $ by the Schur-square of  $W^{(n)}$ coincides by Lemma \ref {L:3.4}  with the Schur-square of  $W^{(n+1)}=\big(I_{n}\oplus V(m_{n+1}, t_{n+1})\oplus I_\infty \big) W^{(n)}$.
\item[(ii)] Obvious since the first $n$ rows of $W(\{m_k,t_k\})$ (resp. $Q(\{m_k,t_k\})$) coincide with the first n rows of $W^{(n)}$ (resp. $Q^{(n)}$.) 
\ep

To every majorization $\xi \prec \eta$  with $\xi_n>0$ for all $n$, the following construction associates a sequence  $\{m_k,t_k\}$ with $m_k\in \mathbb N$ and  $0< t_k\le1$ and hence associates the  corresponding co-isometry $W(\xi, \eta):=W(\{m_k,t_k\})$.

 \bT{T:3.7}
Let $\xi, \eta \in \co*$ with $\xi_n > 0$ for every $n$. 
If $\xi \prec \eta $, then there is a canonical co-isometry $W(\xi, \eta)$ with real entries whose Schur-square  $Q(\xi, \eta)$ satisfies $\xi= Q(\xi, \eta) \eta$.
 \eT
 \bp
We construct the following sequence $\{m_k,t_k\}$ where $m_k\in \mathbb N$ and  $0<t_k\le1$. Set $\rho(0):=\eta$ and choose   $m_1 \in \mathbb N$ for which $\eta_{m_1+1}< \xi_1  \le \eta_{m_1}$. 
Since $\xi\prec \eta$ and hence $\xi_1\le\eta_1$, and since  $\xi_1>0$ and $\eta_j\to 0$, such an integer exists and  by the monotonicity of $\eta$, it is unique. Express $\xi_1$ as a convex combination of $\eta_{m_1}$ and $\eta_{m_1+1}$, that is, choose $t_1$ for which 
$\xi_1= t_1\eta_{m_1}+(1-t_1)\eta_{m_1+1}$. Thus  $0<t_1\le 1$ and also $t_1$ is uniquely determined. 
Set $\delta_1:=(1- t_1)\eta_{m_1}+t_1\eta_{m_1+1} $
and hence  $\delta_1= \eta_{m_1}+\eta_{m_1+1} -\xi_1$. Define the sequence $
\rho(1): = \,<\eta_1, \eta_2, \ldots, \eta_{m_1-1},  \delta_1, \eta_{m_1+2},\eta_{m_1+3}, \ldots>$  
where if $m_1=1$, then the first entry of $\rho(1)$ is  $\delta_1$.
Since $\eta_{m_1+2}\le \eta_{m_1+1}\le \delta_1< \eta_{m_1}\le \eta_{m_1-1},
$
we  see that $\rho(1)$ is monotone nonincreasing and $\rho(1)\le \eta$.  
Let $R^{(1)}$ be the Schur-square of $V(m_1,t_1)\oplus I_\infty$, i.e., $R^{(1)}_{ij}= \big(\big(V(m_1,t_1)\oplus I_\infty\big)_{ij}\big)^2$ for all $i,j$. Then $R^{(1)}\eta = \,<\xi_1, \rho(1)_1,\rho(1)_2, \ldots >.$
Moreover, 
\be{e:11}
\xi^{(1)}\prec \rho(1).
\ee
(Recall the notation $\xi^{(1)}:= \,<\xi_2, \xi_3, \ldots>$.) Indeed, for every $1\le n < m_1$,
\[
 \sum_{j=1}^n(\rho(1)_j- \xi^{(1)}_j) =  \sum_{j=1}^n\eta_j-  \sum_{j=2}^{n+1}\xi_j =  \sum_{j=1}^n(\eta_j-\xi_j) +\xi_1-\xi_{n+1}\ge 0,
 \]

\noindent and for every $n \ge m_1$

\[
 \sum_{j=1}^n(\rho(1)_j- \xi^{(1)}_j) =  \sum_{j=1}^{m_1-1}\eta_j + \delta_1+ \sum_{j=m_1+2}^{n+1}\eta_j
 -\sum_{j=2}^{n+1}\xi_j = \sum_{j=1}^{n+1}(\eta_j-\xi_j) \ge 0.
\]

\noindent Repeat the construction applying it to the pair  $\xi^{(1)}\prec \rho(1)$, and so on. By the assumption that $\xi_k > 0$ for all $k$, the process can be iterated  providing an infinite sequence of pairs $\{m_k,t_k\}$ with $m_k\in \mathbb N$  and  $0< t_k\le1$ and from these,  of   sequences $ \rho(k)$ and scalars  $\delta_k$  satisfying for all $k$ the relations:

\begin{align} 
&\rho(k-1)_{m_k+1}< \xi_k\le \rho(k-1)_{m_k}, \quad \xi_k=t_k\rho(k-1)_{m_k}+(1-t_k)\rho(k-1)_{m_k+1},\label{e:12}\\
&\delta_k:=(1-t_k)\rho(k-1)_{m_k}+t_k\rho(k-1)_{m_k+1}= \rho(k-1)_{m_k}+\rho(k-1)_{m_k+1}-\xi_k,\label{e:13} \\
&\rho(k):= \,<\rho(k-1)_1,\ldots, \rho(k-1)_{m_k-1},\delta_k, \rho(k-1)_{m_k+2},\ldots>,\quad \text{i.e., }\notag\\
&\rho(k)_j=\begin{cases}
\rho(k-1)_j\quad&\text{for all } j<m_k\\
\delta_k &\text{for  } j=m_k\\
 \rho(k-1)_{j+1} &\text{for all } j>m_k,
 \end{cases}\label{e:14}\\
&\xi^{(k)}\prec \rho(k), \quad \eta=\rho(0) \ge \rho(1)\ge \rho(2)\ge\cdots,\label{e:15}\\
& \sum_{j=1}^n(\rho(k)_j- \xi^{(k)}_j) =  \sum_{j=1}^{n+1}(\rho(k-1)_j- \xi^{(k-1)}_j) \quad\text{for all } n\ge m_k.\label{e:16}
\end{align}
Let $R^{(n)}$ be the Schur-square of $V(m_n,t_n)\oplus I_\infty$, and let $Q^{(n)}:=  (I_{n-1}\oplus R^{(n)})\cdots (I_1\oplus R^{(2)})R^{(1)}$  as in \eqref {e:9}. Then  for all $n$,

\be{e:17}
R^{(n)}\rho(n-1)=\,<\xi_n, \rho(n)_1, \rho(n)_2, \ldots>, \quad Q^{(n)} \eta = \,<\xi_1,\xi_2, \ldots,\xi_n,\rho(n)_1, \rho(n)_2, \ldots>
\ee
Let 
\[
W(\xi, \eta):= W(\{m_k,t_k\})\quad\text{and}\quad  Q(\xi, \eta):= Q(\{m_k,t_k\}).
\]

Then by Propositions \ref {P:3.1} and \ref {P:3.6}, $W(\xi, \eta)$ is a co-isometry and $Q(\xi, \eta)$ is its Schur-square.
Finally,  by Remark \ref {R:2.2}, $Q(\xi, \eta)\eta$ is defined and is a sequence in $\text{c}_{\text{o}}$.
In fact,  by (\ref{e:17}), 
\[
P_nQ(\xi, \eta)\eta= P_nQ^{(n)}\eta= \,<\xi_1,\xi_2, \ldots,\xi_n,0,0,\ldots> \,\to \xi ~\text{pointwise},
\]
and hence $Q(\xi, \eta)\eta=\xi$.

\ep
\pagebreak

 \bR{R:3.8}
\item[(i)]The construction of the sequence $\{m_k,t_k\}$  and the associated sequence of matrices $Q^{(n)}$ follows the Markus construction in  \cite [Lemma 3.1] {aM64}. A minor difference is that  while Markus  chose $m_k$ to be an index for which $\eta_{m_k+1}\le \xi_k\le \eta_{m_k}$ so to treat at the same time  also the case when $\xi$ is finitely supported, here we consider only the case of infinitely supported  $\xi$ and then request that $\eta_{m_k+1}< \xi_k\le \eta_{m_k}$, which makes the construction canonical. 
The main difference is that Markus's analysis is at the level of the action of the matrices $Q^{(n)}$ on $\eta$, and  thus yields only that their limit $Q$ is row-stochastic. It is by introducing the underlying matrices $W^{(n)}$ and analyzing their properties that we can obtain  that $Q$ is co-isometry stochastic.
\item[(ii)] As a  consequence of   \cite [Proposition III, pg 205] {GiMa64} obtained by Gohberg and Markus with different methods, is that  if $\xi \prec \eta $,  then $\xi = Q\eta$ for some co-isometry stochastic matrix $Q$ (see Remark \ref{R:6.5}(ii) for more details.)
\eR
 
When the majorization  is strong,  we obtain the following extension of  the Horn Theorem \cite [Theorem 4]{Horn}(see Theorem \ref{T:1.1}(ii)).  In the nonsummable case strong majorization will not be required, as we will see in Theorem \ref  {T:5.3}.

\bT{T:3.9}
If $\xi, \eta \in \co*$ and $\xi \preccurlyeq \eta $, then  $\xi = Q\eta$ for some orthostochastic matrix $Q$.
\eT
\bp
If $\eta$ has finite support, then the conclusion follows from  Lemma \ref {L:2.7}.  If $\eta$ has infinite support, then  $\xi$ too has  infinite support. Indeed, if  otherwise $\xi_n =0$ for some  $n$, then $\sum_{j=1}^\infty \eta_j=  \sum_{j=1}^{n-1} \xi_j\le \sum_{j=1}^{n-1} \eta_j$ which implies that $\eta_n=0$, a contradition. But then, by Theorem   \ref {T:3.7}, $\xi = Q(\xi, \eta)\eta$ where $Q(\xi, \eta)$ is the Schur-square of the co-isometry $W(\xi, \eta).$ By Lemma \ref {L:2.8}  and Remark \ref {R:2.9}(ii), $Q(\xi, \eta)$  is column-stochastic  and hence by Lemma \ref {L:2.4} , $W(\xi, \eta)$ is also an isometry and hence unitary. Since by construction $W(\xi, \eta)$ has real entries, it is orthogonal, hence  $Q(\xi, \eta)$ is orthostochastic.
\ep

\bR{R:3.10}
\item [(i)]  The above proof shows that if  $\xi \preccurlyeq \eta $ and $\eta$ has infinite support, then any co-isometry stochastic matrix $Q$ for which  $\xi = Q\eta$  must be unistochastic, i.e., the Schur product of a unitary matrix by its complex conjugate.
\item [(ii)] In the case when $\xi \preccurlyeq \eta $ and $\xi$ has infinite support but $\eta$ does not, we cannot invoke Lemma \ref {L:2.8}  to conclude that $W(\xi,\eta)$ is orthogonal, so for simplicity's sake, we have chosen in lieu of $Q(\xi, \eta)$ the orthostochastic matrix provided by Lemma \ref {L:2.7}.  However, in the next section we will prove that  if  $\xi \preccurlyeq \eta$, then $Q(\xi, \eta)$ itself is orthostochastic (Theorem \ref {T:4.7}).
 \item[(iii)] If $\xi$ has finite support, say $\{1, \ldots, N\}$, then the construction of Theorem  \ref {T:3.7}  can still be carried on for the first $N$ steps and it provides  yet another proof of the Horn Theorem (see Remark \ref {R:4.3}).  
\eR

 \section{Properties of the canonical matrix $Q(\xi, \eta)$ of a majorization} \label{S:4} 
 In this section, on which the following ones do not depend, we further analyze the construction in Theorem \ref {T:3.7} to relate the properties of the majorization $\xi\prec \eta$ to those of the canonical co-isometry stochastic matrix  $Q(\xi, \eta)$ via the properties of the set $\{t_k \mid m_k=1\}$. In the next lemmas we collect the additional needed properties of the sequences  $m_k$, $t_k $, $\delta_k$, $\rho(k)$,  $W^{(n)}$, etc. that were introduced in Theorem \ref {T:3.7}

\bL{L:4.1}
Let $\xi, \eta \in \co*$,  $\xi \prec \eta$, and assume that $\xi_n>0$ for all $n$. Then for every $k\in \mathbb N$

\item[(i)] $m_k \ge m_{k-1}-1$;
\item[(ii)] $\rho(k)_j= \eta_{j+k}$ for every  $j > m_k$;
\item[(iii)] if $n\ge m_k$, then $\sum_{j=1}^n (\rho(k)_j-\xi^{(k)}_j ) = \sum_{j=1}^{n+k}(\eta_j- \xi_j);$ 
\item[(iv)] if  $m_k=1$, then  $\delta_k= \eta_{k+1}+ \sum_{j=1}^k(\eta_j-\xi_j)$;
\item[(v)] if $t_k=1$, and $m_k= m_{k-1}-1$, then  $t_{k-1}=1$;
\item[(vi)] if  $\sum_{j=1}^n(\eta_j-\xi_j)=0$ and  $\sum_{j=1}^{n-1}(\eta_j-\xi_j)>0$  for some $n>1$,  then $m_n=t_n=1$.
\eL

\bp
\item[(i)] 
Assume by contradiction that $m_k < m_{k-1}-1$, then 
\begin{alignat*}{2}
\xi_k &> \rho(k-1)_{m_k+1}& &\quad \text{(by (\ref{e:12}))}\\
&= \rho(k-2)_{m_k+1}& &\quad \text{(by (\ref {e:14}) since  $m_k+1 < m_{k-1}$)}\\
&\ge \rho(k-2)_{m_{k-1}}& &\quad \text{(by the monotonicity of $ \rho(k-2)$,   since  $m_k+1 < m_{k-1}$)}\\
&\ge \xi_{k-1}& &\quad \text{(by (\ref{e:12}))}.
\end{alignat*}
This is a contradiction because of the monotonicity of $\xi$.

\item[(ii)] The proof is by induction on $k$. The property holds by (\ref {e:14}) for $k=1$ since by definition $\rho(0)=\eta$. Assume it holds for some $k$ and let $j>m_{k+1}$. Then, 
\begin{alignat*}{2}
\rho(k+1)_j&= \rho(k)_{j+1}& &\quad \text{(by (\ref {e:14}))}\\
&=\eta_{j+1+k} & &\quad\text{(by the induction hypothesis, since by (i), $j+1 > m_k.$)}
\end{alignat*}
\item[(iii)]  If $n\ge m_k$, then by (i), $n+p\ge m_{k-p}$ for all $0\le p< k$. Thus iterating (\ref{e:16})
\[
 \sum_{j=1}^n(\rho(k)_j- \xi^{(k)}_j) =   \sum_{j=1}^{n+k}(\rho(0)_j- \xi^{(0)}_j) =  \sum_{j=1}^{n+k}(\eta_j- \xi_j).
\]

\item[(iv)]  
Since $\delta_k= \rho(k)_1$ by (\ref {e:14}), setting $n=1$ in (iii) we obtain
\[
\delta_k= \xi_{k+1} +  \sum_{j=1}^{k+1}(\eta_j - \xi_j) =\eta_{k+1} +  \sum_{j=1}^{k}(\eta_j - \xi_j).
\]

\begin{flalign*}
(\text{v}) \qquad \xi_{k-1}&\le \rho(k-2)_{m_{k-1}}& & \text{(by (\ref{e:12}))}\\
 &\le \rho(k-2)_{m_k}& & \text{(by the monotonicity of $\rho(k-2)$, since $m_k< m_{k-1}$)}\\
 &= \rho(k-1)_{m_k} & &\text{(by (\ref {e:14}), since $m_k< m_{k-1}$)}\\
&= \xi_k & & \text{(by (\ref{e:12}), since $t_k=1$)}\\
 &\le \xi_{k-1} & & \text{(by the monotonicity of $\xi$.)}
\end{flalign*}
But then,  $\xi_{k-1}= \rho(k-2)_{m_{k-1}}$ and hence by (\ref{e:12}), $t_{k-1}=1$.

\item[(vi)]  We reason by induction on $n$ and first prove the property for $n=2$. If  $\eta_1+\eta_2= \xi_1+\xi_2$ and $\eta_1 > \xi_1$, then $\eta_2 < \xi_2 \le \xi_1 < \eta_1$. Thus $m_1=1$, $\rho(1)_1=\delta_1 = \eta_1+\eta_2- \xi_1 = \xi_2$ and hence $m_2=t_2=1$. Assume now that the property (vi) holds  for some $n\ge 2$ for every pair of sequences and  that $\sum_{j=1}^{n+1}(\eta_j-\xi_j)=0$ and $\sum_{j=1}^{n}(\eta_j-\xi_j)>0$. Then  $\eta_{n+1} < \xi_{n+1} \le \xi_1 \le \eta_{m_1}$ implies that $n+1> m_1$ and   from (\ref {e:16}) we obtain that
$ \sum_{j=1}^{n}(\rho(1)_j- \xi^{(1)}_j) =  \sum_{j=1}^{n+1}(\eta_j- \xi_j) = 0.$
We claim that $ \sum_{j=1}^{n-1}(\rho(1)_j- \xi^{(1)}_j)>0.$ 
If  $n > m_1$, the claim holds because
\begin{alignat*}{2}
 \rho(1)_n&=\eta_{n+1}& &\quad \text{(by (\ref {e:14}))}\\
 &< \xi_{n+1}\\
 &=\xi^{(1)}_n& &\quad \text{(by definition)}\\
 &=  \rho(1)_n+  \sum_{j=1}^{n-1}(\rho(1)_j- \xi^{(1)}_j)& &\quad \text{(since } \sum_{j=1}^{n}(\rho(1)_j- \xi^{(1)}_j) =0.)
\end{alignat*}
If  $n = m_1$, i.e., $\eta_{n+1}< \xi_1\le \eta_{n}$, then
\begin{alignat*}{2}
\sum_{j=1}^{n-1}(\rho(1)_j- \xi^{(1)}_j) &=\xi^{(1)}_n-\rho(1)_n& &\quad \text{(since } \sum_{j=1}^{n}(\rho(1)_j- \xi^{(1)}_j) =0)\\
&= \xi_{n+1}  -\delta_1& &\quad \text {(by (\ref{e:14}))}\\
&= \xi_{n+1}  - \eta_{n+1}  + \xi_1- \eta_{n} & &\quad \text {(by (\ref {e:13}))} \\
&=  \sum_{j=1}^{n-1}\eta_j-  \sum_{j=2}^{n}\xi_j& &\quad \text{(since }\sum_{j=1}^{n+1}(\eta_j- \xi_j) = 0.)
\end{alignat*}
Thus, if  $\xi_1= \eta_{n}$, then $\sum_{j=1}^{n-1}(\rho(1)_j- \xi^{(1)}_j) =  \xi_{n+1}  - \eta_{n+1} > 0$. If on the other hand 
$\xi_1< \eta_{n}$, then by the monotonicity of $\eta$ and $\xi$, we have $
 \sum_{j=1}^{n-1}\eta_j>   \sum_{j=1}^{n-1}\xi_j\ge \sum_{j=2}^{n}\xi_j,$
thus completing the proof of the claim.  Therefore the sequences $\xi^{(1)}\prec \rho(1)$  satisfy the hypotheses of (vi) for $n$ and hence, by the induction hypothesis, satisfy the thesis of (vi). But by definition, the pair $\{m_n,t_n\}$ for $\xi^{(1)}\prec \rho(1)$ coincides with the pair $\{m_{n+1},t_{n+1}\}$ for  $\xi \prec \eta$, which concludes the induction proof.

\ep

Without the assumption that  $\sum_{j=1}^{n-1}(\eta_j-\xi_j)>0$, the conclusion of (vi) may fail:  consider for instance $\xi =<1,1,*,\ldots >$ and $\eta = <1,1, 1,0, \ldots>$ where $m_2 = 2$.

\bL{L:4.2}
Let $\xi, \eta \in \co*$ with $\xi_n>0$ for all $n$ and  $\xi \prec \eta$.   
\item[(i)] 
$W^{(n)}=P_{n+m_n}W^{(n)}P_{n+m_n}+P_{n+m_n}^\perp$ for every $n$.  
\item[(ii)] If  $P_n$ commutes with $W(\xi,\eta)$, then $\sum_{j=1}^n(\eta_j-\xi_j)=0.$
\item[(iii)]  If $m_n=t_n=1$, then $P_n$ commutes with $W^{(n)}$ and with  $W(\xi,\eta)$ and $\rho(n)=\eta^{(n)}$.
\eL

\bp
\item[(i)] By Lemma \ref {L:4.1}(i), the sizes $k+m_k$ of the matrices $I_{k-1}\oplus V(m_k,t_k)$  are nondecreasing. Thus  for every $1\le k\le n$,
\[
P_{n+m_n}^\perp \big(I_{k-1}\oplus V(m_k,t_k)\oplus I_\infty\big)= 
 \big(I_{k-1}\oplus V(m_k,t_k)\oplus I_\infty\big)P_{n+m_n}^\perp
=P_{n+m_n}^\perp.
\]
By (\ref {e:8}),  $P_{n+m_n}^\perp W^{(n)}= W^{(n)} P_{n+m_n}^\perp= P_{n+m_n}^\perp$ and hence
the claim.  
     
\item[(ii)] If $W(\xi, \eta)$ commutes with $P_n$, i.e., $W(\xi, \eta)_{ij}=0$ when $1\le i\le n$ and $j>n$ and when $1\le j\le n$ and $i>n$, then so does its Schur-square $Q(\xi, \eta)$. But then, the $n \times n$ matrix  $Q_n:=P_nQ(\xi, \eta)P_n \bigr  | _{P_n H}$ is also orthostochastic. Since $Q(\xi, \eta)\eta=\xi$, it follows that $Q_n<\eta_1, \ldots, \eta_n> = <\xi_1, \ldots, \xi_n>$ and hence  $\sum_{j=1}^n(\eta_j-\xi_j)=0.$

\item[(iii)] 
By Lemma \ref {L:4.1} (i), for every $k$, either $m_{k-1}\le m_k$ or  $m_{k-1}=  m_k+1$. 
Let $j$ be the largest index $i\le n$ for which $m_{i-1}\le m_{i}$ and if there is none, set $j=1$. Then $m_i= n+1-i$ for all $j\le i \le n$. By applying recursively Lemma \ref {L:4.1}(v) we obtain that  $t_i= 1$  for all $j\le i \le n$. But then the size of all the matrices $I_{i-1}\oplus V(n+1-i,1)$ is constant and equal to $n+1$, hence
\[
W^{(n)} =
\bigg(\Big(\big(I_{n-1}\oplus V(1,1))(I_{n-2}\oplus V(2,1))\dots (I_{j-1}\oplus V(n+1-j,1))\Big)\oplus I_\infty\bigg)W^{(j-1)} 
\]
where we set $W^{(0)}=I$ if $j=1$. All the matrices $I_{i-1}\oplus V(n+1-i,1)$ for $j\le i\le n$ are $n+1\times n+1$ permutation matrices that leave the $n+1$ position fixed and hence they commute with $P_n$. If $j=1$, then  $W^{(0)}=I$ commutes trivially with $P_n$, while if $j>1$, then $m_{j-1}=m_j= n+1-j$, hence $n= j-1+m_{j-1}$, and thus by (i), $W^{(j-1)}$ also commutes with $P_n$. Thus $P_nW^{(n)}=W^{(n)}P_n$. 
As $P_nW(\xi, \eta)= P_nW^{(n)}$  by Proposition \ref {P:3.1}, it follows that  $P_nW(\xi, \eta)P_n^\perp=0$. On the other hand, $W(\xi, \eta)$ is a co-isometry and $W^{(n)}$ is unitary, hence
\begin{align*}
P_n^\perp W(\xi, \eta)P_n&= P_n^\perp W(\xi, \eta)(W^{(n)})^*W^{(n)}P_n
= P_n^\perp W(\xi, \eta)(W^{(n)})^*P_nW^{(n)}\\
&= P_n^\perp W(\xi, \eta)W(\xi, \eta)^*P_nW^{(n)}
=P_n^\perp P_nW^{(n)}
=0
\end{align*}
which proves that $W(\xi, \eta)$ commutes with $P_n$. Moreover,
\begin{alignat*}{2}
\rho(n) &= <\delta_n, \eta_{n+2}, \ldots>&& \quad\text{(by (\ref {e:14}) and Lemma \ref {L:4.1}(ii), since $m_n=1$)}\\
&= <\eta_{n+1}+\sum_{j=1}^n(\eta_j-\xi_j),  \eta_{n+2}, \ldots> && \quad\text{(by Lemma \ref{L:4.1} (iv), since  $m_n=1$)}\\
&= \eta^{(n)}&& \quad\text{(by (ii))}
\end{alignat*}
\ep

\bR{R:4.3}[\textbf{Proof of the Horn Theorem}]
If the sequence $\xi$ has finite support, say $\{1, \ldots N\}$, then, as mentioned in Remark \ref {R:3.10} (iii),  the construction in Theorem \ref {T:3.7} can be carried for the first $N$ steps, and the properties obtained in Lemmas \ref {L:4.1} and  \ref {L:4.2} hold for $1\le n\le N$. Thus  $Q^{(N)}$ is an infinite orthostochastic matrix and $Q^{(N)}\eta = < \xi_1, \ldots, \xi_N, \rho(N)_1, \rho(N)_2, \ldots>$. 
If furthermore $\xi \preccurlyeq \eta$, then also $\eta_j=0$ for all $j>N$ and $\rho(N)\equiv 0$. It is then easy to verify that then the upper left $N \times N$ block $Q_N$ of $Q^{(N)}$ is also orthostochastic and that $<\xi_1, \ldots \xi_N > = Q_N <\eta_1, \ldots, \eta_N >$.

Thus if we start with two finite (monotone) sequences $\xi, \eta\in \mathbb R^N$ with  $\xi \preccurlyeq \eta$ we can obtain the required orthostochastic matrix $Q_N$ by applying the construction in Theorem \ref {T:3.7} to $N\times N$ matrices, thus providing an algorithmic proof of the Horn Theorem. For  the reader's convenience we summarize this adaptation.

\eR
\bp
For every integer $1\le m \le n-1<N$ and $0< t\le1$, define the $n \times n$ orthogonal matrix

\[
V(m,t, n):=
\begin{pmatrix}
0&0&\dots&0&\sqrt {t}&-\sqrt{1-t}&0&0&\cdots &0\\
1&0&\dots&0&0&0&0&0&\cdots &0\\
0&1&\dots&0&0&0&0&0&\cdots &0\\
\vdots&\vdots&&\vdots&\vdots&\vdots&\vdots&\vdots&&\vdots\\
0&0&\dots&1&0&0&0&0&\cdots &0\\
0&0&\dots&0&\sqrt {1-t}&\sqrt{t}&0&0&\cdots &0\\
0&0&\dots&1&0&0&1&0&\cdots &0\\
0&0&\dots&1&0&0&0&1&\cdots &0\\
\vdots&\vdots&&\vdots&\vdots&\vdots&\vdots&\vdots&&\vdots\\
0&0&\dots&1&0&0&0&0&\cdots &1\\
\end{pmatrix}.
\]
where the first nonzero entry on the first row occurs in position $m$. Construct the sequence $\{m_k,t_k\}_1^{N-1}$ with 
$1\le m_k\le N-k$ and $0<t_k\le1$ for which $\xi_k=t_k\tilde\rho(k)_{m_k}+(1-t_k)\tilde\rho(k)_{m_{k+1}}$ where $\tilde\rho(k)$ is defined inductively by 
\[
\tilde\rho(k):=<\tilde\rho(k-1)_1, \cdots, \tilde\rho(k-1)_{m_{k-1}}, (1-t_k) \tilde\rho(k-1)_{m_{k}}+ t_k\tilde\rho(k-1)_{m_{k+1}}-\xi_{k-1}, \tilde\rho(k-1)_{m_{k+2}},\cdots,  \tilde\rho(k-1)_{N-k}>
\]
starting with  $\tilde\rho(0):=\eta=< \eta_1, \cdots, \eta_N>$.
Then 
\[
W_N:=  \big(I_{N-2}\oplus V(m_{N-1},t_{N-1}, 2)\big)\big(I_{N-3}\oplus V(m_{N-2},t_{N-2}, 3)\big)\cdots,  V(m_1,t_1, N)
\]
is an orthogonal matrix and its Schur square $Q_N$ satisfies $\xi=Q_N\eta.$
\ep

Now we return to infinite sequences and apply Lemmas \ref {L:4.1} and \ref {L:4.2} to  show that  the set $\{t_k\mid m_k=1\} $  encodes key information about $W(\xi, \eta)$ and $Q(\xi, \eta)$. 

First, we characterize block-majorization (see Definition \ref {D:1.2}), both because it might be of independent interest and because it provides a key step in the characterization of strong majorization. Recall from (\ref{e:6}) that an immediate consequence of the Horn Theorem   is that $\xi\prec_b  \eta$ if and only if  $\xi = Q\eta$ for some block-orthostochastic matrix $Q$.  The next proposition  states  that if $\xi\prec_b  \eta$,  then $Q(\xi,\eta)$ itself must be block-orthostochastic (equivalently, $W(\xi,\eta)$ is the direct sum of finite orthogonal matrices) and characterizes when this occurs in terms of the sequence $\{m_k,t_k\}$. 

\bP{P:4.4} Let  $\xi \prec \eta $  for some $\xi, \eta \in \co*$ with $\xi_n >0$ for every $n$. Then the following conditions are equivalent.
\item[(i)] $\xi\prec_b  \eta$.
\item[(ii)] The set $\{k\mid m_k=t_k=1\}$ is infinite.
\item[(iii)] $Q(\xi,\eta)$ is block orthostochastic.
\eP
\bp
\item[(i)]$ \Rightarrow$ (ii)  By definition, $\sum_{j=1}^{n_k}(\eta_j-\xi_j)=0$ for some strictly increasing sequence $\{n_k\}$.  Then either there is an infinite sequence $p_k$ for which $\sum_{j=1}^{p_k}(\eta_j-\xi_j)=0$ and $\sum_{j=1}^{p_k-1}(\eta_j-\xi_j)>0$ or there is some $N\in \mathbb N$  for which $\sum_{j=1}^{n}(\eta_j-\xi_j)=0$ for all $n\ge N$  and hence  $\eta_j=\xi_j$ for all $j>N$. In the first case $m_{p_k}= t_{p_k} =1$ for all $k $ by Lemma \ref  {L:4.1} (vi). In the second case,  choose  the smallest $N$ with this property.  If  $N>1$,  then $\sum_{j=1}^{N-1}(\eta_j-\xi_j)>0$ and hence $m_N=t_N=1$ by Lemma \ref {L:4.1} (vi) and 
\begin{alignat*}{2}
\rho(N)&= <\delta_N, \eta_{N+2}, \ldots > & & \quad \text{(by (\ref {e:14}) and Lemma \ref {L:4.1} (ii))}\\
&= <\eta_{N+1}, \eta_{N+2}, \ldots > & &\quad ( \text{by Lemma \ref {L:4.1} (iv), since } \sum_{j=1}^{}(\eta_j-\xi_j)=0)\\
&= \eta^{(N)} =  \xi^{(N)}.
\end{alignat*}
If $N=1$ then we see directly that $\xi=\eta$.  It is easy to see now that whether $N=1$ or $N>1$, $\rho(j)= \eta^{(j)}$ for all $j\ge N$. Since $\eta\to 0$ and $\eta$ has infinite support since by assumption and so has $\xi$, there is an infinite collection of indices $j > N$ for which $\rho(j-1)_2=\eta_{j+1} < \xi_j=  \eta_j=\rho(j-1)_1$ and thus for those indices $m_j=1=t_j$.

\item[(ii)]$ \Rightarrow$ (iii) By Lemma \ref {L:4.2}(ii), $W(\xi,\eta)$ commutes with every $P_k$ for which $m_k=t_k=1$. Thus $W(\xi, \eta)$ is block diagonal with each (finite) block an orthogonal matrix and hence its Schur-square $Q(\xi, \eta)$ is block orthostochastic.

\item[(iii)]$ \Rightarrow$ (i) Obvious (see (\ref{e:6})).
\ep
Next, we proceed to characterize  strong majorization  $\xi\preccurlyeq  \eta$.  To do so, we will first need to further analyze the property obtained in Proposition \ref {P:3.6} (i) that for every $n$ the orthogonal matrix $W^{(n)}$ has in each column at most one non-zero entry below row $n$.  For a given $j$, define 
\begin{alignat}{2}\label{e:18}
&q(n,j)= \gamma(n,j)=0\quad \quad\quad &&\text{if all the entries of $W^{(n)}_{ij}$ for $i>n$ are zero}\\
&W^{(n)}_{n+q(n,j), j}=\gamma(n,j) &&\text{is the unique nonzero entry. }\notag
\end{alignat}
Reformulating \eqref{e:18}  in vector form, 
\be{e:19} \begin{pmatrix} 
W^{(n)}_{n+1,j}&\\
W^{(n)}_{n+2,j}&\\
\dots &
\end{pmatrix}  = 
\begin{cases}
0\quad &\text{if } q(n,j) = 0\\
\gamma(n,j)  e_{q(n,j)} & \text{if } q(n,j) \ne 0\
\end{cases}.
\ee
and thus we obtain the recurrence relation 
\begin{align*}
\begin{pmatrix} 
W^{(n+1)}_{n+1,j}&\\
W^{(n+1)}_{n+2,j}&\\
\dots &
\end{pmatrix}
&= \big(V(m_{n+1}, t_{n+1})\oplus I_\infty\big) \begin{pmatrix} 
W^{(n)}_{n+1,j}&\\
W^{(n)}_{n+2,j}&\\
\dots &
\end{pmatrix} \\
&= \begin{cases}
0\quad &\text{if } q(n,j) = 0\\
\gamma(n,j) \big(V(m_{n+1}, t_{n+1})\oplus I_\infty\big)e_{q(n,j)} & \text{if } q(n,j) \ne 0
\end{cases}.\notag
\end{align*}

We leave  to the reader to  verify the following lemma.

\bL{L:4.5}
Given a sequence $\{m_k,t_k\}$ where $m_k\in \mathbb N$  and  $0< t_k\le1$, and the   co-isometry  \linebreak $W:=W(\{m_k,t_k\})$, for every $n, j$, let  $q(n,j)$ and $\gamma(n,j)$ be the sequences  defined by (\ref{e:18})
\item [(i)] 
\bq
q(1,j) = \begin{cases} j\quad &\text{for } j < m_1\\
j  \quad&\text{for }  j=m_1, t_1\ne 1\\
0&\text{for }  j=m_1, t_1= 1\\
j-1&\text{for }  j>m_1
\end{cases},\quad  \quad 
\gamma (1,j) = \begin{cases} 
1\quad &\text{for } j < m_1\\
\sqrt{1-t_1} &\text{for }  j=m_1\\
\sqrt{t_1} &\text{for }  j=m_1+1\\
1&\text{for }  j>m_1+1
\end{cases},
\eq
and 
\begin{align*}
q(n+1,j) &=\begin{cases} q(n,j)\quad\quad \quad \quad \,\, &\text{for } q(n,j) < m_{n+1}\\
q(n,j) &\text{for } q(n,j) = m_{n+1},  t_{n+1}\ne 1\\
0&\text{for }  q(n,j) = m_{n+1},  t_{n+1}= 1\\
q(n,j)-1 &\text{for } q(n,j) > m_{n+1}
\end{cases},\\
\gamma(n+1,j) &=\begin{cases} \gamma(n,j)\quad &\text{for } q(n,j) < m_{n+1}\\
\sqrt{1-t_{n+1}}\gamma (n,j) &\text{for } q(n,j) = m_{n+1}\\
\sqrt{t_{n+1}}\gamma (n,j)&\text{for } q(n,j) = m_{n+1}+1\\
\gamma(n,j) &\text{for } q(n,j) > m_{n+1}+1
\end{cases}.\notag
\end{align*}
\item[(ii)] $0\le q(n+1,j)\le q(n,j)\le j$ and $0\le \gamma(n+1,j)\le \gamma(n,j)\le 1$ for every $n$ and $j$. 
\item[(iii)] $||P_n^\perp W^{(n)}e_j||= \gamma(n,j)$ for every $n$ and $j$.
\item[(iv)] For all $n>1$ and all $j$,
\[
W_{nj}= W^{(n))}_{nj}= \begin{cases}0\quad &\text{if } q(n-1,j) = 0\\
 \gamma(n-1,j)( V(m_{n},t_{n})\oplus I_\infty\big)_{1, q(n-1,j)} &\text{if } q(n-1,j)\ne0
\end{cases}. 
\]
In particular, all the entries of  $W$ are either $0$, $1$, or products of a finite number of the factors 
$\sqrt{t_k}$, $\sqrt{1-t_k}$ and $-\sqrt{1-t_k}$, but not more than one for each $k$.
\eL

The case $j=1$ is of special use. 

\bL{L:4.6} Given a sequence $\{m_k,t_k\}$ where $m_k\in \mathbb N$  and  $0< t_k\le1$, and the   co-isometry   \linebreak$W:=W(\{m_k,t_k\})$, let $g_n:=\gamma(n,1)^2$ and set  $g_\infty:=\lim g_n$.  Then
\item[(i)] $
g_\infty= \begin{cases}1& \quad   m_k> 1\text{ for all  }  k \\
\prod \{(1-t_k)\mid m_k=1\} & \quad\text{otherwise}
\end{cases}$

\item[(ii)]$
W_{n1}= \begin{cases}
\sqrt{t_ng_{n-1}}  &\text{if } m_n= 1\\
0\quad &\text{if } m_n >  1
\end{cases}$
\item[(iii)] 
$||We_1||^2= 1- g_\infty.$
\eL

 \bp
 \item[(i)] It is straightforward to solve the recurrence relation in Lemma \ref {L:4.5} (i) for $j=1$ and obtain 
\begin{align*}
q(n,1)&=\begin{cases}0\qquad\qquad\qquad \qquad\qquad& \text{if  } m_k=1, t_k=1 \,\text { for some } 1\le k \le n\\
1& \text{otherwise}
\end{cases},\\
 \gamma(n,1)&=\begin{cases}1 \qquad  \qquad&   \text{if  } m_k>1\,\text { for all } 1\le k \le n\\
\prod\{\sqrt{1-t_k}\mid m_k=1, 1\le k\le n\} &\text{otherwise}.
\end{cases}
\end{align*} 
Now (i) follows immediately.
\item[(ii)] By the proof of (i),  $q(n,1)\in \{0,1\}$ for all $n$ and thus

\begin{alignat*}{2}
W_{n1}&= \begin{cases}0&q(n-1,1)=0\\
\gamma(n-1,1) \big(V(m_{n},t_{n})\oplus I_\infty )_{1, 1}& q(n-1,1)=1
\end{cases} 
& & \qquad \text{(by  Lemma \ref {L:4.5}(iv))}\\
&= \begin{cases}0&q(n-1,1)=0\\
\sqrt{g_{n-1}}\big(V(m_{n},t_{n}) )_{1, 1}& q(n-1,1)=1\\
\end{cases}
\\
&= \sqrt{g_{n-1}}\begin{cases}\sqrt{t_n} & m_n=1\\
0&m_n>1
\end{cases}
& &  (q(n-1,1)=0\Leftrightarrow g_{n-1}=0).
\end{alignat*}
\item[(iii)]
Assume first that  the set $\{k\mid m_k=1\}$ is non-empty and order it into a strictly increasing, possibly finite, sequence $\{k_n\}_{1\le n \le N \le \infty}$. If $N=\infty$, then $g_\infty= \lim_n g_{k_n}$. If $N<\infty$, then 
 $g_k= g_{k_N}$ for every $k\ge k_N$, and hence $g_\infty=g_{k_N}$. Furthermore, for every $1\le n \le N$, $g_{ k_n-1}= g_{ k_{n-1}}$, where we set $n_0=0$ and $g_o=1$. Thus
$t_ {k_n}g_{ k_n-1}= g_{ k_{n-1}}-g_{ k_n}$. But then, from (ii) we have
\[
||We_1||^2= \sum_{n=1}^N t_ {k_n}g_{ k_n-1} = 
 \sum_{n=1}^N (g_{ k_{n-1}}-g_{ k_n})= g_{ k_0}-\lim g_{ k}
=1-g_\infty.
\]
Finally, if the set $\{k\mid m_k=1\}$ is empty,  then $g_k=1$ for all $k$ by (i) and hence $g_\infty=1$. By (ii), $W_{n1}=0$ for all $n$ and hence $||We_1||^2= 0$, also satisfying (iii).

 \ep

\bT{T:4.7} Let  $\xi \prec \eta $  for some $\xi, \eta \in \co*$ with $\xi_n >0$ for every $n$. Then the following conditions are equivalent.
\item[(i)] $\xi \preccurlyeq \eta $. 
\item[(ii)]  $\sum\{t_k \mid  m_k=1\}=\infty$.
\item[(iii)] $Q(\xi,\eta)$ is orthostochastic.

\eT
\bp
Notice that by Lemma \ref {L:2.4}  (see also Remark \ref{L:2.6} (ii)), $Q(\xi,\eta)$ is orthostochastic if and only if   $W(\xi,\eta)$ is unitary, in fact, orthogonal, since it has real entries.
\item [(iii)]  $\Rightarrow$ (ii)
In the case that there are infinitely many indices $k$ for which $m_k=1$ and $t_k=1$, then  \linebreak$\sum\{t_k \mid  m_k=1\}=\infty$ holds trivially, thus assume that there is an integer $N$ for which there are no $k>N$ with $m_k=1$ and $t_k=1$.  Then 
\begin{align*}
W(\xi, \eta)(W^{(N)}&)^*= w\text{-}\lim_n \Big( (I_{n-1}\oplus V(m_n,t_n)\oplus I_\infty)\cdots (I_N\oplus V(m_{N+1}, t_{N+1})\oplus I_\infty)\Big)\\
&=I_N\oplus   w\text{-}\lim_j \Big( (I_{j-N-1}\oplus V(m_j,t_j)\oplus I_\infty)\cdots ( V(m_{N+1}, t_{N+1})\oplus I_\infty)\Big)\\
&=I_N\oplus W(\{m_k,t_k\}_{k>}).
\end{align*}
By construction, $W(\{m_k,t_k\}_{k>})= W(\xi^{(N)}, \rho(N))$ is a co-isometry, however, since  $W^{(N)}$ and $W(\xi, \eta)$ are orthogonal,  the former by construction, the latter by hypothesis, it follows that $W(\{m_k,t_k\}_{k>})$ too is orthogonal. But then, by  Lemma \ref {L:4.6} applied to the majorization $\xi^{(N)}\prec \rho(N)$, we have 
\[
1= ||W(\xi^{(N)}, \rho(N))e_1||^2=  1-g_\infty=1- \prod \{ (1-t_k) \mid m_k=1, k > N \}.
\]
Thus $ \prod \{ (1-t_k) \mid m_k=1, k > N \} = 0$, hence  $\sum \{t_k \mid  m_k=1, k >N \}=\infty$ and therefore  \linebreak $\sum\{t_k \mid  m_k=1\}=\infty$.

\item [(ii)]  $\Rightarrow$ (iii)
\begin{alignat*}{4}
W(\xi, \eta)\,\text{is unitary} \, &\Leftrightarrow ||W(\xi, \eta)e_j||=1 \quad \text{for all } j&& \quad  (\text{as $W(\xi, \eta)$ is a co-isometry})\\
&  \Leftrightarrow ||P_nW(\xi, \eta)e_j||\to 1 \quad\text{for all } j \\
&  \Leftrightarrow ||P_nW^{(n)}e_j|| \to 1 \quad\text{for all } j & & \quad (\text{by Proposition \ref {P:3.1}})\\
&  \Leftrightarrow ||P_n^\perp W^{(n)}e_j|| \to 0 \quad\text{for all } j& & \quad(\text{as  $W^{(n)}$ is unitary })\\
&  \Leftrightarrow \gamma(n,j)\to 0 \quad\text{for all } j& & \quad(\text{by Lemma \ref {L:4.5}(iii) })
\end{alignat*}
For a fixed $j$, by Lemma \ref {L:4.5}(ii) and (i), the integer sequence $q(n,j)$ is monotone noincreasing in $n$ and it decreases by 1 for every $n$ for which $m_{n+1}=1$ and $q(n,j)>1$. Since there are infinitely many integers $k$ for which $m_k=1$, the sequence $q(n,j)$ must stabilize to either $0$ or $1$. If it is the former, since  $\gamma(n,j)=0$ whenever $q(n,j)=0$  we are done. If $q(n,j)=1$ for all $n\ge N$ for some $N\in \mathbb N$, then we obtain
 from the recurrence relation in Lemma \ref {L:4.5}(i)
 \[
\gamma(n,j) = \Big(\prod\{\sqrt{1-t_k} \mid m_k=1, n\ge k > N \}\Big) \gamma(N,j).
\]
Thus
\[
\lim_n(\gamma(n,j))^2= \Big(\prod\{(1-t_k) \mid m_k=1, k> N \}\Big) (\gamma(N,j))^2 = 0
\]
because  $\sum\{t_k \mid  m_k=1, k > N\}=\infty$.

\item [(ii)]  $\Rightarrow$ (i)
If there are infinitely many indices $k$ for which $m_k=t_k=1$, then $\xi\prec_b  \eta$ by Proposition \ref {P:4.4}  and hence $\xi\preccurlyeq  \eta$. If there are only finitely many $k$ for which $m_k=t_k=1$ and $K$ is the largest one, then by Lemma \ref{L:4.2}(ii) and (iii), $\sum_{j=1}^N(\eta_J-\xi_j)=0$ and hence $\xi^{(K)}\prec \eta^{(K)}$.  Thus it is sufficient (and necessary) to prove that $\xi^{(K)}\preccurlyeq \rho(K)$. Since by Lemma \ref {L:4.2} (iii), $\rho(K)= \eta^{(K)}$, and hence  $\{m_k,t_k\}_{k>K}$ is the sequence generated by  $\xi^{(K)}\prec \eta^{(K)}$,  we can assume without loss of generality that  $t_k<1$ whenever $m_k=1$. 
Order the indices $k$ for which $m_k=1$ into a strictly increasing sequence $\{k_n\}$ and set $q_n:= \prod_{j=1}^n(1-t_{k_n})$.  By the assumption that $0< t_{k_n}<1$,  it follows that $q_n>0$ for all $n$ and $q_n$ is strictly decreasing. The condition $\sum_{n=1}^\infty t_{k_n}=\infty$ guarantees that $q_n\to 0$. Since  $\delta_{k_n}= \eta_{k_n+1}+ \sum_{j=1}^{k_n}(\eta_j-\xi_j)$ by Lemma \ref{L:4.1}(iv), in order to show that $\xi\preccurlyeq \eta$ it is sufficient to prove that  $\lim_n \delta_{k_n}=0$.
For every $n>1$, 
\begin{alignat*}{2}
\delta_{k_n}&=  (1-t_{k_n})\rho(k_n-1)_1+t_{k_n}\rho(k_n-1)_2 & & \quad \text {(by definition, see (\ref {e:13}))}\\
&\le (1-t_{k_n})\rho(k_{n-1})_1+t_{k_n}\rho(k_{n-1})_2 & & \quad \text  {(by  (\ref{e:15}))}\\
&=  (1-t_{k_n})\delta_{k_{n-1}}+t_{k_n}\rho(k_{n-1})_2& & \quad \text  {(since  $m_{k_{n-1}}=1$, see (\ref{e:14}))}\\
&=  (1-t_{k_n})\delta_{k_{n-1}}+t_{k_n}\eta_{k_{n-1}+2}& & \quad \text  {(by Lemma \ref{L:4.1}(ii))}\\
&\le  (1-t_{k_n})\delta_{k_{n-1}}+t_{k_n}\eta_{k_{n-1}} & & \quad \text  {(by the monotonicity of $\eta$)}.
\end{alignat*}
Also $\delta_{k_1} =(1-t_{k_1})\rho(k_1-1)_1+t_{k_1}\rho(k_1-1)_2\le (1-t_{k_1})\eta_1+t_{k_1}\eta_2$. For convenience, set $k_o:=0$ and $\eta_{k_0}:=\eta_2$. Iterating, we obtain

\begin{align*}
 \delta_{k_n}&\le\big( \prod_{j=2}^n(1-t_{k_j})\big)\big((1-t_{k_1})\eta_1+t_{k_1}\eta_2\big)+ \sum_{j=3}^{n}\bigg( t_{k_{j-1}}\prod_{i=j}^n(1-t_{k_i})\eta_{k_{j-2}}\bigg)+ t_{k_n}\eta_{k_{n-1}}\\ 
&= \big(\prod_{i=1}^n(1-t_{k_i})\big)\eta_1 + \sum_{j=2}^{n}\bigg( t_{k_{j-1}}\prod_{i=j}^n(1-t_{k_i})\eta_{k_{j-2}}\bigg)+ t_{k_n}\eta_{k_{n-1}}\\
&= \big(\prod_{i=1}^n(1-t_{k_i})\big)\eta_1 + \sum_{j=2}^{n}\bigg(\Big(\prod_{i=j}^n(1-t_{k_i})-\prod_{i=j-1}^n(1-t_{k_i}) \Big)\eta_{k_{j-2}}\bigg)+ (1-(1-t_{k_n}))\eta_{k_{n-1}}\\
 &=q_n\eta_1+ \sum_{j=2}^{n}\Big(\frac{q_n}{q_{j-1}}- 
 \frac{q_n}{q_{j-2}} \Big)\eta_{k_{j-2}} + \big(1- \frac{q_n}{q_{n-1}}\big)\eta_{k_{n-1}}\\
 &=q_n\Bigg(\eta_1+  \sum_{j=2}^{n}\Big(\frac{1}{q_{j-1}}- 
 \frac{1}{q_{j-2}} \Big)\eta_{k_{j-2}}- \frac{1}{q_{n-1}}\eta_{k_{n-1}} \Bigg) +\eta_{k_{n-1}}\\
 &=q_n(\eta_1-\eta_2) + q_n \sum_{j=2}^{n}\frac{1}{q_{j-1} }(\eta_{k_{j-2}}-\eta_{k_{j-1}}) + \eta_{k_{n-1}},
 \end{align*}
where the last equality is obtained by ``summation by parts".  We know that $q_n\to0$ and clearly, $\eta_{k_{n-1}}\to 0$.
We claim that also $ q_n \sum_{j=2}^{n}\frac{1}{q_{j-1} }(\eta_{k_{j-2}}-\eta_{k_{j-1}})\to 0$.
 Indeed, for every $\epsilon >0$, choose $m$ for which $\eta_{k_{m}} < \epsilon$ and choose $N \ge m + 2$ so that for all $n \ge N$ 
\[
q_n \sum_{j=2}^{m+1}\frac{1}{q_{j-1} }(\eta_{k_{j-2}}-\eta_{k_{j-1}})< \epsilon.
 \]
 Then by the monotonicity of $q$ and $\eta$
\begin{align*}
q_n &\sum_{j=2}^{n}\frac{1}{q_{j-1}}(\eta_{k_{j-2}}-\eta_{k_{j-1}})
< \epsilon+
q_n \sum_{j=m+2}^{n}\frac{1}{q_{j-1} }(\eta_{k_{j-2}}-\eta_{k_{j-1}})\\
&\le  \epsilon + \sum_{j=m+2}^{n}(\eta_{k_{j-2}}-\eta_{k_{j-1}}) 
< \epsilon +\eta_{k_{m}} < 2\epsilon.
 \end{align*}
 This proves that $\lim_n \delta_{k_n} =0$ and hence that $\xi\preccurlyeq \eta$.\\
 \noindent We split the proof of the implication (i) $\Rightarrow$ (ii) or (iii) in two cases.\\
\noindent If $\eta$ has infinite support, then (i) $\Rightarrow$ (iii). Immediate by Remark \ref {R:3.10}.\\
\noindent If $\eta$ has finite support,  then (i)  $\Rightarrow$ (ii). Let  $\eta_N>0$ and $\eta_{N+1}=0$.  First notice that if $m_h=1$ for some $h\ge N-1$, then  $\rho(h)_2= \eta_{h+2}=0$ by Lemma \ref {L:4.1} (ii). For every  $k\ge h$, $\rho(k)_2\le \rho(h)_2$, hence $\rho(k)_2=0$ and by the definition of $m_k$ we have $m_k=1$. Thus  the sequence $\{m_k\}$ either eventually stabilizes at $1$ or is bounded away from $1$ from $N-1$ on.  We claim that the latter case is impossible. Reasoning by contradiction, assume that $m_k\ge 2$  for all $k\ge N-1$. Then for every $k\ge N-1$ and every $n\ge m_k$ we have
\begin{alignat*}{2}
\xi_N&\le \rho(N-1)_1& & \quad\text{(since $\xi^{(N-1)}\prec \rho(N-1)$ by (\ref {e:15}))}\\
&=\rho(k)_1	& & \quad\text{(by (\ref {e:14})}\\
&\le \sum_{j=1}^n \rho(k)_j\\
&= \sum_{j=1}^{n+k}\eta_j - \sum_{j=1}^k \xi_j & & \quad\text{(by Lemma \ref {L:4.1}(iii))}\\
&=\sum_{j=k+1}^\infty \xi_j \to 0 & & \quad\text{(since $\sum_{j=1}^N\eta_j=\sum_{j=1}^\infty \xi_j$)}
\end{alignat*}
which contradicts the assumption that $\xi$ has infinite support and hence $\xi_N> 0$.
Therefore there is a $K\ge N$ such that $m_k=1$ for all $k \ge K$.  But then, if $k> K$,   $\rho(k-1)= < \delta_{k-1}, 0, \ldots>$  and hence $\xi_k= t_k\delta_{k-1}$, where $\delta_{k-1} = \sum_{j=1}^k\eta_j-\sum_{j=1}^{k-1}\xi_j = \sum_{j=k}^{\infty}\xi_j$ by Lemma \ref{L:4.1}(iv).  Therefore 
$1-t_k= \frac{\sum_{j=k+1}^{\infty}\xi_j} {\sum_{j=k}^{\infty}\xi_j}  \in (0,1)$
and hence for all $M\ge K+1$
\[
\prod_{k={K+1}}^M(1-t_k) = \frac{\sum_{j=M+1}^{\infty}\xi_j} {\sum_{j=K+1}^{\infty}\xi_j}\to 0 \quad \text{for } M\to \infty.
\]
As a consequence,  $\sum_{k=K+1}^\infty t_k = \infty$ and hence $\sum\{t_k\mid m_k =1\}=\infty$.\
\ep
\bR{R:4.8} Given a sequence $\{m_k,t_k\}$ where $m_k\in \mathbb N$ and $0< t_k\le 1$,  Proposition \ref {P:3.1} constructs the co-isometry $W(\{m_k,t_k\})$ and Lemmas  \ref {L:4.5}, and \ref {L:4.6} provide further properties for that construction. It is easy to see that, if in lieu of $Q(\xi, \eta)$ we consider the Schur-square of $W(\{m_k,t_k\})$,  the implications (ii) $\Leftrightarrow$ (iii) in Proposition \ref {P:4.4} and Theorem \ref {T:4.7} still hold for this more general setting. 
\eR

\section{An extension of the Horn Theorem to nonsummable sequences} \label{S: 5} 
In Theorem \ref {T:3.9} we proved that if $\xi\preccurlyeq \eta$, then $\xi=Q\eta$ for some orthostochastic matrix $Q$. While strong majorization is necessary and sufficient in the summable case by Lemma \ref {L:2.10}, in the nonsummable case it is not, as seen in Example \ref {E:2.11}. In fact, we are going to prove that the condition $\xi\prec\eta$ will always suffice when $\xi$ is nonsummable. 

Our strategy will be to decompose any pair of sequences $\xi,\eta\in \co*$ with $\xi \prec \eta$ and $\xi$ nonsummable into  ``direct sums" of pairs of sequences  
$\xi(r) \preccurlyeq \eta(r)$ ($r=1,2, \ldots $). 
The key step in this process is the following ``shift" lemma.

\bL{L:5.1}
Let $\xi,\eta\in \co*\setminus (\ell^1)^*$ and assume that  $\xi\prec \eta $ but $\xi \not\preccurlyeq \eta$. Then there are  integers $p$ and $n$ with $0\le p<n$,   for which $\xi\chi[1,n]\prec \eta\chi[1,n-p]$ and $\xi^{(n)}\prec \eta^{(n-p)}.$
\eL

\bp
By hypothesis, $\alpha:={\varliminf} \sum_{j=1}^n(\eta_j-\xi_j)>0$. 
If there is some $n\in \mathbb N$ for which   \linebreak $\sum_{j=1}^n(\eta_j-\xi_j) \le  \sum_{j=1}^m(\eta_j-\xi_j)$ for every $m\ge n$, which is certainly the case if $\alpha = \infty$,  then $\xi^{(n)}\prec \eta^{(n)}$ and hence the pair $p=0$ and $n$ satisfies the requirement. Assume therefore that there is no such $n$ and hence that  $\sum_{j=1}^n(\eta_j-\xi_j) >\alpha $ for every $n$. 
In particular, $\eta_1>\xi_1 +\alpha> \alpha$. 
Let $N_1$ be an integer for which $\eta_{N_1}< \alpha$ and for every $n\ge N_1$, let $ p(n)$  be the largest integer  in $[1, n)$ for which
\be{e:20}
\sum_{j=n-p(n)+1}^n\eta_j \le \alpha.
\ee
 By the monotonicity of $\eta, $ $
\sum_{j=n-p(n)+2}^{n+1}\eta_j \le \sum_{j=n-p(n)+1}^n\eta_j \le \alpha$ 
and hence by the maximality of $p(n+1)$, it follows that $ p(n+1) \ge p(n) $ for every $n\ge N_1$, i.e., the sequence $p(n)$ is monotone nondecreasing. 
Since $\eta_j\to 0$, it follows that $p(n) \to \infty$,  and since $\eta$ is nonsummable, it follows that $n-p(n) \to \infty$.  
Then $\eta_{n-p(n)} < \frac{\alpha}{2} $ for all $n\ge N_2$ for some $N_2 \ge N_1$ and hence by the maximality of $p(n)$,
\be {e:21}
\sum _{j=n-p(n)+1}^{n}\eta_j > \frac{\alpha}{2}\quad\text {for every }\,  n\ge N_2. 
\ee
Now $
 \sum_{j=1}^{n-p(n)} \eta_j - \sum_{j=1}^{n}\xi_j = \sum_{j=1}^n (\eta_j-\xi_j) -  \sum _{j=n-p(n)+1}^{n}\eta_j > 0$
from which it follows that  \linebreak $\xi\chi[1,n]\prec \eta\chi[1,n-p(n)]$ for every $n\ge N_2$. 
It remains to prove  that there is an $n\ge N_2$ for which  $\xi^{(n)}\prec \eta^{(n-p(n))}$. 
Reasoning by contradiction, assume that for every  $n \ge N_2$ there is an integer $q(n)> n$ for which
$\sum_{j=1}^{q(n)-n}(\eta^{(n-p(n))}_j-\xi^{(n)}_j) < 0$, i.e.,  $\sum_{j=n+1}^{q(n)}\xi_j >  \sum_{j=n-p(n) +1}^{q(n)-p(n)}\eta_j$. Then

\be{e:22}
\sum_{j=n+1}^{q(n)}(\xi_j-\eta_j) >   \sum_{j=n-p(n) +1}^{q(n)-p(n)}\eta_j- \sum_{j=n+1}^{q(n)}\eta_j
=\sum_{j=n-p(n)+1}^n\eta_j-\sum_{j=q(n)-p(n)+1}^{q(n)}\eta_j \ge 0 ,
\ee
for every $n\ge N_2$, where the last inequality follows form the monotonicity of $\eta$.
From this inequality and (\ref{e:21}) 
\[
\frac{\sum_{j=q(n)-p(n)+1}^{q(n)}\eta_j }{ \sum_{j=n-p(n)+1}^n\eta_j}
> 1- \frac{\sum_{j=n+1}^{q(n)}(\xi_j-\eta_j)}{ \sum_{j=n-p(n)+1}^n\eta_j} 
> 1- \frac{2}{\alpha} \sum_{j=n+1}^{q(n)}(\xi_j-\eta_j).
\]
Set $m_1=N_2$ and $m_{k+1}:= q(m_{k})$. The sequence $m_k$ is strictly increasing and for every $k\ge 1$,
\[
\frac{\sum_{j=m_{k+1}-p(m_k)+1}^{m_{k+1}}\eta_j } { \sum_{j=m_k-p(m_k)+1}^{m_k} \eta_j}
>1- \frac{2}{\alpha} \sum_{j=m_k+1}^{m_{k+1}} (\xi_j-\eta_j).
\]
Given that $\eta$ is nonincreasing and that $p(m_k)$ is nondecreasing, the average of $\eta$ over the integer interval  \linebreak $\{m_{k+1}-p(m_{k+1})\le j\le m_{k+1}\}$ must be larger or equal than its average over the integer  interval  \linebreak $\{m_{k+1}-p(m_{k})\le j\le m_{k+1}\}$  and hence 
\be{e:23}
\frac{\frac{1}{p(m_{k+1})}\sum_{j=m_{k+1}-p(m_{k+1})+1}^{m_{k+1}}\eta_j}{\frac{1}{p(m_k)} \sum_{j=m_k-p(m_k)+1}^{m_k}\eta_j}
\ge \frac{\frac{1}{p(m_k)}\sum_{j=m_{k+1}-p(m_k)+1}^{m_{k+1}}\eta_j}{\frac{1}{p(m_k)} \sum_{j=m_k-p(m_k)+1}^{m_k}\eta_j}
> 1- \frac{2}{\alpha} \sum_{j=m_k+1}^{m_{k+1}}(\xi_j-\eta_j).
\ee
Now by (\ref {e:22}), $\sum_{j=m_k+1}^{m_{k+1}}(\xi_j-\eta_j) > 0$ for every $k\ge1$ and by assumption, $\sum_{j=1}^{m_k}(\eta_j-\xi_j)>\alpha>0$. Thus for every $h>1$,
\[
\sum_{j=1}^{m_h}(\eta_j-\xi_j)= \sum_{j=1}^{m_1}(\eta_j-\xi_j) - \sum_{j=m_1+1}^{m_{h}}(\xi_j-\eta_j) 
=  \sum_{j=1}^{m_1}(\eta_j-\xi_j) -\sum_{k=1}^{h-1} \sum_{j=m_k+1}^{m_{k+1}}(\xi_j-\eta_j) 
>0,
\]
whence $\sum_{k=1}^\infty \sum_{j=m_k+1}^{m_{k+1}}(\xi_j-\eta_j) < \infty.$
 Choose $k_o\in \mathbb N$ for which  $ \sum_{k=k_o}^\infty \sum_{j=m_k+1}^{m_{k+1}}(\xi_j-\eta_j) < \frac{\alpha}{2}.$
In particular, for all $k\ge k_o$ we have $0< \frac{2}{\alpha} \sum_{j=m_k+1}^{m_{k+1}}(\xi_j-\eta_j)< 1$ and hence from  equation (\ref {e:23}) we have for every $K\ge k_o$
\[
0 < \prod_{k=k_o}^K\bigg(1- \frac{2}{\alpha} \sum_{j=m_k+1}^{m_{k+1}}(\xi_j-\eta_j) \bigg)
<
\frac {\frac{1}{p(m_{K+1)}} \sum_ { j=m_{K+1}-p(m_{K+1})+1 }^{ m_{K+1} } \eta_j }  { \frac{ 1 }{ p(m_{k_o}) } \sum_{j=m_{k_o}-p(m_{k_o})+1}^{m_{k_o}}\eta_j}
\le 
2\frac { p(m_{k_o}) } {p( m_{K+1} ) },
\]
where the last inequality follows from  the inequalites (\ref {e:20}) and (\ref  {e:21}).
Now, on the one hand,  $p(m_{k}) \to \infty$ and hence $2\frac { p(m_{k_o}) } {p( m_{K+1} ) } \to 0$ for $K\to\infty$. On the other hand, the sequence $  \frac{2}{\alpha}\sum_{j=m_k+1}^{m_{k+1}}(\xi_j-\eta_j) \in (0,1)$ and  is summable, hence $\prod_{k=k_o}^\infty\Bigg(1- \frac{2}{\alpha}\sum_{j=m_k+1}^{m_{k+1}}(\xi_j-\eta_j) \Bigg) > 0,$ 
 a contradiction.
 
\ep

\bL{L:5.2} Let $\xi,\eta\in \co*\setminus (\ell^1)^*$ and assume that  $\xi\prec \eta $ but $\xi \not\preccurlyeq \eta$. 
Then there there are two partitions of $\mathbb N$ into  sequences, $\mathbb N = \{n_j^{(1)}\} \, \dot \cup \, \{n_j^{(2)}\}$ and  $\mathbb N = \{m_j^{(1)}\} \, \dot{ \cup }\, \{m_j^{(2)}\}$  with $n_1^{(1)}= m_1^{(1)}=1$ for which, 
if $\xi':=\{\xi_{n_j^{(1)}}\}$, $\eta':=\{\eta_{m_j^{(1)}}\}$, $\xi'':=\{\xi_{n_j^{(2)}}\}$, and
$\eta'':=\{\eta_{m_j^{(2)}}\}$ are the corresponding subsequences of $\xi$ and $\eta$, 
then $\xi'\preccurlyeq \eta'$, $\xi''\prec \eta''$,  $\xi''\not\preccurlyeq \eta''$, and $\xi'\in  (\ell^1)^*$.
\eL
\bp 
By Lemma \ref {L:5.1}, $\xi\chi[1,N]\prec \eta\chi[1,N-p]$ and $\xi^{(N)}\prec \eta^{(N-p)}$ for some pair of  integers $p$ and $N$ with $0\le p < N$. Let 
\[
\alpha:= {\varliminf}\Big (\sum_{j=1}^k (\eta_j- \xi_j)\Big), \quad \beta:=\sum_{j=1}^{N-p}\eta_j - \sum_{j=1}^N \xi_j, \quad \text{and }\quad \gamma:= {\varliminf} \Big (\sum_{j=1}^k (\eta^{(N-p)}_j- \xi^{(N)}_j)\Big). 
\]
By hypothesis, $\alpha > 0$ and $\beta\ge 0, \gamma \ge 0$. Since for $k > N-p$
\[
\sum_{j=1}^k(\eta_j-\xi_j)= \beta + \sum_{j=1}^{k-N+p}(\eta^{(N-p)}_j-\xi^{(N)}_j) + \sum_{j=k+1}^{k+p}\xi_j
\]
and $\sum_{j=k+1}^{k+p}\xi_j\to 0$ for $k\to\infty$, it follows that  $0<\alpha = \beta+ \gamma$, so $\beta$ and $\gamma$ cannot both vanish.

Assume first that $\beta> 0$.  The strategy for the construction of the sequences $\xi'$, $\xi''$, $\eta'$, and $\eta''$ is to first move a finite number of entries from the infinite sequence $\xi^{(N)}$ to the finite sequence $\xi\chi[1,N]$, i.e., delete them from the first sequence and insert them after the last nonzero term of the second one,  and do so while while controlling the sum and preserving the majorization by $\eta\chi[1,N-p]$ of the new finite sequence. This will automatically preserve majorization of the new infinite sequence by $\eta^{(N-p)}$. At the next step, move a single entry from the sequence $\eta^{(N-p)}$ to the sequence $ \eta\chi[1,N-p]$, so to preserve majorization of the two infinite sequences and still control the sums, while majorization of the two finite ones is automatically preserved. And then  iterate the process.
 
Now we make this strategy precise. We construct three strictly increasing sequences of integers $k_j$, $h_i$ and $q_i$ with $N< k_{q_i}< h_i \le h_i +p < k_{q_i+1} < k_{q_i+2} < \dots < k_{q_{i+1}}$ so that
\begin{align}
&\beta + \sum_{j=1}^{i-1}\eta_{h_j}-\frac{1}{i} <  \sum_{j=1}^{q_i} \xi_{k_j} < \beta +  \sum_{j=1}^{i-1}\eta_{h_j},\label{e:24}\\
&\delta_i:=\sum_{j=1}^{q_i}\xi_{k_j} - \sum_{j=k_{q_i}+1}^{k_{q_i}+q_i}\xi_j   >\delta_{i-1},\label{e:25}\\
&\eta_{h_i}< \min\{\frac{1}{2^i}, \delta_i-\sum_{j=1}^{i-1}\eta_{h_j}\},\label{e:26}
\end{align}
where for $i=1$ we take $0$ in place of $\sum_{j=1}^{i-1}\eta_{h_j}$ and of  $\delta_{i-1}$.

To start the construction, use the fact that $\xi_j\to 0$ and is nonsummable to choose $N < k_1 < \dots < k_{q_1-1}$ for which $\beta-1< \sum_{j=1}^{q_1-1} \xi_{k_j} < \beta$. Since   $\xi$ has infinite support,  it has an infinite  subsequence for which $\xi_{p_n}> \xi_{p_n+1}$. Choose  $k_{q_1}\in \{p_n\}$ large enough so that $ \sum_{j=1}^{q_1} \xi_{k_j} < \beta$. By the monotonicity of $\xi$, it follows that $\delta_1:= \sum_{i=1}^{q_1}\xi_{k_i}-\sum_{j=k_{q_1}+1}^{k_{q_1}+q_1}\xi_j    >0 $ and conditions   (\ref {e:24}) and (\ref {e:25}) are thus satisfied for $i=1$. To satisfy  also (\ref {e:26}) it is enough to choose $h_1> k_{q_1}$ so that $\eta_{h_1}< \min\{\frac{1}{2}, \delta_1\} $, which is always possible since $\eta_j\to 0$ and $\delta_1>0$.   Assume now the construction of the three integer sequences up to some $i-1$ and choose  $h_{i-1}+p < k_{q_{i-1}+1} <  k_{q_{i-1}+2}< \dots < k_{q_i-1}$ for which
\[
 \beta+ \sum_{j=1}^{i-1}\eta_{h_j} - \sum_{j=1}^{q_{i-1}} \xi_{k_j}- \frac{1}{i} 
<\sum_{j=q_{i-1}+1}^ {q_{i}-1}\xi_{k_j}
<  \beta + \sum_{j=1}^{i-1}\eta_{h_j} - \sum_{j=1}^{q_{i-1}} \xi_{k_j}.
\]
 Choose  $k_{q_i}\in \{p_n\}$, $k_{q_i}> k_{q_i-1}$ large enough so that $\sum_{j=q_{i-1}+1}^ {q_{i}}\xi_{k_j}
<  \beta + \sum_{j=1}^{i-1}\eta_{h_j} - \sum_{j=1}^{q_{i-1}} \xi_{k_j}$, i.e., so to satisfy (\ref {e:24}).
Now
\begin{alignat*}{2}
\delta_i-\delta_{i-1} &= \sum_{ j=q_{i-1}+1 } ^{ q_i } \xi_{k_j}  - \sum_{ j=k_{q_i}+1}^{ k_{q_i}+q_i }\xi_j 
 + \sum_{j=k_{ q_{i-1}}+1 }^{ k_{ q_{i-1} }+q_{i-1} }\xi_j & &\\
&=\sum_{j=q_{i-1} +1}^ {q_i} \xi_{k_j} - \sum_{j=k_{q_i}+q_{i-1}+1}^{ k_{q_i}+q_i }\xi_j 
+ \sum_{ j=k_{q_{i-1} }+1 }^ {k_{ q_{i-1} }+q_{i-1} }\xi_j -  \sum_{ j=k_{q_i}+1 }^{ k_{ q_i }+q_{i-1} }\xi_j \\
&\ge \sum_{j=q_{i-1} +1}^ {q_i} \xi_{k_j} - \sum_{j=k_{q_i}+q_{i-1}+1}^{ k_{q_i}+q_i }\xi_j & & \text{(by the monotonicity of $\xi$)}\\
&> 0.  & &\text{(because $ \xi_{k_{q_i} } > \xi_{k_{q_i}+1 }$)}
\end{alignat*}
Thus (\ref{e:25}) is satisfied. By the induction assumption that $\eta_{h_n}< \min\{\frac{1}{2^n}, \delta_n-\sum_{j=1}^{n-1}\eta_{h_j}\}$ for all $1\le n\le i-1$, we see that $\delta_{i} >  \delta_{i-1} > \sum_{j=1}^{i-i}\eta_{h_j}$ and since $\eta_n \to 0$ we can choose $h_i > k_{q_i}$ so to satisfy also  (\ref{e:26}).  

Now define 
\be{e:27}
n^{(1)}:= < 1, \ldots, N, k_1, k_2, \ldots>\quad \text {and} \quad
m^{(1)}:= < 1, \ldots, N-p, h_1, h_2, \ldots>
\ee
and $n^{(2)}$, $m^{(2)}$ are the complementary sequences of $n^{(1)}$, $m^{(1)}$ respectively. Explicitly, 
\begin{alignat}{4}\label{e:28}
\xi':&=\, <\xi_1, \ldots, \xi_N, \xi_{k_1}, \xi_{k_2}, \ldots> & &\quad\text{and }
\eta':&&=\, <\eta_1, \ldots, \eta_{N-p}, \eta_{h_1}, \eta_{h_2}, \ldots>\\
\xi'':&= \,< \xi_{N+1}, \xi_{N+2},\ldots,\xi_{k_1-1}, \xi_{k_1+1},\ldots >& &\quad\text{and }
\eta'': &&=\, < \eta_{N-p+1}, \ldots,\eta_{h_1-1}, \eta_{h_1+1},\ldots >.\notag
\end{alignat}
First we verify  that  $\xi'\preccurlyeq \eta'$.\\
If $m\le N-p$, then $\sum_{j=1}^m(\eta'_j-\xi'_j)=\sum_{j=1}^m(\eta_j-\xi_j)\ge 0$. \\
If $N-p < m \le N$, then $\sum_{j=1}^m(\eta'_j-\xi'_j)\ge  \sum_{j=1}^{N-p}\eta_j- \sum_{j=1}^N\xi_j = \beta >0.$\\
Finally, if  $m>N$, let $q_{i-1}< m-N\le q_i$, where we set $q_o=0$ for convenience. Then $m> N+q_{i-1}\ge N+i-1\ge N-p +i-1$ and hence
 \begin{alignat*}{2}
 \sum_{j=1}^m\eta'_j &\ge  \sum_{j=1}^{ N-p +i-1}\eta'_j 
 =  \sum_{j=1}^{N-p}\eta_j +\sum_{j=1}^{i-1}\eta_{h_j} &&\\
 &\ge \sum_{j=1}^{N-p}\eta_j-\beta+\sum_{j=1}^{q_i}\xi_{k_j} &&\quad \text{(by (\ref{e:26}))}\\
&=  \sum_{j=1}^N\xi_j+\sum_{j=1}^{q_i}\xi_{k_j}
\ge  \sum_{j=1}^N\xi_j+\sum_{j=1}^{m-}\xi_{k_j} = \sum_{j=1}^m\xi'_j.
 \end{alignat*}
Thus $\xi' \prec \eta'$. For every $i>1$, 
 \begin{alignat*}{2}
 \sum_{j=1}^{N-p+i-1}\eta'_j&= \sum_{j=1}^{N-p}\eta_j + \sum_{j=1}^{i-1}\eta_{h_j} 
=  \sum_{j=1}^N\xi_j+\beta + \sum_{j=1}^{i-1}\eta_{h_j}&&\quad \text{(by the definition of $\beta$)}\\
&\le \sum_{j=1}^N\xi_j+ \sum_{j=1}^{q_i} \xi_{k_j}+ \frac{1}{i}  &&\quad \text{(by (\ref{e:26}))}\\
&= \sum_{j=1}^{N+q_i}\xi'_j+ \frac{1}{i}  < \sum_{j=1}^\infty \xi'_j + \frac{1}{i}  .
 \end{alignat*}
Therefore $\sum_{j=1}^\infty \eta'_j\le \sum_{j=1}^\infty \xi'_j$ and since $\xi' \prec \eta'$ and  by (\ref{e:26}), $\eta'\in(\ell^1)^*$ and hence  $\xi'\in(\ell^1)^*$, 
equality follows, i.e., $\xi'\preccurlyeq \eta'$.

Next, we verify that $\xi''\prec \eta'' $.  We start with the following two inequalities. \\
If $h_i+p\le N+m< h_{i+1}+p$, then
\be{e:29}
\sum_{j=1}^m \eta''_j = \sum_{j=1}^{m+i} \eta^{(N-p)}_j- \sum_{j=1}^i\eta_{h_j}
\ge \sum_{j=1}^{m} \eta^{(N-p)}_j - \sum_{j=1}^i\eta_{h_j}
\ee
If $k_{q_{i}}\le N +m< k_{q_{i+1}}$ and $N+m \ge h_1+p$, then 
\begin{align}\label{e:30}
\sum_{j=1}^m \xi''_j&=   \sum_{j=1}^{m+q_{i}} \xi^{(N)}_j -  \sum_{j=1}^{q_{i}}\xi_{k_j}
= \sum_{j=1}^{m} \xi^{(N)}_j +  \sum_{j=N+m+1}^{N+m+q_{i}} \xi_j -  \sum_{j=1}^{q_{i}}\xi_{k_j}\\
&\le  \sum_{j=1}^{m} \xi^{(N)}_j +  \sum_{j=k_{q_{i}}+1}^{k_{q_{i}}+q_{i}} \xi_j -  \sum_{j=1}^{q_{i}}\xi_{k_j}
= \sum_{j=1}^{m} \xi^{(N)}_j-\delta_i,\notag
\end{align}
where the inequality follows from the monotonicity of $\xi$.\\
Since $
N < h_1+p < k_{q_2}< \dots k_{q_i}< h_i \le h_i+p < k_{q_{i+1}}< \dots,$
to prove that $\sum_{j=1}^m (\eta''_j-\xi''_j) \ge 0$, we need to consider three cases:  $N+m < h_1 + p$, $h_i+p\le N+m< k_{q_{i+1}}$ for some $i\ge 1$, and  $k_{q_i} \le N+m < h_i+p$ for some $i\ge 2$.\\
In the first case, since $\xi^{(N)}\prec \eta^{(N-p)}$,
\[
\sum_{j=1}^m (\eta''_j-\xi''_j) = \sum_{j=1}^m (\eta^{(N-p)}_j -\xi''_j)  \ge  \sum_{j=1}^{m} (\eta^{(N-p)}_j - \xi^{(N)}_j)\ge0.
\]
In the second case $h_i+p\le N+m< h_{i+1}+p$ and $k_{q_{i}}\le N +m< k_{q_{i+1}}$ while $N+m \ge h_1+p$, hence combining (\ref {e:29}) and (\ref{e:30}) yields
\[
\sum_{j=1}^m (\eta''_j-\xi''_j) \ge \sum_{j=1}^m (\eta^{(N-p)}_j -\xi''_j)+ \delta_i- \sum_{j=1}^{i}\eta_{h_j}>0
\]
where the last inequality follows from $\xi^{(N)}\prec \eta^{(N-p)}$ and (\ref {e:26}).\\
In the third case  $h_{i-1}+p\le N+m < h_i+p$ and $k_{q_{i}}\le N +m< k_{q_{i+1}}$ while $N+m \ge k_{q_2}>h_1+p$, hence combining (\ref {e:29}) and (\ref{e:30}) yields
\[
\sum_{j=1}^m (\eta''_j-\xi''_j) \ge \sum_{j=1}^m (\eta^{(N-p)}_j -\xi''_j)+ \delta_i- \sum_{j=1}^{i-1}\eta_{h_j}>0.
\]
This proves that $\xi''\prec \eta'' $.
Finally, since $\sum_{j=1}^\infty \eta'_j=\sum_{j=1}^\infty  \xi'_j\infty$, it follows that 
\[
\varliminf \sum_{j=1}^k (\eta''_j- \xi''_j)  = \varliminf \sum_{j=1}^k (\eta_j- \xi_j)  = \alpha >0,
\]
  i.e.,  $\xi''\not\preccurlyeq \eta''$. This completes the proof of the lemma for the case $\beta >0$. 

The proof for the case where $\gamma >0$ is similar, but simpler. To further slightly simplify the proof, assume without loss of generality that $\beta=0$, i.e., $ \sum_{j=1}^N\xi_j =\sum_{j=1}^{N-p}\eta_j$.
Choose an $M\in \mathbb N$ and $\gamma_o>0$ so that $\sum_{j=1}^m(\eta^{(N-p)}_j - \xi^{(N)}_j) >\gamma_o$ for all $m \ge M$ and select a strictly increasing infinite sequence of integers $\{h_j\}$ with $h_1\ge N-p+ M $ for which $\sum_{j=1}^\infty\eta_{h_j} \le \gamma_o$. Then choose  strictly increasing sequences of integers $\{k_j\}$ and $\{q_i\}$ so that for all $i\ge 1$
\be{e:31}
 \sum_{j=1}^{i}\eta_{h_j}-\frac{1}{i} <  \sum_{j=1}^{q_i} \xi_{k_j} <   \sum_{j=1}^{i}\eta_{h_j}
\ee
which can be achieved as in the proof of the case $\beta>0$, by using the fact that $\xi_i\to0$ and $\xi$ is nonsummable. Define the sequences $n^{(1)}$,  $m^{(1)}$ and their complements $n^{(2)}$,  $m^{(2)}$ as in (\ref {e:27}) and hence $\xi'$, $\eta'$, $\xi''$, and $\eta''$  as in (\ref{e:28}), i.e., by ``moving" the entries $\xi_{k_j}$ (resp. $\eta_{k_j}$) from $\xi^{(N)}$ to $\xi\chi[1, N]$ (resp. from $\eta^{(N-p)}$ to $\eta\chi[1, N-p]$.)

First we show that $\xi'\preccurlyeq \eta'$.  If $1\le n\le N$, then from $\xi\chi[1, N]\prec \eta\chi[1, N-p])$ we have
\[
\sum_{j=1}^n \xi'_j = \sum_{j=1}^n \xi_j \le \sum_{j=1}^n (\eta\chi[1, N-p])_j \le \sum_{j=1}^n \eta'_j.
\]
Set $q_o=0$. If  $q_{i-1}<m\le q_i$ for some $i\ge1$, then $m \ge q_{i-1}+1\ge i$ and by (\ref {e:31}),
\[
\sum_{j=1}^{N+m}\xi'_j
\le \sum_{j=1}^{N+q_i}\xi'_j
=\sum_{j=1}^N\xi_j + \sum_{j=1}^{q_i}\xi_{k_j} 
\le \sum_{j=1}^{N-p}\eta_j+ \sum_{j=1}^{i}\eta_{h_j}
= \sum_{j=1}^{N+i}\eta'_j
\le \sum_{j=1}^{N+m}\eta'_j,
\]
which proves that  $\xi'\prec \eta'$. By construction, $\eta'$ is summable and hence $\sum_{j=1}^\infty \xi' \le \sum_{j=1}^\infty \eta'$. On the other hand, from (\ref {e:31}), for every $i$, 
\[
\sum_{j=1}^{N-p+i}\eta'_j
= \sum_{j=1}^{N-p}\eta_j+ \sum_{j=1}^{i}\eta_{h_j}
<  \frac{1}{i} +\sum_{j=1}^N\xi_j + \sum_{j=1}^{q_i}\xi_{k_j}
=  \frac{1}{i} + \sum_{j=1}^{N+q_i}\xi'_j
< \frac{1}{i} + \sum_{j=1}^{\infty}\xi'_j.
\]
Thus $\sum_{j=1}^{\infty}\eta'_j= \sum_{j=1}^{\infty}\xi'_j$ and hence $\xi'\preccurlyeq \eta'$.

Next we prove that $\xi^{(N)}\prec\eta''$, and hence, since $\xi'' \le \xi^{(N)}$, that  $\xi''  \prec\eta''$. Indeed,   if $N-p +m< h_1$, then $\sum_{j=1}^m \eta''_j = \sum_{j=1}^{m} \eta^{(N-p)}_j \ge \sum_{j=1}^{m} \xi^{(N)}_j $. If  $h_i\le N-p +m< h_{i+1}$ for some $i\ge1$, then, as in (\ref{e:29}),
$\sum_{j=1}^m \eta''_j \ge \sum_{j=1}^{m} \eta^{(N-p)}_j - \sum_{j=1}^i\eta_{h_j}$. But $\sum_{j=1}^i\eta_{h_j} < \gamma_o$ and since $m > h_1-N+p \ge M$ it follows also that $\sum_{j=1}^m(\eta^{(N-p)}_j - \xi^{(N)}_j) >\gamma_o$. Thus we conclude that $\sum_{j=1}^m \eta''_j>  \sum_{j=1}^{m} \xi^{(N)}_j .$
Finally, as in the case of $\beta >0$ it is now immediate to see that  $\xi''\not \preccurlyeq \eta''$.
\ep

Having thus prepared the groundwork, we can present  an infinite-dimensional extension of the Horn theorem \cite[Theorem 4]{Horn} (see Theorem \ref {T:1.1}(ii)) for nonsummable sequences.

\bT{T:5.3}
If $\xi, \eta\in \co*$ and $\xi\not \in (\ell^1)^*$, then $\xi\prec\eta$ if and only if $\xi=Q\eta$ for some orthostochastic matrix $Q$. 
\eT
\bp
The sufficiency is well-known since Markus (see \eqref{e:1} or the proof of Lemma \ref{L:2.8} ).  
For the necessity, if $\xi\preccurlyeq \eta$, then the result follows  from Theorem \ref{T:3.9}.  
If $\xi\not \preccurlyeq \eta$, then  Lemma \ref {L:5.2} applied iteratively partitions $\mathbb N$ into the union of infinitely many sequences $\{n_k^{(r)}\}$, and sequences $\{m_k^{(r)}\}$, i.e.,  \linebreak
$ \mathbb N = \dot {\bigcup}_{k=1}^\infty \{n_k^{(r)}\} = \dot {\bigcup}_{k=1}^\infty \{m_k^{(r)}\} $,
so that $\{\xi_{n_k^{(r)}} \} \preccurlyeq \{\eta_{m_k^{(r)}} \}$, the exhaustion of $\mathbb N$ being guaranteed by the condition in Lemma \ref {L:5.2}  that  \mbox{$n_1^{(1)}=m_1^{(1)}=1$.} Then by Theorem  \ref{T:3.9} we can find for each $r$ an orthogonal matrix  $W(r)$ with Schur-square $Q(r)$, for which $\{\xi_{n_k^{(r)}}\}=Q(r) \{\eta_{m_k^{(r)}} \}$. ``Direct summing" the matrices $Q(r)$ with respect to the two decompositions of the basis yields the required orthostochastic matrix $Q$. Explicitly, let $R_r, S_r\in B(H)$ be the diagonal projections on the subspaces with bases $\{e_{m_k^{(r)}}\}$ and $\{e_{m_k^{(r)}}\}$ respectively. Then $\sum_{r=1}^\infty R_r=\sum_{r=1}^\infty S_r=I$. We can identify the matrices $Q(r)$ and $W(r)$ with  corresponding operators in $B(R_rH, S_rH)$ and then define $Q:\sum_{r=1}^\infty S_rQ(r)R_ r$ and  $W:\sum_{r=1}^\infty S_rW(r)R_ r$. It is clear that $W, Q\in B(H)$, $W$ is unitary and its matrix is
\[
W_{ij}:=\begin{cases} (W(r))_{hk} &\quad \text{if }i = n^{(r)}_h, \, j = m^{(r)}_k \,\,\text{for some $r$, $h$, $k$}\\
0& \quad \text{otherwise}
\end{cases}.
\]
Similarly, 
\[
Q_{ij}:=\begin{cases} (Q(r))_{hk} &\quad \text{if }i = n^{(r)}_h, \, j = m^{(r)}_k \,\,\text{for some $r$, $h$, $k$}\\
0& \quad \text{otherwise}
\end{cases}.
\]
Thus  the matrix $Q$ is the Schur-square of the matrix $W$ and hence $Q$ is orthostochastic.  For every $i\in\mathbb N$, $i= n_h^{(r)}$ for a unique $r$ and $h$ and
\[
(Q\eta)_i=\sum_{j=1}^\infty Q_{ij}\eta_j=  
\sum_{k=1}^\infty Q(r)_{hk}\eta_{m_k^{(r)}} = \xi_{n_h^{(r)}}= \xi_i
\]
i.e., $\xi=Q\eta$.
\ep
Thus combining the above theorem with Theorem \ref {T:3.9},  we have the following infinite dimensional extension of the Horn Theorem
\bR{R:5.4}
 If  $\xi, \eta \in \co*$ then $\xi=Q\eta $ for some orthostochastic matrix $Q \Longleftrightarrow \begin{cases} \xi \prec \eta &\text{ if } \xi \not\in \ell^1\\
 \xi \preccurlyeq \eta &\text{ if } \xi \in \ell^1.
\end{cases}$
\eR
\section{Diagonals of compact operators}\label{S: 6} 

Recall that having fixed an orthonormal basis for the Hilbert space $H$, we denote by $\mathscr D$ the masa of  diagonal operators,
by $E: B(H) \to \mathscr D$ the operation of taking the main diagonal, i.e., the normal faithful and trace preserving conditional expectation from $B(H)$ onto $\mathscr D$, 
and by $\diag: \ell^\infty \to \mathscr D$ the isometric isomorphism that maps a sequence $\eta$ to the diagonal operator having diagonal $\eta$.

For the readers' convenience, let us first collect here some of the equivalent conditions of the majorization relations that we have obtained explicitly or that are immediate consequences of results obtained in the previous sections. 

\bC{C:6.1} Let $\xi, \eta\in \co*.$\\
\noindent If $\xi\not \in (\ell^1)^*$, then  the following conditions are equivalent:
\item [NS(i)] $\xi \prec \eta$
\item [NS(ii)] $\xi=Q\eta$ for some orthostochastic matrix $Q$.
\item [NS(ii$'$)] $\diag \xi = E(U\diag\eta \1 U^*)$ for some orthogonal matrix $U$.
\item [NS(iiid)] $\xi=Q\eta$ for some doubly stochastic matrix $Q$.
\item [NS(iiis)] $\xi=Q\eta$ for some substochastic matrix $Q$.
\item [NS(iii$'$is)]  $\diag \xi = E(W\diag\eta \1 W^*)$ for some isometry $W$.
\item [NS(iii$'$cs)]  $\diag \xi = E(W\diag\eta \1 W^*)$ for some co-isometry $W$.
\item[NS(iii$'$cn)] $\diag \xi = E(L\diag\eta \1 L^*)$ for some contraction $L$.\\
\noindent  If $\xi \in (\ell^1)^*$, then  the following conditions are equivalent:\
\item[S(i)] $\xi \preccurlyeq \eta$
\item [S(ii)] $\xi=Q\eta$ for some orthostochastic matrix $Q$.
\item [S(ii$'$)] $\diag \xi = E(U\diag\eta \1 U^*)$ for some orthogonal matrix $U$.
\item [S(iiid)] $\xi=Q\eta$ for some doubly stochastic matrix $Q$.
\item [S(iiic)] $\xi=Q\eta$ for some column stochastic matrix $Q$.
\item [S(iii$'$is)]  $\diag \xi = E(W\diag\eta \1 W^*)$ for some isometry $W$.\\
\noindent In general the following conditions are equivalent:
\item[(i)] $\xi\prec \eta.$
 \item[(ii)] $\xi= DQ\eta$ for some orthostochastic matrix $Q$ and some $D\in \mathscr D$ with $0\le D \le I$.
 \item[(iii)]  $\xi= QD\eta$ for some orthostochastic matrix $Q$ and some $D\in \mathscr D$ with $0\le D \le I$.
 \item [(ii$'$)] $\diag \xi = E(DU\diag\eta \1 (DU)^*)= $ for some orthogonal  $U$ and some $D\in \mathscr D$ with $0\le D \le I$.
\item [(iii$'$)] $\diag \xi = E(UD\diag\eta \1 (UD)^*)$ for some orthogonal   $U$ and some $D\in \mathscr D$ with $0\le D \le I$.
\item [(iv)]  $\diag \xi = E(L\diag\eta \1 L^*)$ for some contraction $L$.
\eC

\bR {R:6.2}
\item [(i)] In lieu of orthogonal  matrices it suffices to take unitary matrices and conversely, we can always ask that the isometries, co-isometries, and contractions have real entries. 
\item[(ii)] The equivalence of (i) and (iv) was proven by Arveson and Kadison in \cite [Theorem 4.2] {AK02}
\eR 

In the case of finite rank positive operators, the ``sequence" formulation S(i)$\Leftrightarrow$S(ii$'$) in the above corollary is obviously equivalent to the operator theory formulation of the classical  Schur-Horn Theorem (see Theorem \ref {T:1.1} (iii)).   

For infinite rank positive compact operators however, the  above corollary leads to a somewhat different  operator theory reformulation. We find convenient to introduce the notion of partial isometry orbit of an operator.

For every operator $A\in B(H)$, denote by  
\ba 
  \mathscr U(A):&= \{UAU^* \mid U\text { unitary in  $B(H)\}$ \quad  \textit{the unitary orbit} of $A$.}\\  
\mathscr V(A):&= \{VAV^* \mid V\text { partial isometry  in $B(H), ~V^*V= R_A\vee R_{A^*}\}$ \quad \textit{the  partial isometry orbit} of $A.$}
\end{align*} 

\bL{L:6.3}  $ \mathscr U(A) \subset  \mathscr V(A)$  for every  $A\in B(H)$  and  $\mathscr U(A)=  \mathscr V(A)$ if and only if $A$ has finite rank. 
\eL
\bp
The inclusion is obvious since if $U$ is unitary, then $V:=U(R_A\vee R_{A^*})$ is  a partial isometry with $V^*V= R_A\vee R_{A^*} $ and $UAU^*= VAV^*$.  

If $A$ has finite rank and $V$ is a partial isometry with $V^*V= R_A\vee R_{A^*}$, then  both $(V^*V)^\perp$ and $ (VV^*)^\perp $ are infinite  and hence  equivalent projections. Thus  $V$ can be extended to a unitary $U$ and  $VAV^*= UAU^*$. Thus  $\mathscr U(A)=  \mathscr V(A)$.

If $A$ has infinite rank, then choose a partial isometry $V$  with $V^*V= R_A\vee R_{A^*}$ but $(VV^*)^\perp \not \sim (V^*V)^\perp$ and let $B=VAV^*$. Since $VR_AV^*B= VR_AV^*VAV^*= B$, it follows that $VV^*\ge VR_AV^*\ge  R_B$. Similarly, $VR_{A^*}V^*B^*=B^*$ and hence $VV^*\ge VR_{A^*}V^*\ge R_{B^*}$. Thus $VV^*\ge  V(R_A\vee R_{A^*})V^*= VR_AV^*\vee VR_{A^*}V^* \ge  R_B\vee R_{B^*}$. So far, we have only used the fact that $V^*V\ge R_A\vee R_{A^*}$, and since $A=V^*BV$, the same argument shows that $V^*(R_B\vee R_{B^*})V\ge  R_A\vee R_{A^*}$. But then
\[
R_B\vee R_{B^*}\le VV^*= V(R_A\vee R_{A^*})V^*\le VV^*(R_B\vee R_{B^*})VV^*= R_B\vee R_{B^*}
\]
whence $VV^*=R_B\vee R_{B^*}$.
But then $R_B\vee R_{B^*}$ is not unitarily equivalent to $R_A\vee R_{A^*}$, hence $B\not \in \mathscr U(A)$.
\ep
Denote by $s(A)$ the sequence of s-numbers of  $A$. In particular, if $A$ is a positive compact operator,  $s(A)$ is  the eigenvalue list of $A$ in monotone non-increasing order, with repetition according to multiplicity, and with infinitely many zeros added in case $A$ has finite rank. Notice that  the eigenvalue list of $A$ ``ignores"  the null space of $A$, i.e., $s(A)= s(A\oplus 0_n)$  where $0_n$ denotes the zero operator on a space of dimension $n \in \mathbb N\cup\{\infty\}$.

Since for all $A\in K(H)^+$ there is an isometry $V$ for which  $A=V\diag s(A)V^*$ (and $V$ can be chosen unitary when $A$ has finite rank or when $R_A=I$), we see that   
\be{e:32}
\mathscr V(A)= \mathscr V(\diag s(A))= \{ B \in K(H)^+ \mid s(B)=s(A) \} =\bigcup_{0\le n \le \infty}\mathscr U(\diag s(A)\oplus 0_n).
\ee

In \cite {AK02} the set $\mathscr V(A)$ is denoted by  $\mathscr O_{s(A) }$ and 
 it was shown that 
$\mathscr V(A)=\overline{\mathscr U(A))}^{||.||_1}$ when $A$ is trace-class,\cite [Proposition 3.1] {AK02}.

The  infinite dimensional extensions  of the Horn Theorem obtained in Theorems  \ref {T:3.9} and \ref{T:5.3} provide a characterization  of ``the diagonals"  of the partial isometry orbit of a positive compact operator.

 \bP{P:6.4}
Let $A\in K(H)^+$. Then 
\[
E(\mathscr V(A))=\begin{cases}
  \{B\in \mathscr D\cap K(H)^+\mid s(B) \prec s(A) \}\setminus \mathscr L_1 &\text{if } \tr(A)= \infty\\
   \{B\in \mathscr D\cap K(H)^+ \mid s(B) \preccurlyeq s(A)\} &\text{if }  \tr(A)< \infty.\\
  \end{cases}
  \]
\eP
\bp
By \eqref{e:32}, assume without loss of generality that $A= \diag s(A)$.

To prove that the left hand set is  contained in the right hand set, let  $V$ be a partial isometry with $V^*V=R_{\diag s(A)}$  and let $B= E(V\diag s(A) V^*)$. Then $\tr(B)=\tr(V\diag s(A) V^*)=\tr(\diag s(A))=\tr(A)$ and we only need to prove that $s(B) \prec s(A)$.   Let $W$ be the isometry that rearranges the sequence $s(B)$ into the sequence of the diagonal entries of $B$, i.e.,  $B= W\diag s(B) W^*$. Then $W$ commutes with the expectation  $E$ and hence
\[
\diag s(B)=  W^* B W= W^* E(V\diag s(A) V^*) W= E(W^*V\diag s(A) V^*W).
\]
If  $\diag s(A)$ has infinite rank then $V$ is an isometry and if $\diag s(A)$ has finite rank, we can extend $V$ to an   isometry. Since $R_W= R_B=R_{E(VV^*)}\ge VV^*= R_V$, $W^*V$ is an isometry. By Lemma \ref {L:2.3},  $s(B) = Q s(A)$ where $Q_{ij}= |(W^*V)_{ij}|^2$  and hence by  \eqref{e:1} and Lemma  \ref {L:2.4}  it follows that  $s(B) \prec s(A)$.

To prove the opposite inclusion, let  $B\in \mathscr D\cap K(H)^+$ with $s(B) \prec s(A)$ and if $\tr(A)=\infty$  assume that $\tr(B)=\infty$ while if $\tr(A)<\infty$ assume that $\tr(B)=\tr(A)$ and hence $s(B) \preccurlyeq s(A)$. Then by Theorems   \ref {T:3.9}  and \ref{T:5.3}, there is a unitary $U$ for which $\diag s(B)= E(U\diag s(A) U^*)$.  
As above, let $W$ be the isometry that rearranges the sequence $s(B)$ into the sequence of the diagonal entries of $B$. Then $WU$ is an isometry and hence
\[
B= W\diag s(B) W^* =WE(U\diag s(A) U^*)W^*= E(WU\diag s(A) U^*W^*)\in E(\mathscr V(\diag s(A)) =E(\mathscr V(A)).
\]

 \ep

\bR{R:6.5}
\item[(i)] The proof that $s(E(C))\prec s(C)$ for all $C\in K(H)^+$ (i.e., the inclusion  of the left-hand set into the right hand set in Proposition \ref {P:6.4})  is usually attributed to Ky Fan \cite{Fa49}.  See also \cite{Fa51} and see an elegant proof in \cite [Theorem 4.2]{AK02} of the more general fact that $s(E(LAL^*))\prec s(A)$ for every contraction $L$.
\item[(ii)] Gohberg and Markus have proven in \cite [Proposition III, pg 205] {GiMa64} that if $A\in K(H)^+$,   $\xi \in \co*$, and $\xi \prec s(A)$,  then there is an orthonormal sequence $f_n\in H$ for which $(Af_n,f_n)=\xi_n$ for all $n$. Thus setting  $W^*e_n=f_n$ for  the fixed orthonormal basis $\{e_n\}$, defines a co-isometry  $W$  for which $\diag \xi = E(WAW^*)$. Applied to the case of $A=\diag \eta$, their result proves that if $\xi \prec \eta$, then $\xi=Q\eta$ for some co-isometry stochastic matrix $Q$ (cf. Theorem \ref {T:3.7}). 

In the case that $A$ is of  trace-class  and  $\xi \preccurlyeq s(A)$,  Gohberg and Markus have furthermore proven in \cite[Theorem 1] {GiMa64} that $A$ vanishes on span$\{f_n\}^\perp$, i.e., that $W^*W\ge R_A$. Thus  $V:=WR_A$ is a partial isometry  and $\diag \xi \in E(\mathscr V(A))$, which proves  the inclusion of the right hand set into the left hand set in Proposition \ref {P:6.4}.
 \item[(iii)] 
In the trace-class case, Proposition \ref{P:6.4}   was  derived from the (finite dimensional) Schur-Horn Theorem  by Arveson and Kadison  using compactness arguments (see \cite[Theorem 4.1]{AK02}). 
\eR

Since  $s(A\oplus 0)=s(A)$ but in general  $E(\mathscr U(A\oplus 0))\ne  E(\mathscr U(A))$,  it is clear that we cannot characterize   $E(\mathscr U(A))$  only in terms of the sequence $s(A)$.  

\bP{P:6.6} Let $A\in K(H)^+$. Then  
\item[(i)]  $E(\mathscr U(A)) \subset \ E(\mathscr V(A)) \cap \{B\in \mathscr D \mid \tr (R_B^\perp ) \le  \tr (R_A^\perp)\}.$
\item[(ii)] If  $R_A=I$, then 
\begin{multline}
E(\mathscr U(A)) = E(\mathscr V(A)) \cap \{B\in \mathscr D \mid R_B=I\} \\=\begin{cases}
  \{B\in \mathscr D\cap K(H)^+\mid s(B) \prec s(A), R_B=I \}\setminus \mathscr L_1 &\text{if } \tr(A)= \infty\\
   \{B\in \mathscr D\cap K(H)^+ \mid s(B) \preccurlyeq s(A), R_B=I\} &\text{if }  \tr(A)< \infty.\\
  \end{cases}\notag
\end{multline}
\item[(iii)] The inclusion in (i) is proper unless $R_A=I$ or $A$ has finite rank.
\eP

\bp
\item[(i)] That $E(\mathscr U(A)) \subset \ E(\mathscr V(A))$ follows from Lemma \ref {L:6.3}.  Let $B= E( UAU^*)$ for some unitary $U$. Then 
$E( R_B^\perp UAU^*R_B^\perp)=  R_B^\perp E( UAU^*)R_B^\perp   =0$ since $R_B\in \mathscr D$.  By the faithfulness of the expectation, it follows that $ R_B^\perp UAU^*R_B^\perp=0$, hence  $R_B^\perp UAU^*=0$ and thus  $ R_B^\perp\le (R_{UAU^*})^\perp = UR_A^\perp U^*$, whence   $\tr(R_B^\perp)\le 
\tr(R_A^\perp)$.
\item[(ii)] The second set equality is given by Proposition \ref{P:6.4}  . By (i), $E(\mathscr U(A))\subset E(\mathscr V(A)) \cap \{B\in \mathscr D \mid R_B=I\} $. To prove the opposite inclusion, let $B\in \mathscr D\cap K(H)^+$ with $R_B=I$ and $s(B) \prec s(A)$ and assume that $s(B) \preccurlyeq s(A)$ if $\tr(A)< \infty$ and $\tr(B)= \infty$ if $\tr(A)= \infty$. 
By Theorems \ref {T:3.9}  and \ref {T:5.3} (or see for convenience Corollary \ref {C:6.1}), $\diag s(B) \in E(\mathscr U(\diag s(A))$.  Since $R_B=I$, there is a permutation matrix $\Pi$ for which $B= \Pi \diag s(A) \Pi^*$.  As $\Pi$ commutes with the expectation,  $B\in  E(\mathscr U(\diag s(A))$. 
But $\mathscr U(\diag s(A)= \mathscr U(A)$ since $R_A=I$, and thus  $B\in E(\mathscr U(A)).$ 
\item[(iii)] If $R_A=I$, the equality holds by (ii) and if $A$ has finite rank, then $\mathscr U(A)= \mathscr V(A)$ by Lemma \ref {L:6.3}  and hence 
\[
E(\mathscr U(A)) =  E(\mathscr V(A))= E(\mathscr V(A)) \cap \{B\in \mathscr D \mid \tr (R_B^\perp ) \le \infty\}.
\]
Assume now that $A$ has infinite rank but $R_A\ne I$. Set $n:= \tr (R_A^\perp)\in \mathbb N \cup \{\infty\}$ and $\eta:= s(A)$. Then  $\eta_j\ne 0$ for all $j$  and $n\ne 0$.
Let  $\Lambda$ be an infinite subset of $\mathbb N$ with $\card \mathbb N\setminus \Lambda = n$,  let  $\pi:  \mathbb N \to \Lambda$ be a bijection, $\tilde{\eta_k}:= \begin{cases} 0 & \text {if }k\not \in \Lambda\\
\eta_j &\text {if } k= \pi(j) 
\end{cases}$.
Then 
$
\mathscr U(A) = \mathscr U(\diag\tilde{\eta}).
$
Clearly, $B:=\diag \eta \in E(\mathscr V(A))= E(\mathscr V(\diag \eta))$ and $R_B^\perp =0$.

We claim that $B\not\in E(\mathscr U(A))$. Reasoning by contradiction, assume that  $\diag \eta = E(U\diag \tilde{\eta}U^*)$  for some unitary $U$. By Lemma \ref {L:2.3}, $\eta = Q\tilde{\eta}$ for the unistochastic and hence doubly stochastic matrix $Q$ given by $Q_{ij}= |U_{ij}|^2$. But then, for every $i\in \mathbb N$
\[
\sum _{j=1}^\infty Q_{i\pi(j)} \eta_j = \sum _{j=1}^\infty Q_{i\pi(j)}\tilde{\eta}_{\pi(j)} =\sum_{k=1}^\infty Q_{ik}\tilde{\eta_k}=   \eta_i = \sum_{k=1}^\infty Q_{ik}\eta_i =\sum_{k\not \in \Lambda}Q_{ik}\eta_i + \sum _{j=1}^\infty Q_{i\pi(j)}\eta_i
\]
since $Q$ is column stochastic.  Hence for all $i\in \mathbb N$,
\be{e:33}
\sum_{k\not \in \Lambda}Q_{ik}\eta_i + \sum _{j=1}^\infty Q_{i\pi(j)}(\eta_i-\eta_j)=0\
\ee
Let $n_p$ be the strictly increasing sequence of integers starting with $n_0=0$ for which $\eta_j= \eta_{n_p}$ for all $n_{p-1} < j \le  n_{p}$. Applying the identity \eqref{e:33} to any $0< i \le n_1$ we have 
\[
\sum_{k\not\in \Lambda}Q_{ik}\eta_i + \sum _{j=n_1+1}^\infty Q_{i\pi(j)}(\eta_i- \eta_j)=\sum_{k\not\in \Lambda}Q_{ik}\eta_i + \sum _{j=1}^\infty Q_{i\pi(j)}(\eta_i- \eta_j)= 0.
\]
Since $\eta_i=\eta_{n_1}> \eta_j> 0$ for all $j> n_1$,  and $ Q_{ik}\ge 0$ for all $k$, we see that $Q_{ik}=0$ for all $k\not \in \{\pi(1),\cdots, \pi(n_1)\}$ and in particular  for all $k \not \in \Lambda$.  But then  $\sum_{j=1}^{n_1}Q_{i\pi(j)}=1$ since $Q$ is row stochastic. Hence $\sum_{i,j=1}^{n_1}Q_{i\pi(j)}=n_1$. Since $Q$ is also column-stochastic, it follows that
$\sum_{i=1}^{n_1}Q_{i\pi(j)}=1$ for every $0< j\le n_1$. Thus $Q_{i\pi(j)}=0$  for every pair $i> n_1$ and $0< j\le n_1$.
Now applying the  identity \eqref{e:33} to  $n_1< i \le n_2$ we obtain
\[
\sum_{k\not\in \Lambda}Q_{ik}\eta_i +  \sum _{j=n_2+1}^\infty Q_{i\pi(j)}(\eta_i- \eta_j)=\sum_{k\not\in \Lambda}Q_{ik}\eta_i + \sum _{j=1}^\infty Q_{i\pi(j)}(\eta_i- \eta_j)= 0.
\]
Thus again we obtain  for  all $n_1< i \le n_2$ that  $Q_{ik}=0$ for all $k\not\in \{\pi((n_1+1), \cdots \pi(n_2)\} $ and in particular for all $k\not \in \Lambda$. Iterating, we obtain that  $Q_{ik}=0$ for all $k\not\in \Lambda$ and all $i$. Since $\Lambda \ne \emptyset$ we conclude that $Q$ is not column stochastic, a contradiction.
\ep

Notice that if $\xi, \zeta \in (\text{c}_{\text{o}})^+$, then  $(\xi+\zeta)^*\prec \xi^*+\zeta^*$. Thus the majorization condition in  Proposition \ref {P:6.4} is preserved by convex combinations, which yields the following simple conclusions.

\bC{C:6.7} Let $A\in K(H)^+$. Then
\item[(i)]  $E(\mathscr V(A))$ is convex.
\item[(ii)] If $A$ has finite rank or $R_A=I$, then $E(\mathscr U(A))$ is convex.
\eC
\bp
\item[(i)]  Let $B_i\in \mathscr D\cap K(H)^+$ and $s(B_i) \prec s(A)$ for $i=1, 2$,  let $t\in [0,1]$, and let $B:=  tB_1+(1-t)B_2$. Then $B\in \mathscr D\cap K(H)^+$. Clearly,  $s(tB_1) \prec s(tA)$, $s((-t)B_2) \prec s((1-t)A)$, and hence $s(B) \prec  s(tA)+ s((1-t)A)= s(A)$. Furthermore, if $\tr(A)=\infty$ and we choose $B_i\not \in \mathscr L_1$, then $B\not \in  \mathscr L_1$, while if we $\tr(A)<\infty$ and we choose $B_i$ with $\tr(B_i)=\tr(A))$, then also $\tr(B)=\tr(A)$. Thus the right hand set in Proposition \ref {P:6.4} is convex  and hence so is $E(\mathscr V(A))$.
\item[(ii)] If in the proof of (i) we choose the diagonal operators $B_i$ so that $R_{B_i}=I$ for $i=1,2$, i.e., all the diagonal entries of $B_i$ don't vanish, then $R_B=I$ and hence the conclusion follows from 
Proposition \ref {P:6.6} 
\ep
When $A$ has infinite rank but $R_A\ne I$, we can identify a distinguished subset of $E(\mathscr U(A))$ in terms of the following stronger notion of sequence majorization.

\bD{D:6.8} Let $ \xi\in \co*, \eta\in (\text{c}_{\text{o}})^+$, $p\in \mathbb Z_+$, and $ N\in \mathbb N$. Then we say that  
\item[(i)] $\xi \overset{p}{\underset{}{\prec}} \eta$ if $\xi \prec \eta^*$ and $\sum_{j=1}^{n+p}\xi_j\le \sum_{j=1}^n\eta_j $ for all $n\ge N$.
\item[(ii)] $\xi \overset{p}{\prec} \eta$ if $\xi \overset{p}{\underset{}{\prec}} \eta$ for some $ N\in \mathbb N.$
\eD

 \bL {L:6.9}
 Let $\xi, \eta \in \co*$ and  $\xi \overset{p}{\prec} \eta$ for some $p \in \mathbb N$,  and assume also that  $\xi\preccurlyeq \eta $  if  $\xi\in \ell^1$.  Then $\xi = Q< \overbrace{0, 0, \cdots, 0}^p, \eta>$ for some  orthostochastic matrix $Q$.
 \eL
\bp
We start by  disposing of the case where $\xi$ or $\eta$ or both have finite support.  If $\xi$ has finite support, by hypothesis   $\xi \preccurlyeq \eta$ and hence $\eta$ too has finite support.  If $\eta$ has  finite support, then there is a permutation matrix $\Pi$ for which $\Pi< \overbrace{0, 0, \cdots, 0}^p, \eta>= \eta$. By Lemma \ref {L:2.7}, there is an  orthostochastic matrix $Q'$ for which $\xi =  Q'\eta= Q'\Pi < \overbrace{0, 0, \cdots, 0}^p, \eta>$. By Lemma \ref {L:3.4}  $Q:= Q'\Pi$ is also orthostochastic. 

Thus assume henceforth that both $\xi_n$ and $\eta_n$ never vanish. Let  $N$ be an integer for which $\xi \overset{p}{\underset{}{\prec}} \eta$ and let $\tilde \eta: =<\eta_1, \eta_2, \cdots, \eta_N, \overbrace{0, 0, \cdots, 0}^p, \eta^{(N)}>$ (recall the notation $\eta^{(N)}=< \eta_{N+1},  \eta_{N+2} , \cdots>.$) As in the first part of the proof, since $\tilde \eta$ is the permutation of the original sequence $< \overbrace{0, 0, \cdots, 0}^p, \eta> $, it suffices to find an  orthostochastic matrix $Q$ for which $ \xi = Q\tilde \eta$.   

Consider first the case when $\xi_1> \eta_N$ and  apply the first step of the construction  in the proof of Theorem \ref{T:3.7} to the sequences $\xi\prec \eta$. We have proven in  \eqref {e:11} that $\xi^{(1)} \prec \rho(1)$. The assumption  that   $\xi_1> \eta_N$ guarantees that $m_1\le N-1$ and hence for every $n\ge N-1$ we see from  Lemma \ref{L:4.1} (iii) or from a direct elementary computation that 
\[
\sum_{j=1}^n\rho(1)_j= \sum_{j=1}^{n+1}\eta_j-\xi_1\ge  \sum_{j=1}^{n+1+p}\xi_j-\xi_1= \sum_{j=1}^{n+p}\xi^{(1)}_j.
\] 
This shows that  $\xi^{(1)}\overset{p}{\underset{N-1}{\prec}} \rho(1)$. Furthermore, if $\xi\in \ell^1$ and thus $\xi \preccurlyeq \eta$, then $\xi^{(1)} \preccurlyeq \rho(1)$. 

Recall that $R^{(1)}$ is the Schur-square of $V(m_1,t_1)\oplus I_\infty$, where $V(m_1,t_1)$ is the orthogonal $m_1+1\times m_1+1$ matrix  defined in \eqref {e:7}. Then $R^{(1)}$ is orthostochastic  and  for every $j$ there is at most one index $i>1$ for which $R^{(1)}_{ij}\ne 0$. Since $R^{(1)}\eta = < \xi_1, \rho(1)>$ and   $m_1+1\le N$, we also have
\[
R^{(1)}\tilde \eta =  <\xi_1, \rho(1)_1, \rho(1)_2, \cdots, \rho(1)_{N-1}, \overbrace{0, 0, \cdots, 0}^p, \rho(1)^{(N-1)}>.
\]

Now, consider the case when  $\xi_1\le \eta_N$ and define $t_1:= \frac{\eta_N}{\xi_1}$,  $\tilde{R}^{(1)}$ to be the Schur-square of  $V(N,t_1)\oplus I_\infty$, and $\rho(1):=<\eta_1, \eta_2, \cdots, \eta_{N-1}, \eta_N-\xi_1, \eta^{(N)}>^*$.  Then $\rho(1)_j=\eta_j$ if $j\le N-1$ and $\sum_{j=1}^n\rho(1)_j\ge \sum_{j=1}^{n}\eta_j -\xi_1$ for all $n\ge N$. Thus 
\[
\sum_{j=1}^n\rho(1)_j\quad \begin{cases}=\quad  \sum_{j=1}^n\eta_j \ge \sum_{j=1}^n\xi_j\ge  \sum_{j=1}^n\xi^{(1)}\quad&\text{if }n\le N-1\\
\ge \quad \sum_{j=1}^{n}\eta_j -\xi_1\ge \sum_{j=1}^{n+p}\xi_j -\xi_1=  \sum_{j=1}^{n+p-1}\xi^{(1)}&\text{if } n\ge N
\end{cases}.
\]
This shows that $\xi^{(1)}\overset{p-1}{\underset{}{\prec}} \rho(1)$  and  if $\xi\in \ell^1$, then clearly, $\xi^{(1)}\preccurlyeq \rho(1)$. Furthermore, since  $t_1\eta_N=\xi_1$ and $(1-t_1)\eta_N= \eta_N-\xi_1$, we have
\[
\tilde{R}^{(1)}\tilde\eta=
<\xi_1, \eta_1, \eta_2, \cdots, \eta_{N-1}, \eta_N-\xi_1, \overbrace{0, 0, \cdots, 0}^{p-1},\eta^{(N)}>.
\]
Let $\Pi_1$ be the permutation matrix for which 
\[
\Pi_1<\eta_N-\xi_1, \overbrace{0, 0, \cdots, 0}^{p-1},\eta^{(N)}>= < \rho(1)_{}, \overbrace{0, 0, \cdots, 0}^{p-1},\rho(1)^{(N)}>
\]
and let $R^{(1)}:= \big(I_N\oplus\Pi_1\big)\tilde{R}^{(1)}$. Then
\[
R^{(1)}\tilde\eta= 
<\xi_1, \rho(1)_1, \rho(1)_2, \cdots, \rho(1)_{}, \overbrace{0, 0, \cdots, 0}^{p-1}, \rho(1)^{(N)}>.
\]
Furthermore, by Lemma \ref {L:3.4},  $R^{(1)}$ is orthostochastic and by Lemma \ref {L:3.5} (ii) for every $j$ there is at most one index $i>1$ for which $R^{(1)}_{ij}\ne 0$. 

Iterate this construction. At every step we decrease by one unit either $N$ or  $p$. Thus  we end the process when after $k\le p+N-1$ steps we reach $p=0$. Define as in  \eqref{e:9} 
\[
Q^{(k)}:=\big(I_{k-1}\oplus R^{(k)}\big)\big(I_{k-2}\oplus R^{(k-1)}\big)\cdots \big(I_1 \oplus R^{(2)}\big)R^{(1)}
\]
where for every $1\le h\le k$ and every $j$  there is at most one index $i>1$ for which $R^{(h)}_{ij}\ne 0$.   By Lemma \ref  {L:3.5}~(iii), as in the proof of Proposition \ref {P:3.6},  $Q^{(k)}$ is orthostochastic and for every $j$ there is at most one index $i>k$ for which $Q^{(k)}_{ij}\ne0$.
Then $Q^{(k)}\tilde\eta = < \xi_1, \xi_2, \cdots, \xi_k, \rho(k)>$ and $\xi^{(k)}\prec \rho(k)$, with $\xi^{(k)}\preccurlyeq \rho(k)$ when $\xi\in \ell^1$.
By Theorem \ref {T:3.9} there is an orthostochastic matrix $Q'$ for which $\xi^{(k)}= Q'\rho(k)$. Let $Q:=\big(I_{k}\oplus Q'\big) Q^{(k)}$. Then $\xi =Q\tilde\eta$.  Finally, by Lemmas \ref  {L:3.5}(i) and  \ref {L:3.4}, $Q$ is also orthostochastic.

\ep

 \bP{P: 6.10}  
 Let $A\in \mathscr K(H)^+$ and assume that $A$ has infinite rank.  Then
 \item[(i)]
\[
\bigcup_{p\in \mathbb Z_+ , ~0\le p \le  \tr (R_A^\perp)}\{E(\mathscr V(A)) \cap \{B\in \mathscr D\cap K(H)^+ \mid s(B) \overset{P}{\prec } s(A),~ \tr  (R_A^\perp)-p\le \tr(R_B^\perp ) \le  \tr (R_A^\perp)\}\\  \subset E(\mathscr U(A))  \]
\item[(ii)] If $R_A\ne I$, then the inclusion in (i) is proper.
 \eP
 
\bp
\item [(i)] Let $B\in \mathscr D\cap  K(H)^+$ belong  to the left-hand set and by passing  if necessary  to a smaller integer, assume without loss of generality that  $s(B) \overset{P}{\prec } s(A)$ with $\tr(R_B^\perp)= \tr  (R_A^\perp)-p$. Since $\tr(B)= \tr(A)$, if $s(A)\in \ell^1$, i.e., $\tr(A)< \infty$, then $s(B)\preccurlyeq s(A)$ (see also Proposition \ref{P:6.4}). But then by Lemma  \ref  {L:6.9}, there is an orthostochastic matrix $Q$ for which  $s(B) = Q< \overbrace{0, 0, \cdots, 0}^p, s(A)>$. Thus in $B(R_BH)$, $\diag s(B)=E_{R_B}(V(\diag s(A)\oplus 0_p)V^*)$ where $V\in B(R_BH)$ is  unitary  and $E_{R_B}$ denotes the conditional expectation with respect to the orthonormal basis of $R_BH$, which is the subset of the fixed basis on $H$ since $R_B\in \mathscr D$. Extend $V$ to a unitary $U\in B(H)$ that commutes with $R_B$. Then
\[
\diag s(B) \oplus 0_{\tr(R_B^\perp)}= E(V(\diag s(A)\oplus 0_p)V^*\oplus 0_{\tr(R_A^\perp)-p})= E(U(\diag s(A)\oplus 0_{\tr(R_A)^\perp})U^*)\in E(\mathscr U(A)).
\]
Finally, $B=\Pi\big(\diag s(B) \oplus 0_{\tr(R_B^\perp)}\big)\Pi^*$ for some permutation matrix $\Pi$. But then 
\[
B=\Pi E(\mathscr U(A))\Pi^*= E(\Pi(\mathscr U(A)\Pi^*)=E(\mathscr U(A)).
\]
\item [(ii)] Assume that $R_A\ne I$ and let $N:= \tr(R_A^\perp)$. Set $\eta:=s(A)$ and $\tilde \eta:= <0, \eta_1, \eta_2, \cdots>$. Then $A$ is unitarily equivalent to the diagonal operator
\[
D:= \begin{cases}\diag \tilde  \eta \quad &\text{if } N=1\\
0_{N-1}\oplus \diag \tilde  \eta, &\text{if } N>1.
\end{cases}
\]
Let $U$ be an  orthogonal matrix for which $U_{ij}=0$ for all $ i>j>1$, all other entries being nonzero.  (such matrices exist, see Example \ref {E:6.11}  below) and  let $ Q$ be the Schur square of $U$.  Let 
\[
\tilde U:= \begin{cases} U&\text{if } N=1\\ I_{N-1}\oplus  U&\text{if } N>1\end{cases} 
\quad \text {and} \quad   B:=E(\tilde UD\tilde U^*)
\] 
and let $\xi:=  Q\tilde \eta$. Then 
 \[
 B=\left \{ \begin{array}{lll} E(U \diag \tilde \eta~ U^*) &=  \diag \xi   \quad &\text {if } N=1\\
E(0_{N-1}\oplus U \diag \tilde \eta~ U^*)  &=  0_{N-1}\oplus \diag \xi    &\text{if  } N>1.
 \end{array}\right. 
  \]
It is immediate to verify that $\xi_n > 0$ for all $n$ and hence $s(B)$ is the monotone rearrangement $\xi^*$ of $\xi$. Now  $B\in E(\mathscr U(A))$ and $ \tr(R_B^\perp)= N-1=  \tr  (R_A^\perp)-1 $ but $\xi^*\not  \overset{1}{\prec} \eta$  and hence $s(B)\not  \overset{p}{\prec} s(A)$ for all  $p>0$. Indeed, since $\eta_j>0$ for all $j$,  $ Q_{ij} \begin{cases} =0& \text{ for }i> j  >1 \\ > 0& \text{ otherwise }\end{cases}$, and  $ Q$ is doubly stochastic and hence column stochastic,  we have for every $n>1$

\ba
\sum_{i=1}^n\xi^*_i&\ge \sum_{i=1}^n\xi_i= \sum_{i=1}^n\sum_{j=1}^\infty  Q_{ij}\tilde \eta_j\\
&= \sum_{i=1}^n\sum_{j=2}^\infty  Q_{ij}\eta_{j-1}
=\sum_{j=2}^n  \sum_{i=1}^n Q_{ij}\eta_{j-1} + \sum_{j=n+1}^\infty  \sum_{i=1}^n Q_{ij}\eta_{j-1} \\
&= \sum_{j=1}^{n -1}\eta_j+  \sum_{j=n+1}^\infty  \sum_{i=1}^n Q_{ij}\eta_{j-1}
 >  \sum_{j=1}^{n -1}\eta_j.
\end{align*}
\ep

\begin{example}\label{E:6.11} An orthogonal matrix  $U$ for which $ U_{ij}=0$ for all $ i>j>1$, all other entries being nonzero.  Given a sequence $\{a_n\}$  with $a_n>0$  for all $n$ and $\sum_{j=1}^\infty a_n=1$, set $b_n:=\big(\sum_{j=1}^n a_j\big)^{-1}$ and 
\[
 U:=
\begin{pmatrix}
\sqrt{a_1}&\sqrt{a_1(b_1-b_2)}  & \sqrt{a_1(b_2-b_3)} & \cdots &\sqrt{a_1(b_{n-1}-b_n)} &\sqrt{a_1(b_n-b_{n+1})}&\cdots\\
\sqrt{a_2}&-\sqrt{1-a_2b_2}&  \sqrt{a_2(b_2-b_3)} & \cdots &\sqrt{a_2(b_{n-1}-b_n)} &\sqrt{a_2(b_n-b_{n+1})}&\cdots\\
\sqrt{a_3}&0&-\sqrt{1-a_3b_3}& \cdots &\sqrt{a_3(b_{n-1}-b_n)} &\sqrt{a_3(b_n-b_{n+1})}&\cdots\\
\vdots&\vdots&\vdots&\vdots&\vdots&\vdots\\
\sqrt{a_n}&0&0&\cdots&-\sqrt{1-a_{n}b_{n}}&\sqrt{a_n(b_n-b_{n+1})}&\cdots\\
\sqrt{a_{n+1}}&0&0&\cdots&0&-\sqrt{1-a_{n+1}b_{n+1}}&\cdots\\
\vdots&\vdots&\vdots&\vdots&\vdots&\vdots&\vdots
\end{pmatrix}
\]
A direct computation shows that $U$ is indeed  unitary. 
\end{example}

Thus our characterization of $E(\mathscr U(A))$ is complete for the cases when $R_A=I$ or when $A$ has finite rank and it points to an interesting and delicate role of $R_A$ in the general case.  

The role of $R_A$ disappears when one takes the  closure of $E(\mathscr U(A))$ under the operator or the trace-class norm (in the trace class case).  Indeed from the more general work by A. Neumann  in \cite {Na99} follows that 
\[
\overline{E(\mathscr U(A))}^{||.||}= \{B\in \mathscr D\cap K(H)^+ \mid s(B) \prec s(A)\} \quad \text{ \cite [Corollary 2.18, Theorem 3.13]{Na99}}
\]

In case that  $\tr (A) < \infty$,  Antezana,  Massey,  Ruiz, and Stojanoff extended the work of Neumann and proved that 
\[
\overline{E(\mathscr U(A)) }^{||.||_1} = \{B\in \mathscr D\cap K(H)^+ \mid s(B) \preccurlyeq s(A)\} \quad \quad \text{\cite[Proposition 3.13] {AMRS} }
\]

For the readers' convenience we collect here below the relations between the various sets and their closures in the special case when $R_A=I$ or $A$ has finite rank  (see also the introduction.) 

\vspace{0.5 cm}

\begin{itemize}
\item $\text{If }A\in K(H)^+, R_A=I, \text{and }\tr(A)=\infty ,\text{ then}$
\vspace{-0.1cm}
\begin{multline*}
E(\mathscr U(A) = \{B\in \mathscr D\cap K(H)^+ \mid s(B) \prec s(A), R_B = I\} \setminus \mathscr L_1 \\
 \subsetneq E(\mathscr V(A))= \{B\in \mathscr D\cap K(H)^+ \mid s(B) \prec s(A)\} \\
  \subsetneq   \overline{E(\mathscr U(A))}^{||.||}= \{B\in \mathscr D\cap K(H)^+ \mid s(B) \prec s(A)\} 
\end{multline*}
\item $\text{If }0\ne A\in K(H)^+, R_A=I, \text{and }\tr(A)<\infty, \text{ then}$
\begin{multline*}
E(\mathscr U(A) = \{B\in \mathscr D\cap K(H)^+ \mid s(B) \preccurlyeq s(A), R_B = I\} \\
 \subsetneq E(\mathscr V(A))= \{B\in \mathscr D\cap K(H)^+ \mid s(B) \preccurlyeq s(A)\}  \\
=E\big(  \overline{\mathscr U(A))}^{||.||_1}\big) =E\big(  \overline{\mathscr U(A))}^{||.||_1}\big)\\ 
\subsetneq \overline{E(\mathscr U(A))}^{||.||}= \{B\in \mathscr D\cap K(H)^+ \mid s(B) \prec s(A)\} 
\end{multline*}
\item $\text{If }0\ne A\in K(H)^+\text{ has finite rank, then}$
\begin{multline*}
E(\mathscr U(A) =E(\mathscr V(A) = \{B\in \mathscr D\cap K(H)^+ \mid s(B) \preccurlyeq s(A)\} \\
=E\big(  \overline{\mathscr U(A))}^{||.||_1}\big) =E\big(  \overline{\mathscr U(A))}^{||.||_1}\big)\\ 
\subsetneq \overline{E(\mathscr U(A))}^{||.||}= \{B\in \mathscr D\cap K(H)^+ \mid s(B) \prec s(A)\} 
\end{multline*}
\end{itemize}

\noindent So, in particular, $E(\mathscr U(A)) $ is closed in the trace norm if and only if $A$ has finite rank (see also \cite[Remark 4.8] {AMRS} and it is never closed in the operator norm but for the trivial  case $A\ne 0$. This  answers in the negative the question by A. Neumann \cite [pg 447]{Na99} on whether $E(\mathscr U(A))$ must be norm closed if $A$ is of trace class or Hilbert-Schmidt class and provides an alternative to his \cite [Remark 3.7] {Na99} .

The inclusions in the finite rank case are illustrated by the following simple example.

\begin{example}\label{E:6.12}Let $\eta:=\,<1, 0, \ldots >$. Then $\mathscr U(\diag \eta)) =\mathscr V(\diag \eta))$ consists of all the rank-one projections, $E(V(\diag \eta))$ consists of all trace class positive diagonal operators with trace $=1$, while  $\overline{E(\mathscr U(\diag \eta))}^{||.||} $ consists of all trace class positive diagonal operators with trace $\le1$. That the latter set contains $0$ is immediately clear since if $F_n$ denotes the infinite matrix with the upper left $n \times n$ corner having entries all equal to $\frac{1}{n}$ and all other entries zero, then $F_n$ is a rank-one projections and hence   $F_n\in \mathscr U(\diag \eta)$ and $||E(F_n)||= \frac{1}{n}\to 0$.
\end{example}

Finally, notice that  in the summable case, if we pass from the unitary orbit of  to the bounded orbit, we obtain:

\bC{C:6.13} If $0\ne\eta\in(\ell^1)^*$, then $\mathscr L_1^+\cap \mathscr D= E\{T\diag \eta\1 T^*\mid T\in B(H)\}$.
\eC
\bp
Let $B\in \mathscr L_1^+\cap \mathscr D$. Then $B= V\diag s(B) \1V^*$ for some isometry $V$ that commutes with $E$. Choose  $c\ge \frac{\tr B}{\eta_1}$. Then  $s(B) \prec c\eta$  and by Corollary \ref {C:6.1}, $\diag s(B)= E(DU\diag c\eta\1 U^*D)$ for some $D\in \mathscr D$ with $0\le D\le I$ and some unitary $U$. Then
\[
B= VE(DU\diag c\eta\1 U^*)V^*= E((c^{1/2}VDU) \diag \eta \1(c^{1/2}VDU)^*).
\]
The opposite inclusion is trivial since  $\tr ( E( T\diag \eta\1 T^*) )= \tr (T\diag \eta\1 T^*) <\infty$ for every $T\in B(H)$.
\ep

\section{Intermediate sequences} \label{S: 7}
The main goal of this section is to extend to infinite sequences the  finite sequence results we list in the following proposition (see \cite[5.A.9 and 5.A.9a]{MO79}). 
As mentioned in Section \ref{S:2}, we use for finite sequences the same notations denoting majorization orderings like $\preccurlyeq$ etc.  introduced in  Definition \ref{D:1.2} for infinite sequences.

\bP{P:7.1} Let $\xi, \eta \in \mathbb (R^N)^+$ be two nonincreasing finite sequences. 
\item[(A)] If $\xi \prec \eta$, then there are finite nonincreasing sequences $\zeta, \rho \in \mathbb (R^N)^+$ with
 \item[(i)]  $\xi\preccurlyeq \zeta \le \eta$;
  \item[(ii)] $\xi \le \rho\preccurlyeq \eta$.
 \item[(B)] If $\xi \prec_\infty\eta$, then there are finite nonincreasing sequences $\zeta, \rho \in \mathbb (R^N)^+$ with
 \item[(i)]   $\xi \preccurlyeq_\infty \zeta \le \eta$, (i.e., $\zeta\preccurlyeq \xi$ and $\zeta \le \eta$);
 \item[(ii)] $\xi \le \rho   \preccurlyeq_\infty \eta$, (i.e., $\xi \le \rho$ and $\eta\preccurlyeq \rho$).
\eP
\noindent(A)(i) is due to Mirsky \cite {Mi60} and the one-line proof given for  $\mathbb R^N$ in \cite [5.A.9]{MO79} easily extends to $\mathbb (R^N)^+$ (see also \cite [Lemma 4.3]{AK02}) and to $(\ell^1)^*$ (see the first step in the proof of Theorem \ref {T:7.4}(i) below). \\
(A)(ii) is due to Fan \cite {Fa51} (see also \cite [5.A.9]{MO79}).\\
(B)(i)-(ii) are stated without proof or attribution in \cite [5.A.9a]{MO79}. 
The proof of (B)(ii) is immediate as it is enough to choose $\rho:=\,< \xi_1 + \sum_{j=1}^N(\eta_j-\xi_j), \xi_2, \ldots,\xi_N>$. 
The proof of (B)(i) is similar to that of (A)(ii) and since we need this result in Proposition  \ref{P:7.6}, for completeness's sake and the readers'  convenience, we present it next.
\bp[Proof of (B)(i)]  Assume that $\sum_{j=1}^N(\eta_j-\xi_j) \ne 0$, otherwise  choose $\zeta= \eta$. The proof is by induction on $N$. If $N=1$ then $\xi_1 <  \eta_1$, thus $\zeta:=\xi$ satisfies the condition. 
Assume now that the property is true for all $k \le N-1$ for some $N > 1$ and let $ \alpha:= \underset{1\le n \le }{\min} \sum_{j=n}^N(\eta_j-\xi_j).$  
In particular, $\eta_N -\alpha\ge \xi_N \ge 0.$ By passing if necessary to $\tilde{\eta}:= <\eta_1, \eta_2, \ldots, \eta_{N-1}, \eta_N-\alpha>$ and since $\xi\prec_\infty \tilde{\eta}\le \eta$ and $\underset{1\le n \le }{\min} \sum_{j=n}^N(\tilde{\eta}_j-\xi_j)=0$, we can assume without loss of generality that $\alpha = 0$.
 Let $m$ be an integer for which the minimum $\alpha=0$ is attained, i.e.,  $ \sum_{j=m}^N(\eta_j-\xi_j)=0$. Then $m>1$ since  by assumption $\sum_{j=1}^N(\eta_j-\xi_j)\ne 0$.  
For every $1\le n \le m-1$ we have $
 \sum_{j=n}^{m-1}(\eta_j-\xi_j) = \sum_{j=n}^N(\eta_j-\xi_j)\ge 0$,  i.e., $ <\xi_1, \xi_2, \ldots, \xi_{m-1}> \prec_\infty<\eta_1, \eta_2, \ldots, \eta_{m-1}>.$
By the induction hypothesis there is a  monotone nonincreasing sequence $<\zeta_1, \zeta_2, \ldots, \zeta_{m-1}>$  such that
\[
<\xi_1, \xi_2, \ldots, \xi_{m-1}> \;\preccurlyeq_\infty \; <\zeta_1, \zeta_2, \ldots, \zeta_{m-1}> \;\le \; <\eta_1, \eta_2, \ldots, \eta_{m-1}>.
\] 
In particular,  $\sum_{j=n}^{m-1}(\zeta_j-\xi_j) \ge 0$ for $1\le n \le m-1$, with equality for $n=1$. Now define 
\[
\zeta:=\,<\zeta_1, \zeta_2, \ldots, \zeta_{m-1}, \eta_m, \eta_{m+1}, \ldots,  \eta_N>. 
\]
Then $\zeta \le \eta$,
$
  \sum_{j=n}^N(\zeta_j-\xi_j) = 
  \begin{cases}
   \sum_{j=n}^{m-1}(\zeta_j-\xi_j) \ge 0 \quad &\text{for $1\le n \le m-1$}\\
   \sum_{j=n}^N (\eta_j-\xi_j) \ge 0 & \text{for $m\le n \le N$}
     \end{cases},
$
and thus $\xi  \prec_\infty \zeta$.  Furthermore,  $ \sum_{j=1}^N(\zeta_j-\xi_j)=  \sum_{j=1}^{m-1}(\zeta_j-\xi_j) = 0$ and so $\xi  \preccurlyeq_\infty \zeta$.
 It remains to prove that $\zeta$ is nonincreasing, which requires only proving that $\zeta_{m-1} \ge \zeta_m$.  Since  $  \sum_{j=n}^{m-1}(\zeta_j-\xi_j)\ge 0$ for all $n \le m-1$, it follows that $\zeta_{m-1}\ge \xi_{m-1} \ge \xi_m$. Since $  \sum_{j=m}^{}(\zeta_j-\xi_j)= 0$ while $  \sum_{j=m+1}^{}(\zeta_j-\xi_j)\ge 0$ it follows  that $\xi_m \ge \zeta_m$, whence $\zeta_{m-1}\ge \zeta_m$, which completes the proof.
 \ep

 \bP{P:7.2}
 If $\xi, \eta \in \co*$, then the following conditions are equivalent.
 \item[(i)] There is a strictly increasing sequence of indices $\mathbb N \ni n_k \uparrow \infty$ with $n_o=0$ for which
 $\xi^{(n_k)} \prec \eta^{(n_k)}.$
  \item[(ii)] $\xi \prec \eta$ and $\{ \sum_1^n(\eta_j-\xi_j)\} _{ n \ge m}$ attains a minimum for every $m\in \mathbb N$.
   \item[(iii)]  $\xi\prec_b  \zeta \le \eta$ for some $\zeta \in \co*$.
 \item[(iii$'$)]  $\xi \le \rho\prec_b  \eta$ for some $\rho \in \co*$.
  \eP
 \bp
\item[(i)] $\Leftrightarrow $ (ii)
$\xi^{(m)}\prec \eta^{(m)}$ means by definition that $\sum_{j=m+1}^n(\eta_j -\xi_j)\ge 0$   for every $n \ge m+1$, hence   \linebreak $\sum_{j=1}^n (\eta_j -\xi_j) \ge \sum_{j=1}^m (\eta_j -\xi_j) $ for every $n \ge m$, i.e, 
$
\sum_{j=1}^m (\eta_j -\xi_j)= \min_{ n \ge m} \sum_1^n(\eta_j-\xi_j). $
Thus if  $\xi^{(n_k)} \prec \eta^{(n_k)}$ for a strictly increasing sequence of indices $n_k$ with $n_o=0$, then for every $n_{k-1}< m\le n_k$ we have
\[ 
\inf_{ n \ge m} \sum_1^n(\eta_j-\xi_j)= \min\{\min_{m\le n< n_k}  \sum_1^n(\eta_j-\xi_j),  \inf_{n\ge n_k}\sum_1^n(\eta_j-\xi_j)\} \\
= \min_{m\le n\le n_k}  \sum_1^n(\eta_j-\xi_j)
\]
and hence the minimum is attained. Conversely, if (ii) holds, then set $n_o=0$ and construct $n_k$ recursively by setting $n_{k+1}$ to be an index for which the minimum of $\{ \sum_1^n(\eta_j-\xi_j)\}_{ n > n_k}$ is attained. Then by the above equivalence,  $\xi^{(n_k)} \prec \eta^{(n_k)}$ for all $k>0$, while the same relation holds by hypothesis for $k=0$.
\item[(iii)] $\Rightarrow $ (i) Let  $\sum_{j=1}^{n_k} \xi_j = \sum_{j=1}^{n_k} \zeta_j$ for a sequence of indices $\mathbb N \ni n_k \uparrow \infty$ and set $n_0=0$. Then   for every $m > n_k$,
 \[
  \sum_{j={n_k+1}}^m \xi_j =   \sum_{j=1}^m \xi_j - \sum_{j=1}^{n_k} \xi_j \le \sum_{j=1}^m \zeta_j - \sum_{j=1}^{n_k} \zeta_j = \sum_{j={n_k+1}}^m \zeta_j
  \]
  i.e.,  $\xi^{(n_k)} \prec \zeta^{(n_k)}.$ But $\zeta^{(n_k)}\le \eta^{(n_k)}$, hence $\xi^{(n_k)} \prec \eta^{(n_k)}$ for all $k\ge1$. The same relation holds for $k=0$ and $n_o=0$ as $\xi^{(0)} = \xi \prec \eta = \eta^{(0)}$. 
\item[(iii$'$)] $\Leftrightarrow $ (i) Similar argument.
\item[(ii)] $\Rightarrow $ (iii) If $\eta$ is finitely supported, condition (ii) implies that so is $\xi$. If $\xi$ is finitely supported and $\xi_N=0$, then $<\xi_1, \ldots, \xi_{N-1}>\prec <\eta_1, \ldots, \eta_{N-1}>$. In either case the result follows from the finite dimensional one by completing with infinitely many zeros the finite sequence $\zeta$ provided by  Proposition \ref {P:7.1} (A)(i). So assume henceforth that both $\xi_n>0$ and $\eta_n>0$ for all $n$. 
For every $t \in [0,\eta_1] $, define the sequence $\co*\ni \eta(t):=<\min \{t, \eta_n\}>$. 
Then $\eta(t) \le \eta$, $\eta (t)$ is monotone decreasing in $t$, $\xi\prec \eta(\eta_1)=\eta$, and $\xi\not\prec \eta (t)$ for $t<\xi_1$, since $\xi_1 \le \eta_1$. Thus it is easy to see that 
\[
\{t \in [0,\eta_1] \mid \xi \prec   \eta(t)\} = \bigcap_{n=1}^\infty \{\xi_1\le t\le \eta_1\mid \sum_{j=1}^n (\min \{t, \eta_j\}-\xi_j)\ge0\} =[t_1, \eta_1] 
\]
for some $t_1 \in [\xi_1,\eta_1]$. Let $N_1$ be the index for which $\eta_{N_1+1} < t_1 \le \eta_{N_1}$. Then
\be{e:34}
 \sum_{j=1}^n(\eta(t_1)_j-\xi_j) =\begin{cases} \sum_{j=1}^n (t_1- \xi_j) &\quad \text{for } 1\le n \le N_1\\
 \sum_{j=1}^n(\eta_j- \xi_j) - \sum_{j=1}^{N_1}(\eta_j- t_1) &\quad \text{for } n> N_1.
 \end{cases}.
 \ee

But then, $\eta(t)_j= \eta(t_1)_j$ for all $t\in  (\eta_{N_1+1}, t_1]$ and all $j\ge N_1$,  hence the same computation yields for all  $n$
\[
 \sum_{j=1}^n(\eta(t)_j-\xi_j) =  \sum_{j=1}^n(\eta(t_1)_j-\xi_j) - \min\{n,N_1\}(t_1-t).
 \]
By the minimality of $t_1$, it follows that 
$
\inf_n \sum_{j=1}^n(\eta(t_1)_j-\xi_j)=0.
$
Moreover, we see from ( \ref {e:34}) that $\{\sum_1^n(\eta(t_1)_j-\xi_j)\}_{n \ge m} $ attains a minimum for every $m\in \mathbb N$. In particular $ \min_n \sum_{j=1}^n(\eta(t_1)_j-\xi_j)=0$, i.e.,  
 $\sum_{j=1}^{n_1}\xi_j=\sum_{j=1}^{n_1}\eta(t_1)_j$ for some integer  $n_1\ge1$. But then, for every $n > n_1$, we have  \linebreak $\sum_{j=n_1+1}^{n}(\eta(t_1)_j-\xi_j) = \sum_{1}^{n}(\eta(t_1)_j-\xi_j)\ge 0$,
and hence  $\xi^{(n_1)}\prec \eta(t_1)^{(n_1)} \le \eta^{(n_1)}$. Furthermore, by (\ref {e:34}) and the following remark, the sequences $\xi^{(n_1)}$ and $\eta(t_1)^{(n_1)}$ satisfy the condition (ii). Thus by repeating the construction we find 
\[
\xi^{(n_1)}_1\le t_2\le  \eta(t_1)^{(n_1)}_1=  \eta(t_1)_{n_1+1}= \min\{t_1, \eta_{n_1+1}\}
\]
for which $\xi^{(n_1)}\prec \eta(t_1)^{(n_1)}(t_2) = \eta(t_2)^{(n_1)}$ and 
\[\sum_{j=n_1+1}^{n_2}\xi_j=  \sum_{j=1}^{n_2-n_1}\xi^{(n_1)}_j= \sum_{j=1}^{n_2-n_1} \eta(t_2)^{(n_1)}_j= \sum_{j=n_1+1}^{n_2}\eta(t_2)_j
\]
for some integer $n_2> n_1$. Then also $\xi^{(n_2)}$ and $\eta(t_2)^{(n_2)}$ satisfy condition (ii).
Thus iterating the construction, we find a strictly increasing sequence of integers $n_k$ and a nonincreasing sequence $t_k>0$ for which  \linebreak $\xi^{ (n_k) }\prec \eta(t_k)^{ (n_k) }\le \eta ^{ (n_k) }$ and  $ \sum_{j=n_{k-1}+1}^{n_k}\xi_j =\sum_{j=n_{k-1}+1}^{n_k}\eta(t_k)_j$. Let  $\zeta$ be the sequence spliced together from the sequences $\eta(t_k)$, i.e.,   $\zeta_j= \eta(t_k)_j$ for $n_{k-1} < j\le n_k$ and $k\in \mathbb N$. By construction $\zeta\le \eta $ and  $\zeta$ is  nonincreasing since  each $\eta(t_k)$ is nonincreasing and since   $\eta(t_k)\le \eta(t_{k-1})$ for every $k$.  For every $n_{k-1} < n \le n_k$ and $k>1$ (the case $k=1$ being simpler)
\[
\sum_{j=1}^n(\zeta_j-\xi_j)= \sum_{h=1}^{k-1}\sum_{j=n_{h-1}+1}^{n_h}(\eta(t_h)_j-\xi_j) + \sum_{j=n_{k-1}+1}^n(\eta(t_k)_j-\xi_j) = \sum_{j=1}^{n-n_{k-1}}(\eta(t_k)^{(n_{k-1})}_j-\xi^{(n_{k-1})}_j) \ge 0
\]
with equality holding for $n=n_k$. Thus  $\xi\prec_b\zeta$, which completes the proof of (iii).
\item (i) $\Rightarrow $ (iii$'$)   Applying  Proposition \ref{P:7.1}(A)(ii) to each pair of finite sequence
\[
<\xi_{n_k+1}, \xi_{n_k+2}, \ldots, \xi_{n_{k+1}}>\, \prec \,<\eta_{n_k+1}, \eta_{n_k+2}, \ldots, \eta_{n_{k+1}}>
\]
we find  a  nonincreasing finite sequence $<\rho_{n_k+1}, \rho_{n_k+2}, \ldots, \rho_{n_{k+1}}> 
$ for which
 \[
<\xi_{n_k+1}, \xi_{n_k+2}, \ldots, \xi _{n_{k+1}}> \,\, \le \,\,  <\rho_{n_k+1}, \rho_{n_k+2}, \ldots, \rho_{n_{k+1}}> \preccurlyeq <\eta_{n_k+1}, \eta_{n_k+2}, \ldots, \eta_{n_{k+1}}>.
 \]
Let $\rho$ be the sequence spliced together from the finite sequences thus obtained. By construction,  $ \xi \le \rho $  \linebreak and it is immediate to see that  $\rho\prec_b  \eta$.  Since
$
\rho_{n_k} - \eta_{n_k} =  \sum_{j=n_{k-1}+1}^{n_k-1}(\eta_j-\rho_j)\ge 0
$
and hence  \linebreak
$
\rho_{n_k}  \ge \eta_{n_k} \ge \eta_{n_k+1}  \ge \rho_{n_k+1},
$
we see that $\rho$ is monotone nonincreasing, which concludes the proof
\ep

\bR{R:7.3}
\item [(a)] A sufficient condition for (ii) in Proposition \ref {P:7.2} is that $\xi\prec \eta$ and $\sum_{j=1}^{\infty}(\eta_j-\xi_j)=\infty$. 
\item [(b)] The reduction to the finite case that we used to prove the implication (i) $\Rightarrow $ (iii$'$) fails for the  implication (ii) $\Rightarrow $ (iii). Indeed the  sequence created by splicing together the finite sequences  obtained from Proposition \ref {P:7.1}(A)(i) is in general  not monotone nonincreasing.  In other words, to prove (iii) one has to consider the global  behavior of the sequences and not only the behavior  on each of the finite intervals $n_k<j\le n_{k+1}$. 
\eR

As consequence of Proposition \ref{P:7.2} we can extend to infinite sequences Proposition \ref{P:7.1}(A).

\bT{T:7.4}
Let $\xi, \eta \in \co*$ and $\xi \prec  \eta$.
\item [(i)] There is  a  $\zeta \in \co*$ for which $\xi\preccurlyeq \zeta \le \eta$.
\item [(ii)] There is a  $\rho \in \co*$ for which $\xi \le \rho\preccurlyeq \eta$.
\eT
\bp Let $\alpha:={\varliminf} \sum_{j=1}^n ( \eta_j -\xi_j)$.  If $\alpha = 0$, for (i) take $\zeta=\eta$ and for (ii) take $\rho=\xi$. If the conditions of Proposition \ref {P:7.2} are satisfied, e.g., if $\alpha = \infty$ (see Remark \ref {R:7.3}(a)), since block majorization implies strong majorization it is enough to take the sequences $\zeta, \rho \in \co*$ for which $\xi \prec_b \zeta \le \eta$ and $\xi \le \rho\prec_b \eta$  provided by Proposition \ref {P:7.2} (iii) and (iii$'$). Thus assume henceforth that  $0 < \alpha < \infty$ and  that the condition (ii) of Proposition \ref {P:7.2} fails, i.e., that there is some $m\ge 0$ for which $\min_{n>m}  \sum_{j=1}^n ( \eta_j -\xi_j)$ does not exist and hence  $\alpha =\inf _{n> m}\sum_{j=1}^n ( \eta_j -\xi_j)$.

\item[(i)]  If  $\xi \in \ell^*_1$ then either $\xi \preccurlyeq  \eta$, in which case choose $\zeta=\eta$, or  there is a  smallest integer $N$ for which $\sum_{j=1}^N\eta_j-\sum_{j=1}^\infty \xi_j\ge 0$, in which case choose 
\[
\zeta:=\,<\eta_1, \eta_2, \ldots, \eta_{N-1}, \eta_N-\sum_{j=1}^N\eta_j+\sum_{j=1}^\infty \xi_j, 0\ldots>.
\]

Thus assume henceforth that $\xi \not \in \ell^*_1$ and in particular, that both $\xi$ and $\eta$ have infinite support. We proceed as in the proof of the implication of (ii) $\Rightarrow $ (iii) in Proposition \ref {P:7.2}  and construct the nonincreasing sequence $\eta(t_1) :=\min \{t_1,\eta\}$ with $\xi_1\le t_1\le \eta_1$, for which $\xi \prec\eta(t_1) \le \eta$ and $\inf _n\sum_{j=1}^n(\eta(t_1)_j-\xi_j) = 0$.  
If $\varliminf_n\sum_{j=1}^n(\eta(t_1)_j-\xi_j) = 0$, i.e., $\xi\preccurlyeq \eta(t_1)$, then we are done.
If  $\varliminf_n\sum_{j=1}^n(\eta(t_1)_j-\xi_j)> 0$, then there is an integer $n_1$ for which $\sum_{j=1}^{n_1}(\eta(t_1)_j-\xi_j) = 0$. As a consequence, $\xi^{(n_1)}\prec \eta(t_1)^{(n_1)} $ and we can repeat  the construction. We claim that the construction terminates after a finite number of steps. Reasoning by contradiction, assume that the construction can be carried on infinitely many times and thus provides a strictly increasing sequence of integers $n_k$ and a nonincreasing sequence $t_k>0$ starting with $n_o=0$ and $t_o=\eta_1$ for which  
\be{e:35}
\xi^{ ( n_{k-1} ) }\prec (\eta(t_{k-1})^{ ( n_{k-1} ) } )(t_k)= \eta(t_k)^{ ( n_{k-1} ) }
\ee
and 
\be{e:36}
  \sum_{j=1}^{n_k-n_{k-1} }((\eta(t_k)^{(n_{k-1})})_j-(\xi^{(n_{k-1})})_j)=\sum_{j={n_{k-1}+1}}^{n_k}(\eta(t_k)_j-\xi_j)=0.
\ee
Assume furthermore that $n_k> n_{k-1}$ is the smallest integer for which (\ref{e:36}) holds and let $N_k$ be the integer for which   $\eta_{N_k+1}< t_k \le \eta_{N_k}$. Such an integer exists and satisfies the condition $N_k \ge n_{k-1}+1$ since 
\be{e:37}
\eta_{n_{k-1}+1} = ( \eta^{(n_{k-1})})_1 \ge (\eta(t_{k-1})^{(n_{k-1})})_1 \ge t_k =( \eta(t_k)^{(n_{k-1})})_1 \ge  (\xi^{(n_{k-1})})_1= \xi_{n_{k-1}+1} >0
\ee
 and $\eta_j>0$ for all $j$ by the assumption that both $\xi$ and $\eta$ have infinite support.  Moreover, from (\ref{e:34}), we have
 \[
\sum_{j=n_{k-1}+1}^n(\eta(t_k)_j -\xi_j )=
\begin{cases}
\sum_{j=n_{k-1}+1}^n (t_k-\xi_j )\quad & \text{for } n_{k-1}+1\le n \le N_k  \\
\sum_{j=n_{k-1}+1}^n(\eta_j- \xi_j)   - \sum_{j=n_{k-1}+1}^{N_k}(\eta_j- t_k)  &\text{for }  n >  N_k.
\end{cases}
\]

If $n_k > N_k$, then 
\begin{align*}
0=\sum_{j=n_{k-1}+1}^{n_k} (\eta(t_k)_j-\xi_j )&= \sum_{j=n_{k-1}+1}^{n_k} (\eta_j- \xi_j)   - \sum_{j=n_{k-1}+1}^{N_k}(\eta_j- t_k)  \\
&=  \sum_{j=1}^{n_k} (\eta_j- \xi_j) - \sum_{j=1}^{n_{k-1}} (\eta_j- \xi_j)  - \sum_{j=n_{k-1}+1}^{N_k}(\eta_j- t_k)  
\end{align*}
and thus $\{ \sum_{j=1}^{n} (\eta_j- \xi_j) \}_{n>n_{k-1}}$ attains its minimum for $n=n_k$. 
Recall that $m$ is an integer for which  $\{\sum_{j=1}^n ( \eta_j -\xi_j)\}_ {n> m} $ does not  attain a minimum and  let $K$ be an integer for which $n_{k-1}\ge m$. Then  we conclude that $n_k\le N_k$ for every $k\ge K$ and hence
\[
\sum_{j=n_{k-1}+1}^{n_k} \eta(t_k)_j-\xi_j )= \sum_{j=n_{k-1}+1}^{n_k}  (t_k-\xi_j )=0.
\]
It follows that $t_k= \xi_{n_{k-1}+1}=\ldots= \xi_{n_k}$ and since for every $n_{k-1}+1\le n\le n_k$ we also have  \linebreak $\sum_{j=n_{k-1}+1}^{n}( \eta(t_k)_j-\xi_j )= \sum_{j=n_{k-1}+1}^{n}  (t_k-\xi_j )=0$, by the minimality of $n_k$, it follows that $n_k= n_{k-1}+1$. From this and (\ref {e:37}) we have  that $\eta_j\ge \xi_j $ for every $j\ge n_{K-1}+1$. But then  $\sum_{j=n_{K-1}+1}^n(\eta_j-\xi_j)$ is monotone nondecreasing and hence attains a minimum, and therefore $\sum_{j=m}^n(\eta_j-\xi_j)$  also has a minimum,  a contradiction.

Thus we conclude that the construction terminates after $h\ge 1$ steps, i.e., (\ref {e:35}) holds for $1\le k \le h$ and (\ref{e:36}) holds for $1\le k \le h-1$ while  
\[
0=\inf_{n>n_{h-1}}\sum_{j=n_{h-1}+1}^n(\eta(t_h)_j -\xi_j )=\varliminf_n\sum_{j=n_{h-1}+1}^n(\eta(t_h)_j -\xi_j )
\]
and hence $\xi(t_h)^{(n_{h-1})}\preccurlyeq \eta(t_h)^{(n_{h-1})}$.
 Let $\zeta$  be the sequence obtained by splicing together the sequences $\eta(t_k)^{(n_{k-1})}$, i.e.,   
 \[
 \zeta_j= \begin{cases} \eta(t_k)_j \quad &\text{for } n_{k-1} < j\le n_k \quad \text{and } 1\le k\le h-1\\
 \eta(t_h)_j &\text{for } n_{h-1}< j
 \end{cases}
 \]
 Then  it  is easy to verify  that $\xi\preccurlyeq \zeta$ and as in the proof of (ii) $\Rightarrow $ (iii) in Proposition \ref {P:7.2} we see that  $\zeta\le \eta $ and  $\zeta$ is  nonincreasing since  each $\eta(t_k)$ is nonincreasing and since   $\eta(t_k)\le \eta(t_{k-1})$ for every $k$.  
 
\item[(ii)] Let $N$ the smallest  integer $m$ for which the minimum of $\{ \sum_{j=1}^n(\eta_j-\xi_j)\}_{n > m}$ does not exist.  If $N=0$, then $\alpha <   \sum_{j=1}^n(\eta_j-\xi_j)$ for all $n$ and hence $\rho:=\,< \xi_1+ \alpha , \xi_2, \xi_3, \ldots>$ satisfies the required condition.
Assume henceforth that $N\ge 1$. Then 
\be{e:38}
\sum_{j=1}^N(\eta_j-\xi_j) \le \alpha < \sum_{j=1}^n(\eta_j-\xi_j)\quad \text{ for all } n > N,
\ee
By Proposition \ref {P:7.1}(A)(ii) applied to the sequences 
$
<\xi_1, \xi_2, \ldots, \xi_N>\, \prec \,<\eta_1, \eta_2, \ldots, \eta_N>,
$
we can find $<\rho_1, \rho_2, \ldots, \rho_N> \in \mathbb (R^N)^+$ with 
\be{e:39}
<\xi_1, \xi_2, \ldots, \xi_N>\,\le \,<\rho_1, \rho_2, \ldots, \rho_N>  \,\preccurlyeq\,<\eta_1, \eta_2, \ldots, \eta_N>.
\ee
Set 
$\rho:=\,<\rho_1, \rho_2, \ldots, \rho_N, \xi_{N+1}+ \alpha- \sum_{j=1}^N(\eta_j-\xi_j) , \xi_{N+2}, \ldots>$. 
Then $\xi \le \rho$ by (\ref {e:38}) and (\ref {e:39}). Moreover, 
 $\sum_{j=1}^n(\eta_j-\rho_j) \ge 0$  for all  $n \le N$ by (\ref{e:39}), while for $n > N$, 
\begin{align*}
\sum_{j=1}^n(\eta_j-\rho_j) &= \sum_{j=1}^{}(\eta_j-\rho_j)+\eta_{N+1}-\Big(\xi_{N+1}+ \alpha -\sum_{j=1}^{}(\eta_j-\xi_j)\Big) + \sum_{j=N+2}^n(\eta_j-\xi_j)\\
&= \sum_{j=1}^{n}(\eta_j-\xi_j)-\alpha >0
\end{align*}
by (\ref {e:39}) and (\ref {e:38}).
Thus $ \rho \preccurlyeq \eta$.  Finally,  $\rho$  is monotone nonincreasing because

\begin{alignat*}{3}
\rho_N& = \sum_{j=1}^N\rho_j - \sum_{j=1}^{N-1}\rho_j \ge \sum_{j=1}^N\eta_j - \sum_{j=1}^{N-1}\eta_j =\eta_N &\qquad \text{(by (\ref {e:39}))}\\
&\ge \eta_{N+1}> \xi_{N+1}+ \alpha - \sum_{j=1}^N(\eta_j-\xi_j)= \rho_{N+1}&\qquad \text{(by (\ref {e:38}))}
\end{alignat*} 
\ep 

The proof of the case $\xi\in (\ell^1)^*$ in (i) is essentially Mirski's proof \cite {Mi60} of Proposition \ref {P:7.1} (A)(i) (see also \cite[5.A.9]{MO79}).
\bR{R:7.5} 
Let $\xi, \eta,\rho, \zeta \in \co*$ with   $\xi \prec \zeta \le \eta$,  
$\xi \le \rho \prec \eta$, and $\xi\preccurlyeq  \eta$. Then  $\zeta = \eta$  and  $\xi = \rho$. Indeed, for every $m\ge n$, $\eta_n-\zeta_n \le \sum _{j=1}^m(\eta_j-\zeta_j) \le \sum _{j=1}^m(\eta_j-\xi_j)$
and hence $\eta_n-\zeta_n \le {\varliminf}\sum_{j=1}^m ( \eta_j -\xi_j) =0$ for  every $n$. A similar argument shows that $\rho_n-\xi_n \le {\varliminf}\sum_{j=1}^m ( \eta_j -\xi_j) =0$ for every $n$.
\eR

When $  \sum_{j=1}^\infty  \xi_j = \sum_{j=1}^\infty\eta_j < \infty$, then by (\ref{e:3}), $ \eta \prec_\infty\xi \Leftrightarrow \xi \prec \eta$ and so the results in Proposition \ref {P:7.2} and Theorem \ref{T:7.4} can in this case be reformulated in terms of majorization at infinity. However, a more general extension of these results to both summable and nonsummable sequences is given in the next proposition.

\bP{P:7.6} Let $\xi, \eta \in \co*$. Then the following conditions are equivalent.
 \item[(i)] There is a sequence of indices $\mathbb N \ni n_k \uparrow \infty$ with $n_o=0$ for which
 $\sum_{j=n}^{n_k}\xi_j \le  \sum_{j=n}^{n_k}\eta_j$ for all $n_{k-1} < n \le n_k$ and all $k\in \mathbb N$.
 \item[(ii)] Either
 \item[(a)]  $\eta^{(n)} \not \prec \xi^{(n)}$ for  $n=0, 1, 2\ldots$\\
 or 
 \item[(b)] there is an $N\in \mathbb N$ such that $\eta^{(N)} \prec_b  \xi^{(N)}$ and $\sum_{j=n}^N(\eta_j-\xi_j) \ge 0$ for all $1 \le n \le N$. 
\item [(iii)]  $\zeta \le \eta$ and $\zeta\prec_b  \xi$ for some  $\zeta \in \co*$.
\item [(iii$'$)] $\eta\prec_b  \rho$ and $\xi\le \rho$ for some $\rho \in \co*$ .
\eP
\bp
\item (i) $\Rightarrow$ (ii)
If there are infinitely many indices $k$ for which $\sum _{n_{k-1}+1}^{n_k}(\eta_j-\xi_j) >0$, then it is easy to see that (ii)(a) holds. If on the other hand  $\sum _{n_{k-1}+1}^{n_k}(\eta_j-\xi_j) =0$ for every $k> K$ for some $K\in \mathbb N$, then for every $n> n_K$,
\[
\sum _{n_K+1}^n(\xi_j-\eta_j) =\begin{cases} \sum _{n+1}^{n_k}(\eta_j-\xi_j) \ge0 \quad  & n_{k-1}<n<n_k\\
0 & n=n_k
\end{cases},
\]
and thus $\eta^{(n_K)} \prec_b  \xi^{(n_K)}$.  Furthermore, it is clear that $\sum_{j=n}^{n_k}(\eta_j-\xi_j) \ge 0$ for all $0< n \le n_k$ and all $k$ and hence, in particular, for $k=K$, and thus (ii) (b) holds.
\item(ii) (a) $\Rightarrow$ (i) Since $\eta \not \prec \xi $, there is some $n$ for which $\sum_{j=1}^n(\eta_j -\xi_j) >0$. Let $n_1$ be the smallest such integer. Then $
\sum_{j=n}^{n_1}(\eta_j -\xi_j) = \sum_{j=1}^{n_1}(\eta_j -\xi_j) - \sum_{j=1}^{n-1}(\eta_j -\xi_j) >0
$ for every $1< n \le n_1$,
since $ \sum_{j=1}^{n-1}(\eta_j -\xi_j)\le 0$ and the same inequality holds by definition if $n=1$. Iterating the construction we obtain the sequence $n_k$.
\item (ii) (b) $\Rightarrow$ (i) Obvious.
\item(iii) $\Rightarrow$ (i) and (iii$'$) $\Rightarrow$ (i) Obvious.
\item(i) $\Rightarrow$ (iii)
By applying Proposition \ref {P:7.1}(B)(i) to each pair of finite sequences
\[
<\xi_{n_{k-1}+1}, \xi_{n_{k-1}+2},\ldots, \xi_{n_k}> \: \prec_\infty\:  <\eta_{n_{k-1}+1}, \eta_{n_{k-1}+2},\ldots, \eta_{n_k}> 
\]
we find  finite monotone nonincreasing sequences 
\[
<\xi_{n_{k-1}+1}, \xi_{n_{k-1}+2},\ldots, \xi_{n_k}> \: \preccurlyeq_\infty \: 
<\zeta_{n_{k-1}+1}, \zeta_{n_{k-1}+2},\ldots, \zeta_{n_k}> \: \le  \: 
<\eta_{n_{k-1}+1}, \eta_{n_{k-1}+2},\ldots, \eta_{n_k}>. 
\]
Thus  $<\zeta_{n_{k-1}+1}, \zeta_{n_{k-1}+2},\ldots, \zeta_{n_k}>\: \preccurlyeq \: <\xi_{n_{k-1}+1}, \xi_{n_{k-1}+2},\ldots, \xi_{n_k}>$ for every $k$.  Let $\zeta$ be the sequence obtained by splicing together the finite sequences thus obtained for each interval  $n_{k-1} < n \le n_k$. Then  $\zeta\le \eta$ and if $ n_{k-1} < n \le n_k$ for $k>1$, then
\[
 \sum_{j=1}^{n} \zeta_j = \sum_{i=1}^{k-1}\Big(\sum_{j= n_{i-1}+1}^{n_i}\zeta_j  \Big) + \sum_{j= n_{k-1}+1}^n\zeta_j \le  \sum_{i=1}^{k-1}\Big(\sum_{j= n_{i-1}+1}^{n_i}\xi_j  \Big) + \sum_{j= n_{k-1}+1}^n\xi_j =  \sum_{j=1}^{n} \xi_j 
\]
with equality holding for $n = n_k$. The same conclusion holds directly for $k=1$. Thus,  $\zeta\prec_b  \xi$. Finally,
\[
\zeta_{n_k +1} = \sum_{j=1}^{n_k +1}\zeta_j- \sum_{j=1}^{n_k }\zeta_j\le \sum_{j=1}^{n_k +1}\xi_j- \sum_{j=1}^{n_k }\xi_j = \xi_{n_k +1}\le \xi_{n_k}\le \zeta_{n_k},
\]
 which proves that $\zeta$ is monotone nonincreasing.
 \item (ii) (b) $\Rightarrow$ (iii$'$) 
Let $\rho:=\,< \xi_1+\sum_{j=1}^N(\eta_j-\xi_j), \xi_2, \xi_3, \ldots>$. Then $\xi \le \rho \in \co*$ and
\[
 \sum_{j=1}^n(\rho_j - \eta_j)=  \sum_{j=1}^n(\xi_j - \eta_j) +\sum_{j=1}^N(\eta_j-\xi_j)= 
 \begin{cases}
 \sum_{j=n+1}^N(\eta_j-\xi_j)\ge 0 \quad &\text{for } 1\le n < N\\
 0&\text{for } n = N\\
  \sum_{j=N+1}^n(\xi_j - \eta_j)\ge 0 &\text{for } n > N.
  \end{cases}
  \]
Thus it is now immediate to see that $\eta \prec_b  \rho$.

The remaining implication is the crux of the proof.
\item (ii) (a) $\Rightarrow$ (iii$'$)
For all $t\ge0$ and $p\in \mathbb N$ set
 \[
\rho (t,p) =
 \begin{cases}
 <t, t, \ldots> & \text {if $p=1$}\\
 <\xi_1, \xi_2, \ldots, \xi_{p-1}, t, t, \ldots> & \text {if $p>1.$}
 \end{cases}
 \]
Since $\eta = \eta^{(o)} \not \prec   \xi^{(o)} =\xi$, there is an $N \in \mathbb N$ for which $\sum_{j=1}^N\eta_j > \sum_{j=1}^N\xi_j = \sum_{j=1}^N\rho(\xi_N,N)_j$, i.e., $\eta \not\prec \rho(\xi_N,N)$. (Notice that although $ \rho (t,p)$  may fail to be monotone, it is convenient and should cause no confusion to apply to it the majorization notations we use for monotone sequences.) Let  \linebreak $p_1:= \min \{ p \mid \eta \not\prec \rho(\xi_p,p)\};$ thus $p_1\le N$.   
Assume first that $p_1 > 1$.  Then $\eta \prec \rho(t,p_1)$ if and only if 
\be{e:40}
\sum_{j=1}^n(\rho(t,p_1)-\eta_j)=
\begin{cases}
\sum_{j=1}^n(\xi_j-\eta_j)\ge 0\quad &\text{for } 0<n < p_1\\
(n-p_1+1)t - \sum_{j=1}^n\eta_j + \sum_{j=1}^{p_1-1}\xi_j \ge 0 &\text{for  } n \ge p_1.
\end{cases}
\ee
Since  $\eta \prec \rho(\xi_{p_1-1}, p_1-1)= \rho(\xi_{p_1-1}, p_1)$ by the minimality of $p_1$,  the inequalities in (\ref{e:40}) are  satisfied  for all $0<n < p_1$ for any choice of $t\ge0$.  For every $n\ge p_1$ set 
 \[
 t(n):= \frac{1}{n-p_1+1}\big( \sum_{j=1}^n\eta_j -\sum_{j=1}^{p_1-1}\xi_j \big).
 \]
Then the inequalities in (\ref{e:40})  are  satisfied for all $n \ge p_1$ if and only if $t\ge \sup_{n\ge p_1} t(n)$. Since \[ t(n)= \frac{n}{n-p_1+1}(\eta_a)_n - \frac{1}{n-p_1+1}\sum_{j=1}^{p_1-1}\xi_j \]  where $(\eta_a)_n$ is the arithmetic mean of $\eta$ and since $(\eta_a)_n\to 0$ it follows that $t(n)\to 0.$ Since $t_N >0$,  there is some $n_1 \ge p_1$ for which  the supremum is attained, i.e., $t_1:=t(n_1) = \max_{n\ge p_1} t(n)$. Thus  $\eta \prec \rho(t_1,p_1)$ and from  (\ref{e:40}) it follows also that  $ \sum_{j=1}^{n_1}(\rho (t_1, p_1)_j-\eta_j)=0$. 
Since  $\eta \prec \rho(\xi_{p_1-1}, p_1-1) = \rho(\xi_{p_1-1}, p_1)$ and  $\eta \not \prec \rho(\xi_{p_1}, p_1) $, it follows from the monotonicity of $\rho(t,p_1)$ in $t$, that  $\xi_{p_1} <  t_1\le \xi_{p_1-1} $. As a consequence,  $\rho(t_1, p_1)$ is nonincreasing and $\xi \le \rho(t_1, p_1)$.    Furthermore, from $1< p_1\le n_1$ it follows that $\rho(t_1,p_1)_1=\xi_1$ and $ \rho(t_1, p_1)_{n_1}=t_1> \xi_{p_1}\ge \xi_{n_1}\ge \xi_{n_1+1}$ and from $\eta \prec \rho(t_1,p_1)$ it follows that $\rho(t_1,p_1)_1\ge  \eta_1$ and 
\[
\rho(t_1, p_1)_{n_1}\ge  \rho(t_1, p_1)_{n_1+1}=  \sum _{j=1}^{n_1+1}\rho(t_1, p_1)_j-   \sum _{j=1}^{n_1}\rho(t_1, p_1)_j
  \ge  \sum _{j=1}^{n_1+1}\eta_j -   \sum _{j=1}^{n_1}\eta_j = \eta_{n_1+1}. 
  \]
Thus
\begin{align}\label{e:41}
&\xi \le \rho(t_1, p_1), \quad \eta \prec \rho(t_1,p_1), \quad  \sum_{j=1}^{n_1}(\rho (t_1, p_1)_j-\eta_j)=0,\\
& \rho(t_1, p_1)_1= \max \{\xi_1, \eta_1\}, \, \text{and  }  \rho(t_1, p_1)_{n_1}\ge \max\{ \xi_{n_1+1},  \eta_{n_1+1} \}.\notag
\end{align}
   
The case when   $p_1 = 1$, i.e., $\eta \not\prec \rho(\xi_1,1)= \,<\xi_1,\xi_1,\ldots >$ occurs if and only if $\eta_1> \xi_1$.  Then set $t_1:=\eta_1$ and $n_1= 1$ and the relation (\ref {e:41}) follow trivially.

Now we can repeat the same construction applying it to the sequences $\eta^{(n_1)}\,\not \prec\,\xi^{(n_1)}$.  Iterating, we obtain a strictly increasing sequence of indices $n_k$ and  for each $k \in \mathbb N $ a finite monotone nonincreasing sequence   
\[
<\rho_{n_{k-1}+1}, \rho_{n_{k-1}+2},\ldots, \rho_{n_k}> \: \ge \: <\xi_{n_{k-1}+1}, \xi_{n_{k-1}+2},\ldots, \xi_{m_h}> 
\]
for which 
$
\sum_{j= n_{k-1}+1}^{n}\rho_j \ge \sum_{j= n_{k-1}+1}^{n}\eta_j  
$
for all $n_{k-1}<n\le n_k$ with equality  for $n=n_k$ and furthermore for which $\rho_{n_{k-1}+1}= \max\{\xi_{n_{k-1}+1}, \eta_{n_{k-1}+1}\}$  and $\rho_{n_k} \ge \max\{\xi_{n_k+1}, \eta_{n_k+1}\}$. This latter condition shows that $\rho_{n_k} \ge \rho_{n_k+1}$ and hence the sequence $\rho$ obtained by splicing together these finite sequences is monotone nonincreasing. It is now immediate to see that it satisfies (iii$'$).
\ep

As consequence of Proposition \ref{P:7.6} we can extend to infinite summable sequences Proposition \ref {P:7.1}(B).

\bT{T:7.7} Let $\xi,\eta \in (\ell^1)^*$ and   $ \xi \prec_\infty \eta$. 
\item [(i)] $\xi \preccurlyeq_\infty \zeta \le \eta$ for some $\zeta \in \co*$.
\item [(ii)] $\xi \le \rho\preccurlyeq_\infty \eta$ for some $\rho \in \co*$.
\eT

\bp
 \item [(i)] If $\eta^{(n)} \not \prec \xi^{(n)}$ for  all $n$, then by Proposition \ref {P:7.6}there is a  $\zeta \in \co*$ for which $\zeta \le \eta$  and $\zeta\prec_b  \xi$, hence $\zeta\preccurlyeq \xi$, and thus $\xi \preccurlyeq_\infty \zeta $. 
Thus it remains to consider the case when $\eta^{(N)}  \prec \xi^{(N)}$ for some $N$, i.e.,  $\sum_{j=N+1}^n(\xi_j-\eta_j) \ge 0$ for every $n \ge N+1$. Thus $\sum_{j=N+1}^\infty(\xi_j-\eta_j) \ge 0$. By hypothesis,  $\sum_{j=N+1}^\infty(\xi_j-\eta_j) \le 0$, hence we have equality, i.e., $\eta^{(N)} \preccurlyeq \xi^{(N)}$, or, equivalently,  $\xi^{(N)}\preccurlyeq_\infty \eta^{(N)}$. In particular, $\xi_{N+1}\ge \eta_{N+1}$. Moreover,  
$
\sum_{j=n}^N(\eta_j-\xi_j)= \sum_{j=n}^\infty(\eta_j-\xi_j)\ge 0$ for every $1\le n\le N$, 
i.e., $<\xi_1,\xi_2, \ldots, \xi_N>  \prec_\infty <\eta_1,\eta_2, \ldots, \eta_N>$. By Proposition \ref {P:7.1}(B)(i), there is a finite nonincreasing sequence $\zeta \in \mathbb (R^N)^+$ with  $\xi \preccurlyeq_\infty \zeta \le \eta$ on the integer interval $1, 2, \ldots, N$. In particular, $\xi_N \le \zeta_N$. Define the infinite sequence $\zeta:=\,<\zeta_1, \zeta_2, \ldots, \zeta_N, \eta_{N+1}, \eta_{N+2}, \ldots>$. 
Since   $\zeta_{} \ge \xi _{N }\ge \xi _{N +1}\ge \eta_{N +1}= \zeta_{N +1}$, it follows that $\zeta$ is monotone nonincreasing. It is now immediate to see that $\xi \preccurlyeq_\infty \zeta \le \eta$.\\
 \item [(ii)]  It suffices to choose $\rho:=\,<\xi_1+\sum_{j=1}^\infty(\eta_j-\xi_j), \xi_2, \xi_3, \ldots >$.
 \ep

\section{Applications to operator ideals} \label{S: 8} 
Calkin \cite{jC41} established a correspondence between the two-sided 
ideals of $B(H)$ for a  separable infinite-dimensional
Hilbert space $H$ and the \textit{characteristic sets}. 
These are the positive cones of $\text{c}_{\text{o}}^*$
 that are hereditary (i.e., solid) and invariant under \textit{ampliations}, i.e., the maps
\[
\text{c}_{\text{o}}^* \owns \xi \rightarrow 
D_m\xi:=~<\xi_1,\ldots,\xi_1,\xi_2,\ldots,\xi_2,\xi_3,\ldots,\xi_3,\ldots>
\]
where each entry $\xi_i$ of $\xi$ is repeated $m$-times.
The order-preserving lattice isomorphism
$I \rightarrow \Sigma(I)$ maps each ideal to its characteristic set 
$\Sigma(I) := \{s(X) \mid X \in I\}$
and its inverse is the map $\Sigma\rightarrow I(\Sigma)$ where
 $I(\Sigma)$ is the ideal generated by $\{\diag \xi \mid \xi \in \Sigma\}$.
 
Two sequence operations, the arithmetic mean restricted to $\co*$ and the
arithmetic mean at infinity restricted to $(\ell^1)^*$ respectively,
\[
\xi_a := \langle\frac{1}{n}\sum_1^n \xi_j\rangle
  \quad \text{and} \quad \xi_{a_\infty} :=
\langle\frac{1}{n}\sum_{j=n+1}^\infty \xi_j\rangle
 \]
are essential for the study of commutator spaces (i.e., commutator ideals) and hence traces on ideals (see \cite{DFWW} and \cite{vKgW02} -\cite{vKgW07-OT21}.) 
If $I$ is an ideal, then the arithmetic mean ideals $_aI$ and $I_a$, 
called the \textit{pre-arithmetic mean} and \textit{arithmetic mean} 
of $I$, are the ideals with characteristic sets
\[
\Sigma(_aI) := \{\xi \in \text{c}_{\text{o}}^* \mid \xi_a \in \Sigma(I)\},
\]
\[
\Sigma(I_a) := \{\xi \in \text{c}_{\text{o}}^* \mid \xi = 
O(\eta_a)~\text{for some}~ \eta \in \Sigma(I)\},
\]
where the notation $\xi=O(\zeta)$ means $\xi_n\le M\zeta_n$ for some $M>0$ and all $n$. The ideal $I^-:=\,_a(I_a)$ is called the  \textit{am-closure} of $I$ and an ideal $I$ is called am-closed if $I = I^-$  (there is an analogous notion of am-interior and am-open ideals). It is easy to see that $I\subset  I^-$ , the map $I\rightarrow  I^-$ is idempotent, and  $\mathscr L_1$ is the am-closure of the finite rank ideal $F$ and hence the smallest  am-closed ideal (see \cite{vKgW04-Traces}.)
Many of the ideals arising from classical sequence spaces are am-closed and am-closed ideals play a central  role in the study of single commutators in operator ideals. 

Similarly, the arithmetic mean at infinity  ideals $_{a_\infty}I$ and $I_{a_\infty}$
  are the ideals with characteristic sets
\[
\Sigma(_{a_\infty}I) := \{\xi \in (\ell^1)^* \mid \xi_{a_\infty} \in
\Sigma(I)\}
\]
\[
\Sigma(I_{a_\infty}) := \{\xi \in \text{c}_{\text{o}}^{*} \mid \xi =
  O(\eta_{a_\infty})  ~\text{for some}~ \eta \in \Sigma(I \cap \mathscr L_1)\}.
\]
The  am-$\infty$ closure of $I$ is
$I^{-\infty} :=\, _{a_\infty}(I_{a_\infty})$ and an ideal $I$ is called am-$\infty$ closed if $I = I^{-\infty}$  (there is an analogous notion of am-$\infty$ interior and am-$\infty$ open ideals). The map $I\rightarrow  I^{-\infty}$ is also idempotent, but by definition, $ I^{-\infty}\subset  \mathscr L_1$ and hence we only have  $I\cap  \mathscr L_1 \subset  I^{-\infty}.$
For the definition and basic properties of arithmetic mean at infinity ideals, see \cite{vKgW04-Traces} and \cite {vKgW04-Soft}. Papers \cite {vKgW02}-\cite {vKgW07-OT21} develop our investigation of the theory of arithmetic mean and arithmetic mean at infinity ideals.

The am-closure of an ideal is naturally reformulated in terms of majorization by elements of the characteristic set of the ideal since by definition $\xi \prec \eta$ if and only if $\xi_a\le \eta_a$ for any $\xi, \eta \in \co*$ (see the equivalence of conditions (i) and (ii) in the following theorem.)

\bT{T:8.1}
Let $I$ be an ideal and let $ \xi \in \co*$. Then the following conditions are equivalent.
\item[(i)] $\xi \in \Sigma (I^-)$.
\item[(ii)]  $\xi \prec \eta$ for some $\eta \in  \Sigma(I)$.
\item[(ii$'$)] $\xi \preccurlyeq \eta$ for some $\eta \in  \Sigma(I)$.
\item[(iii)] $\xi =P \eta$ for some $\eta \in  \Sigma(I)$ and some substochastic matrix $P$.
\item[(iii$'$)]  $\xi =Q \eta$ for some $\eta \in  \Sigma(I)$ and some orthostochastic matrix $Q$.
\item[(iv)]   $\diag \xi = E(L \diag  \eta\1 L^*)$ for some $\eta \in  \Sigma(I)$ and some contraction $L$.
\item[(iv$'$)]  $\diag \xi = E(U \diag  \eta\1 U^*)$ for some $\eta \in  \Sigma(I)$ and some orthogonal matrix $U$.\\
\noindent
If $I\supset \mathscr L_1$  then these conditions are also equivalent to:
\item[(ii$''$)] $\xi\prec_b  \eta$ for some $\eta \in  \Sigma(I)$.
\item[(iii$''$)]  $\xi =Q \eta$ for some $\eta \in  \Sigma(I)$ and some block orthostochastic matrix $Q$.
\item[(iv$''$)]  $\diag \xi = E(U \diag  \eta\1 U^*)$ for some $\eta \in  \Sigma(I)$ and some matrix $U$ direct sum of finite orthogonal matrices.
\eT
\bp
\item[(i)] $ \Leftrightarrow$ (ii) This is  a reformulation of the definition of $I^-$. Indeed, $\xi\in \Sigma(I^-)=\Sigma(_a(I_a)) \Leftrightarrow \xi_a \in \Sigma(I_a)  \Leftrightarrow  \linebreak \xi_a = O(\eta_a)$ for some $\eta\in \Sigma(I)  \Leftrightarrow \xi_a \le \eta_a$ for some $\eta\in \Sigma(I)   \Leftrightarrow \xi \prec \eta$ for some $\eta \in  \Sigma(I)$.
\item[(ii)] $ \Rightarrow$ (ii$'$) Let $\xi \prec \eta$ for some $\eta \in  \Sigma(I)$. By Theorem \ref {T:7.4}
 there is a $\zeta \in \co*$ for which $\xi\preccurlyeq \zeta \le \eta$. By the hereditariness of $\Sigma$, $\zeta \in  \Sigma(I)$. 
\item [(ii$'$)] $\Rightarrow$ (ii)  Obvious.
\item[(ii)]  $ \Leftrightarrow$ (iii) $\Leftarrow$ (iv) By \eqref{e:1} and Lemmas \ref {L:2.3}  and \ref{L:2.4} . 
\item[(ii$'$)]  $ \Rightarrow$ (iii$'$) $ \Leftrightarrow$ (iv$'$) By Theorem \ref {T:3.9} and Lemma  \ref {L:2.3} .
\item [(iii$'$)]$\Rightarrow$ (iii), (iv$'$)$\Rightarrow$ (iv), and (ii$''$) $\Rightarrow$(ii)  Obvious. 
\item[(ii$''$)]  $ \Leftrightarrow$ (iii$''$) $ \Leftrightarrow$ (iv$''$) By (\ref{e:6}) and Lemma \ref {L:2.3}.\\
\noindent Assume now that $I\supset \mathscr L_1$. 
\item [(ii)] $\Rightarrow$ (ii$''$) If $\xi \in(\ell^1)^*$, then $\xi \in \Sigma (I)$ and hence there is nothing to prove. If $\xi \not \in (\ell^1)^*$, then  \linebreak  $\sum_{j=1}^\infty (2\eta_j-\xi_j)=\infty$.  
As by  Remark \ref {R:7.3}(a),  condition (ii) of Proposition \ref{P:7.2} is satisfied for   $\xi \prec 2\eta\in \Sigma(I)$.  Hence there is a $\zeta \in \co*$ for which $\xi\prec_b  \zeta \le 2\eta$ and then $\zeta \in \Sigma (I)$ by the hereditariness of  $\Sigma (I)$.
\ep

\bR{R:8.2} An example of an ideal $I\not \supset  \mathscr L_1$ for which (i) $\not \Rightarrow$ (ii$''$) is the finite rank ideal $F$. It is easy to see that  $F^-=\mathscr L_1$ (e.g., see \cite {vKgW04-Soft}), but it is clear that every sequence $\xi \in\Sigma(\mathscr L_1) =(\ell^1)^*$ that is not finitely supported cannot be block majorized by a finitely supported $\eta$, i.e., (ii$''$) fails.
\eR

The analogue of Theorem \ref {T:8.1} for am-$\infty$ closure is:

\bT{T:8.3}
Let $I\subset \mathscr L_1$ be an ideal and let $ \xi \in (\ell^1)^*$. Then the following conditions are equivalent.
\item[(i)] $\xi \in \Sigma (I^{-\infty}).$
\item[(ii)]  $ \xi \prec_\infty\eta$ for some $\eta \in  \Sigma(I).$
\item[(ii$'$)]  $\xi \preccurlyeq_\infty \eta$ for some $\eta \in  \Sigma(I).$
\item[(ii$''$)] $\eta \preccurlyeq \xi$ for some $\eta \in  \Sigma(I).$
\item[(iii)]  $P\xi \in \Sigma(I)$ for some column-stochastic matrix $P.$
\item[(iii$'$)]  $Q\xi \in \Sigma(I)$ for some block orthostochastic matrix $Q.$
\item[(iv)]   $E(V \diag  \xi\1 V^*)\in I$ for some isometry $V.$
\item[(iv$'$)]   $E(U \diag  \xi\1 U^*)\in I$ for some matrix $U$ direct sum of finite orthogonal matrices.
\eT
\bp
\item[(i)] $ \Leftrightarrow$ (ii) As in the proof of Theorem \ref{T:8.1}, this equivalence is  a reformulation of the definition of $I^{-\infty}$.
\item[(ii)] $ \Leftrightarrow$ (ii$'$) One direction is by Theorem \ref {T:7.7}
  and the same argument as in the proof of Theorem \ref{T:8.1}; the other direction is obvious.
\item[(ii$'$)] $ \Leftrightarrow$ (ii$''$) By definition, e.g., see (\ref{e:3}).
\item[(ii$'$)] $ \Rightarrow$  (iii$'$) If $\eta_n>0$ for all $n$, then $\sum_{j=n}^\infty(2\eta_j-\xi_j)> \sum_{j=n}^\infty\eta_j >0$, hence $2\eta^{(n)} \not\prec \xi^{(n)}$ for every $n$ and thus condition (ii)(a) in Proposition \ref {P:7.6} applies to the pair of sequences $\xi, 2\eta$. Thus there is  a $\zeta \in \co*$ for which  $\zeta\prec_b  \xi$ and $\zeta \le 2\eta$ and hence $\zeta \in \Sigma(I)$. By (\ref{e:6}),  $\zeta = Q\xi$ for some block orthostochastic matrix $Q$, i.e., $Q\xi\in  \Sigma(I)$. If on the other hand $\eta_{}=0$ for some $N\in \mathbb N$, then also  $\xi_{}=0$ and hence by the Horn Theorem there is an $N\times N$ orthostochastic matrix $Q_0$ mapping $<\xi_1, \ldots, \xi_N>$ onto $<\eta_1, \ldots, \eta_N>$. But then, $Q:=Q_o\oplus I$ is also block orthostochastic and $Q\xi=\eta\in \Sigma(I)$.
\item[(iii$'$)] $ \Rightarrow$  (iv$'$)$\Rightarrow$ (iv)$\Rightarrow$(iii)$\Rightarrow$(ii$''$) Obvious from the definitions and Lemmas \ref {L:2.3} ,  \ref{L:2.4} , and \ref{L:2.10}. 
\ep

\bC{C:8.4}
For every ideal $I$, $E(I) = I^-\cap \mathscr D$.
\eC
\bp
Since ideals are the linear span of their positive parts and the conditional expectation is linear, it is enough to prove that $E(I^+) = (I^-)^+\cap \mathscr D$. 
Let $A\in I^+$. By Proposition \ref{P:6.4} , $s(E(A)) \prec s(A)$ and since $s(A)\in \Sigma(I)$ by definition, by Theorem \ref {T:8.1} it follows that 
$s(E(A))  \in \Sigma (I^-)$, i.e., $E(A)\in I^-$. This proves that $E(I) \subset I^- \cap \mathscr D$.
Let  now $B\in I^- \cap \mathscr D$. Then, again by Theorem \ref {T:8.1}, $s(B) \preccurlyeq \eta$ for some $\eta\in \Sigma(I)$. By  Proposition \ref{P:6.4} , $B\in E(\mathscr V(\diag \eta) )\subset E(I)$, which concludes the proof.
\ep

\bC{C:8.5}
An ideal $I$ is am-closed if and only if $E(I) \subset I.$
\eC
\bp
If $I$ is am-closed, i.e., $I=I^-$, then by Corollary \ref{C:8.4}, $E(I) = I\cap \mathscr D\subset I$. 
Conversely, if $E(I) \subset I$ then, again by Corollary \ref {C:8.4}, $I^-\cap \mathscr D\subset I$.
But since every selfadjoint element of $I^-$ is diagonalizable, and since $I$ and $I^-$ are ideals and hence unitarily invariant, 
this implies that $I^- \subset I$ and hence $I^- = I$, i.e., $I$ is am-closed.
\ep

Theorem  \ref {T:8.1} leads naturally to the following notion of invariance of an ideal under the action of a class of substochastic matrices. Recall from Remark \ref {R:2.2} that  if $\xi\in\co*$ and $P$ is a substochastic matrix, then $(P\xi)^*\in \co*$, where * denotes monotone rearrangement.  

\bD{D:8.6} 
Given a collection $\mathscr P$ of substochastic matrices and an ideal $I$, we say that 
$I$ is invariant under  $\mathscr P$ if $(P\xi)^*\in \Sigma(I)$ for every $P\in \mathscr P$  and every $\xi \in \Sigma(I) $.
\eD

An immediate consequence of Theorem \ref {T:8.1} is the following characterization of am-closed ideals in terms of invariance under various classes of substochastic matrices. 

\bC{C:8.7}Let  $I$ be an ideal. Then the following conditions are equivalent.
\item[(i)] $I$ is am-closed.
\item[(ii)] $I$ is invariant under substochastic matrices.
\item[(ii$'$)] $I$ is invariant under one of the following classes of substochastic matrices: row-stochastic, column-stochastic, double stochastic, isometry stochastic,  co-isometry stochastic, unistochastic, orthostochastic.
\item[(iii)] $I\supset \mathscr L_1$ and $I$ is invariant under block stochastic matrices.
\item[(iii$'$)] $I\supset \mathscr L_1$ and $I$ is invariant under block orthostochastic matrices.
\eC

We consider now invariance under other classes of substochastic matrices. For the finite case, Birkhoff \cite {Bg46} proved that doubly stochastic matrices are convex combinations of permutation matrices. For the infinite case, Kendall \cite {Kd60}  proved that convex combinations of permutation matrices are dense in the class of doubly stochastic matrices for an appropriate topology.

It is thus natural (e.g., see  A. Neumann's remark (\cite [pg 448] {Na99}))  to consider  an intermediate class of matrices, the infinite convex combinations of permutation matrices:
\[
\mathscr C:=\{\sum_{j=1}^\infty t_j\Pi_j \mid \Pi_j \text { permutation matrix, } 0< t_j\le1,~  \sum_{j=1}^\infty t_j=1\}.
\] 

Clearly, all Banach ideals are invariant under $\mathscr C$. However, there are Banach ideals that are not am-closed, and hence, by Corollary \ref {C:8.7} are not invariant under doubly stochastic matrices. Indeed Varga proved  that for any principal ideal  $(\xi)$ generated by $\diag \xi$ for some  irregular and nonsummable sequence $\xi$, 
\be{e:42}
 \cl(\xi) \subsetneq (\xi)^-= ( \cl(\xi))^-
 \ee
  where the closure $\cl$ is taken under the principal ideal norm which is well known to be complete (see \cite [Remark 3]{Vj89}  and see also \cite[Remark 4.8] {vKgW04-Soft}).  Therefore invariance under  $\mathscr C$  is strictly weaker than invariance under doubly stochastic matrices. 

However, for \textit{soft-edged} ideals (ideals $I$ for which $\xi\in \Sigma(I)$ if and only if $\xi = \alpha \eta$ for some $\alpha \in \co*$ and some $\eta \in  \Sigma(I)$) and for \textit{soft-complemented} ideals (ideals $I$ for which $\xi\in \Sigma(I)$ if and only if $ \alpha\xi  \in  \Sigma(I)$ for all $\alpha \in \co*$), invariance under $\mathscr C$  is equivalent to invariance under doubly stochastic matrices. For the properties of soft-edged and soft-complemented ideals we refer the reader to \cite {vKgW04-Soft}) (see also \cite {vKgW02}, \cite{vKgW04-Traces}). 

 To prove  this fact, we introduce first the following class of block substochastic  matrices

\[
\mathscr B:=\{\sum_{k=1}^\infty \oplus t_kP_k \mid P_k \text { finite doubly stochastic matrix, } 0< t_k\le1,~  \sum_{k=1}^\infty t_k=1\}.
\]

\bL{L:8.8}  Let $I\ne\{0\}$ be an ideal.
\item[(i)] If $I$ is invariant under $\mathscr C$, then $I\supset \mathscr L_1$ and $I$ is invariant under $\mathscr B$.
\item[(ii)] If $I$ is soft-edged or soft-complemented and  is invariant under $\mathscr B$, then $I$ is invariant under block stochastic matrices. 
\eL

\bp
\item[(i)] Let $t:=\,<t_k> \in (\ell^1)^*$ and without loss of generality, assume that $\sum_{k=1}^\infty t_k=1$. Let $\Pi_k$ be the permutation matrix corresponding to the transposition $1 \leftrightarrow k$ and let \textbf {1}$:=\,<1,0, \ldots>$. Since \textbf {1} $\in \Sigma(I)$,  by the invariance of $I$ under $\mathscr C$ we have $\big(\sum_{j=1}^\infty t_k\Pi_k \textbf {1}\big)^* = t \in  \Sigma(I)$. This proves that  $\Sigma(I)\supset (\ell^1)^* $ and hence $I\supset \mathscr L_1$. 

Now let $P\in \mathscr B$, i.e., $P:=  \sum_{k=1}^\infty \oplus t_kP_k$ where $0< t_k\le1,~  \sum_{k=1}^\infty t_k=1$, $P_k$ are finite doubly stochastic matrices and the k-th direct sum block corresponds to the indices $n_{k-1} < i,j\le n_k$. Recall that by the Birkhoff Theorem \cite {Bg46} (see also \cite {MO79}), each $P_k = \sum_{h=1}^{m_k} s_{h,k} \pi_{h,k}$  is a convex combination of (finite) permutation matrices $\pi_{h,k}$.
Let  $\Pi_{h,k}:=I\oplus \pi_{h,k}\oplus I$ be the infinite permutation matrix corresponding to $\pi_{h,k}$, i.e., the permutation matrix that agrees with $\pi_{h,k}$ for the indices $n_{k-1}< i \le n_k$  and that leaves all the other indices invariant. For every $\xi\in \Sigma (I)$ and every $h$ and $k$, 
\ba
(0\oplus \pi_{h,k}\oplus 0)\xi &=
< 0, 0, \ldots, 0,\pi_{h,k}\big(<\xi_{n_{k-1}+1}, \xi_{n_{k-1}+2}, \ldots,  \xi_{n_k}>\big), 0, 0, \ldots> \\
&\le <\xi_1, \xi_2\ldots, \xi_{n_{k-1}}, \pi_{h,k}\big(<\xi_{n_{k-1}+1}, \xi_{n_{k-1}+2}, \ldots,  \xi_{n_k}>\big),  \xi_{n_k+1}, \ldots >
= \Pi_{h,k}\xi.
\end{align*}
 But then 
\ba
P\xi&=\big(\sum_{k=1}^\infty \oplus t_k\sum_{h=1}^{m_k} s_{h,k} \pi_{h,k}\big)\xi \\
&=
 \sum_{k=1}^\infty < 0, 0, \ldots,\big(\sum_{h=1}^{m_k} t_k s_{h,k} \pi_{h,k}<\xi_{n_{k-1}+1}, \xi_{n_{k-1}+2}, \ldots,  \xi_{n_k}>\big), 0, 0, \ldots> \\
&\le \sum_{k=1}^\infty\sum_{h=1}^{m_k} t_k s_{h,k}   \Pi_{h,k} \xi =  \Big(\sum_{k=1}^\infty\sum_{h=1}^{m_k} t_k s_{h,k}   \Pi_{h,k}\Big) \xi.
\end{align*}
Now,  $\sum_{k=1}^\infty\sum_{h=1}^{m_k} t_k s_{h,k}   \Pi_{h,k}\in \mathscr C$ and hence
\[
\bigg(\big(\sum_{k=1}^\infty \oplus t_k\pi_k \big)\xi\bigg)^*\le \bigg(\big(\sum_{k=1}^\infty\sum_{h=1}^{m_k} t_k s_{h,k}   \Pi_{h,k}\big) \xi \bigg)^*\in \Sigma (I),
\]
which proves that  $I$ is invariant under $\mathscr B$.

\item[(ii)] Assume first that $I$ is soft-edged, let $\xi\in \Sigma(I)$ and let $P= \sum_{k=1}^\infty \oplus P_k$ be a block stochastic matrix with the k-th direct sum block corresponds to the indices $n_{k-1} < i,j\le n_k$. 
By definition,  $\xi \le \alpha \eta$ for some $\eta\in \Sigma (I)$ and some $\alpha\in \co*$.  
By consolidating finite numbers of blocks and passing if necessary to a larger sequence $\alpha$ still in $\co*$, 
we can assume without loss of generality that $\alpha$ is constant on each block $(n_{k-1}, n_k ]$ and that $\sum_{k=1}^\infty \alpha_{n_k} < \infty$. 
By passing to a scalar multiple of $\eta$, we can furthermore assume that $\sum_{k=1}^\infty \alpha_{n_k} =1$. But then $P\xi \le P(\alpha \eta)= \big (\sum_{k=1}^\infty \oplus \alpha_{n_k}P_k\big)\eta.$
Since $R:=\sum_{k=1}^\infty \oplus \alpha_{n_k}P_k\in \mathscr B$, $(P\xi)^*\le(R\eta)^* \in \Sigma(I)$, which proves the claim

Assume now that $I$ is soft-complemented, let $\xi\in \Sigma(I)$, and let $P= \sum_{k=1}^\infty \oplus P_k$ be a block stochastic matrix  with the k-th direct sum block corresponds to the indices $n_{k-1} < i,j\le n_k$. 
To prove that $(P\xi)^*\in \Sigma (I)$, it is enough to prove that $\alpha(P\xi)^*\in \Sigma (I)$ for every $\alpha\in\co*$.   
The case where $\xi$ and hence $P\xi$ have finite support  being trivial,  assume that  $\xi_n> 0$ for all $n$. Then $((P\xi)^*)_n = (P\xi)_{\Pi (n)}$ for some permutation $\Pi$. 
Let $\gamma:= <\alpha_{\Pi^{-1}(n)}>$, then $\gamma \in \text{c}_o$.  Choose a sequence $\tilde{\gamma}\in\co*$ with $\gamma \le \tilde{\gamma}$ and  a subsequence $n_{k_i}$ for which $r:=\sum_{i=1}^\infty \tilde{\gamma} _{n_{k_i}}< \infty$. Define $\delta_j:=  \tilde{\gamma} _{n_{k_i}}$ for all $n_{k_i}< j\le n_{k_{i+1}}$ and $R:= \sum_{i=1}^\infty \oplus   \tilde{\gamma} _{n_{k_i}} \big(\sum_{h=k_i+1}^{k_{i+1}}\oplus P_h\big)$. Then $\delta\in\co*$, $\gamma \le \delta$, and $\frac{1}{r} R\in \mathscr B$. Since  $(\gamma P\xi)_{\Pi(n)}= (\alpha(P\xi)^*)_n$ is monotone nonincreasing, it follows that $\alpha(P\xi)^*= (\gamma P\xi)^*$ and hence
$
\alpha(P\xi)^* \le ( \delta P\xi )^*= (R\xi)^* = (\frac{1}{r} R(r\xi))^*\in \Sigma(I).
$
which concludes the proof.
 \ep

\bT{T:8.9}
 If $I$ is soft-edged or soft-complemented and it is invariant under $\mathscr C$, then $I$ is am-closed.
\eT
\bp
By combining parts (i) and (ii) of Lemma \ref {L:8.8},  we see that $I$ is invariant under block stochastic matrices and that $I \supset \mathscr L_1$. By  Corollary \ref {C:8.7} it follows that $I$ is am-closed.
\ep

Is is easy to see that Varga's ideal $ \cl(\xi) $ (see (\ref{e:42})) which is not  am-closed, but being 
 Banach is invariant under  $\mathscr C$, is neither soft-edged nor soft-complemented.


\begin{thebibliography}{99}

\bibitem{AGPS87}
Albeverio, S., Guido, D., Ponosov, A., and Scarlatti, S., 
\textit{Singular traces and compact operators,} J. Funct. Anal. \textbf{137} (2) (1996), 281--302. 

\bibitem{AMRS}
Antezana,  J.,  Massey, P.,  Ruiz, M, and Stojanoff, D., \textit{The Schur-Horn Theorem for operators and frames with prescribed norms and frame operator,} Illinois J of Math.  to appear.

\bibitem{AK02}
Arveson, W.  and Kadison, R. V.,  \textit{Diagonals of self-adjoint operators.}  Operator theory, operator algebras, and applications, 
Contemp. Math. \textbf{414}  Amer. Math. Soc., Providence, RI (2006), 247--263.

\bibitem{Bg46}
Birkhoff, G., \textit{Tres observaciones sobre el algebra lineal,} Univ. Nac. Tucuman Rev. Ser. A \textbf{5} (1946), 147--151.

\bibitem{jC41}
Calkin, J. W., \textit{Two-sided ideals and congruences in the ring
of bounded operators in Hilbert space,}
Ann. of Math.  \textbf{42} (2) (1941),  839--873.

\bibitem{CL02}
Casazza, P., Leon, M., \textit{Frames with a given frame operator} (2002) Preprint.

\bibitem{DFWW}
Dykema, K., Figiel, T., Weiss, G.,  and Wodzicki, M., \textit{The
commutator structure of operator ideals,}
Adv. Math.,  \textbf{185} (1) (2004), 1--79.

\bibitem{Fa49}
Fan, K.,  \textit{On a theorem of Weyl concerning eigenvalues of linear transformations,} Proc. Nat. Acad. Sci. U.S.A.  \textbf{35} (1949), 652--655.

\bibitem{Fa51}
Fan, K., \textit{Maximum properties and inequalities for the eigenvalues of completely continuous operators.}  Proc. Nat. Acad. Sci. U.S.A.  \textbf{37} (1951), 760--766.

\bibitem{GiMa64}
Gohberg, I.C., and Markus, A. S., \textit{Some relations between eigenvalues and matrix elements of linear operators} Mat. Sb. \textbf{64} (106) (1964), 481-496 (Russian); Amer. Math. Soc. Transl. (2) \textbf{52} (1966) 201-216  (English)

\bibitem{pH82}
Halmos, P. R., \textit{A Hilbert Space Problem Book,} 2nd Edition, Graduate Texts in Mathematics  \textbf(19), Springer-Verlag (1982). From 1st Edition, D. Van Nostrand Co., Inc., Princeton, N.J.-Toronto, Ont.-London, 1967.

\bibitem{HLP52}
Hardy, G. H., Littlewood, J. E., and P\'olya, G., \textit{Inequalities.} 2d ed. Cambridge University Press, 1952.

\bibitem{Horn}
Horn, R. A. \textit{Doubly stochastic matrices and the diagonal of a rotation matrix.} Amer. J. Math.  \textbf{76} (1954), 620--630.

\bibitem{HJ85}
Horn, R. A. and Johnson, C. R. \textit{Matrix Analysis}. Corrected reprint of the 1985 original. Cambridge University Press, Cambridge, 1996.

\bibitem{HJ91}
Horn, R. A.; Johnson, C. R. \textit{Topics in Matrix Analysis}. Corrected reprint of the 1991 original. Cambridge University Press, Cambridge, 1995.

\bibitem{Kr02a}
Kadison, R.,  \textit{The Pythagorean Theorem I: the finite case,}  Proc. Natl. Acad. Sci. USA \textbf{99} (7) (2002), 4178--4184.

\bibitem{Kr02b}
Kadison, R.,  \textit{The Pythagorean Theorem II: the infinite discrete case,}  Proc. Natl. Acad. Sci. USA \textbf{99}  8  (2002), 5217--5222.

\bibitem{vKgW02}
Kaftal, V. and Weiss, G., \textit{Traces, ideals, and arithmetic
means,} Proc. Natl. Acad. Sci. USA \textbf{99} (11) (2002),
7356--7360.

\bibitem{vKgW04-Soft}
Kaftal, V. and Weiss, G., \textit{Soft ideals and arithmetic mean
ideals,}  Int. Eq. Oper. Theory. \textbf{58} (2007) 363-405 

\bibitem{vKgW04-2nd Order}
Kaftal, V. and Weiss, G., \textit{Second order arithmetic means
 in operator ideals,} Operators and Matrices textbf{1} 2, (2007) 235-256
  
\bibitem{vKgW04-Traces}
Kaftal, V. and Weiss, G., \textit{Traces on operator ideals and
arithmetic means,} to appear in J. Oper. Theory.
  
\bibitem{vKgW07-OT21}
Kaftal, V. and Weiss, G., \textit{A survey on the interplay between arithmetic mean ideals, traces, lattices of operator ideals, and an infinite Schur-Horn majorization theorem}, to appear in OT21, Recent advances in operator theory, operator algebras, and their applications, Oper. Theory Adv. Appl., Birkhauser, Basel.

\bibitem{vKgW04-Density}
Kaftal, V. and Weiss, G., \textit{The $B(H)$ lattices, density and 
arithmetic mean ideals,} preprint.

\bibitem{nK89}
Kalton, N. J., \textit{Trace-class operators and commutators,} J.
Funct. Anal. \textbf{86} (1989), 41--74.

\bibitem{KkLd04}
Kornelson, K. and Larson, D., \textit{Rank-one decompositions of operators and construction of frames,}
Wavelets, frames and operator theory, Contemp. Math \textbf{345}, Amer. Math. Soc. Providence, Ri (2004),  203--214.  
  
\bibitem{Kd60}
Kendall, D., \textit{ On infinite doubly-stochastic matrices and Birkhoff's problem 111}, J. London Math. Soc. \textbf {35} (1960), 81--84.

\bibitem{Lm1905}
Lorenz, M. O.,  \textit{Methods of measuring concentration of wealth,} J. Amer. Statist. Assoc. \textbf {9} (1905), 209--219. 

\bibitem{aM64}
Markus, A. S., \textit{The eigen-and singular values of the sum and product of linear operators,} Uspekhi Mat. Nauk 4  \textbf(118) (1964),  93--123.

\bibitem{MO79}
Marshall, A. W. and Olkin, I., \textit{Inequalities: Theory of 
Majorization and its Applications,} Academic Press Inc. [Harcourt 
Brace Jovanovich Publishers], Mathematics in Science and Engineering.
\textbf{143} (1979).

\bibitem{Mi58}
Mirsky, L., \textit{Matrices with prescribed characteristic roots and diagonal elements.} J. London Math. Soc. \textbf{33} (1958), 14--24.

\bibitem{Mi60}
Mirsky, L., \textit{Symmetric gauge functions and unitarily invariant norms.} Quart. J. Math. Oxford Ser. [2] 
\textbf{11}  (1960), 50--59.

\bibitem{Mu1903}
Muirhead, R.F., \textit{Some methods applicable to identities and inequalities of symmetric algebraic functions of $n$ letters.} Proc. Edinburgh Math. Soc. \textbf{21} {1903}, 144--157.

\bibitem{Na99}
Neumann, A.,  \textit{An infinite-dimensional generalization of the Schur-Horn convexity theorem.} J. Funct. Anal., \textbf{161} (2) {1999},  418--451.

\bibitem{Oa52}
Ostrowski, A. M. , \textit{Sur quelques applications des fonctions convexes et concaves au sens the I. Schur} J. Math. Pures Appl. [9], \textbf{31} (1952), 253--292.

\bibitem{Si23}
Schur, I., \textit{\"{U}ber eine Klasse von Mittelbildungen mit Anwendungen auf der Determinantentheorie,} Sitzungsber. Berliner Mat. Ges., ( \textbf{22})  (1923),  9--29.

\bibitem{Vj89}
Varga, J., \textit{Traces on irregular ideals.} Proc. Amer. Math.
Soc. \textbf{107} (1989), 715--723.


\bibitem{mW02}
Wodzicki, M., \textit{Vestigia investiganda,} Mosc. Math. J.
\textbf{4} (2002), 769--798, 806. 

\end{thebibliography}
\end{document}